\documentclass{amsart}

\usepackage{amsmath, amssymb, mathtools}
\usepackage{amscd}
\usepackage[body={16cm,22cm}]{geometry}
\usepackage{adjustbox}
\usepackage[colorlinks=true,linkcolor=blue, citecolor=cyan, plainpages,urlcolor=magenta]{hyperref}
 
\usepackage[all]{xy}
\usepackage{cleveref}
\usepackage{tikz-cd}

\newcommand{\BA}{{\mathbb {A}}}

\newcommand{\BC}{{\mathbb {C}}}

\newcommand{\BG}{{\mathbb {G}}}

\newcommand{\BM}{{\mathrm{BM}}}
\newcommand{\BN}{{\mathbb {N}}}

\newcommand{\BP}{{\mathbb {P}}}
\newcommand{\BQ}{{\mathbb {Q}}}

\newcommand{\BZ}{{\mathbb {Z}}}

\newcommand{\CF}{{\mathcal {F}}}

\newcommand{\CI}{{\mathcal {I}}}

\newcommand{\CL}{{\mathcal {L}}}

\newcommand{\CO}{{\mathcal {O}}}

\newcommand{\CT}{{\mathcal {T}}}

\newcommand{\RA}{{\mathbf {A}}}

\newcommand{\Br}{{\mathrm{Br}}}

\newcommand{\Pic}{\mathrm{Pic}}

\newcommand{\Spec}{{\mathrm{Spec}}}

 \def\dim{{\mathrm{dim}}}
 \def\top{{\mathrm{top}}}
  
\newcommand{\codim}{{\mathrm{codim}}}
\newcommand{\chara}{{\mathrm{char}}}
\def\Hom{{\mathrm{Hom}}}
\def\Bl{{\mathrm{Bl}}}
\newcommand{\Ker}{{\mathrm{Ker}}}
\newcommand{\CH}{{\mathrm{CH}}}
\newcommand{\Coker}{{\mathrm{Coker}}}
\renewcommand{\Im}{{\mathrm{Im}}}
\def\pr{{\mathrm{pr}}}
\def\tr{{\mathrm{tr}}}
\def\rank{{\mathrm{rank}}}

\newcommand{\et}{ {\rm{\acute et}}}
\newcommand{\iso}{\stackrel{\sim}{\rightarrow} }

\newcommand{\bk}{\bar{k}}

\newcommand{\tSym}{{\rm {\widetilde{Sym}}}}
\newcommand{\Sym}{{\mathrm{Sym}}}
\newcommand{\Gal}{{\mathrm{Gal}}}

\numberwithin{equation}{section}

\theoremstyle{remark}

\newtheorem{defi}{\rm{\textbf{Definition}}}[section]

\newtheorem{rem}[defi]{\rm{\textbf{Remark}}}

\newtheorem{ques}[defi]{\rm{\textbf{Question}}}

\theoremstyle{plain}

\newtheorem*{claim}{\rm{\textbf{Claim}}}
\newtheorem{thm}[defi]{\rm{\textbf{Theorem}}}

\newtheorem{cor}[defi]{\rm{\textbf{\textbf{Corollary}}}}

\newtheorem{lem}[defi]{\rm{\textbf{Lemma}}}

\newtheorem{prop}[defi]{\rm{\textbf{\textbf{Proposition}}}}

\newtheorem {conj}[defi]{\rm{\textbf{Conjecture}}}

\title[Brauer--Manin obstruction for symmetric powers]
{Brauer--Manin obstruction for symmetric powers}

\author[Cao]{Yang Cao}
\address{
School of Mathematics\\
Shandong University\\
250100, Jinan, China}
\email{yangcao1988@gmail.com; yangcao1988@email.sdu.edu.cn}

\author[Zhou]{Lin Zhou}
\address{
Institute of Algebraic Geometry\\
Leibniz University Hannover\\
30167, Hannover, Germany}
\email{zhou@math.uni-hannover.de}
\date{\today}

\begin{document}

\begin{abstract}
For a quasi-projective smooth geometrically integral variety $X$ over a number field, 
we prove that the Brauer--Manin obstruction for its symmetric powers coincides with the obstruction coming from $X$.
As a consequence, Brauer--Manin obstruction for the symmetric powers cannot provide new cohomological obstruction
to weak approximation of zero-cycles for $X$. 
The key intermediate ingredients are a Dold--Thom type result for Brauer groups and a precise descent formula for adelic points.
\end{abstract}

\maketitle

\setcounter{tocdepth}{1}
\tableofcontents

\section{Introduction}

Let $k$ be a number field and $X$ be a smooth geometrically integral over $k$.
We denote by $X(k) $ the set of rational points and $X(\RA_k)$ the topological space of adelic points of $X$; see \cite{Cod12} for details.
It is an important question (e.g., Hasse principle, weak approximation, strong approximation) to study the closure $\overline{X(k)}$ of the diagonal embedding
$$X(k) \subset X(\RA_k) .$$
In his foundational paper \cite{Man71}, Manin combined class field theory and Grothendieck’s theory of Brauer group 
to define a natural closed subset $X(\RA_k)^{\Br}$ of $X(\RA_k)$ which contains $X(k)$.
We say a variety $X$ satisfies \emph{weak approximation with Brauer--Manin obstruction} if $\overline{X(k)}=X(\RA_k)^{\Br}$.
In \cite[p. 174]{CT03}, Colliot-Th\'el\`ene formulated the following programmatic conjecture, which was first proposed by Colliot-Th\'el\`ene and Sansuc in 1979 in the case of surfaces.

\begin{conj}[Colliot-Th\'el\`ene]
Any proper smooth geometrically integral rationally connected variety over a number field satisfies weak approximation with Brauer--Manin obstruction.
\end{conj}

This conjecture has been studied extensively in the last four decades, but it is far from being resolved,
see \cite{CTS}, \cite{Sk}, \cite{W18}, \cite{W23} for a survey. 
On the other hand, for varieties that are not rationally connected, many kinds of cohomological obstructions have been defined,
including the descent obstruction, the \'etale Brauer--Manin obstruction and their iterated obstruction, 
see \cite{Sk99}, \cite[Def. 4.3]{HS02}, \cite[\S 8]{Po} for their definitions
and \cite{H02}, \cite{D09}, \cite{Sk09}, \cite{CDX} and \cite{C4} for their comparisons.

In parallel development, we study the diagonal embedding of the Chow group $\CH_0(X)$ of zero-cycles up to rational equivalence
$$ \CH_0(X) \to \prod_v \CH_0(X_{k_v}) . $$
We denote by $\CH_{0,\RA}(X)$ the canonical quotient of $\prod_v \CH_0(X_{k_v})$ defined in \cite[(1.1)]{HW16}.
As in \cite[\S 15.1]{CTS}, we have the Brauer--Manin pairing
$\langle-,-\rangle_{\BM} \colon \CH_{0,\RA}(X) \times \Br(X) \to \BQ/\BZ  $
and, similar to the rational points case, this Brauer--Manin pairing  and class field theory induce the following complex

\begin{equation}\label{CTKSconj-eq}
   \widehat{\CH_0(X) }\to \widehat{\CH_{0,\RA}(X)} \to \Hom(\Br(X),\BQ/\BZ),
\end{equation}
where we denote $\hat{M}:=\varprojlim_nM/nM$ for any abelian group $M$.
Colliot-Th\'el\`ene, Kato and Saito (in \cite{KS86} and \cite{CT95}) then formulated the following conjecture
about the \emph{local-global principle and weak approximation for zero-cycles}.

\begin{conj}[Colliot-Th\'el\`ene, Kato, Saito]\label{Conj-zero}
For any proper smooth geometrically integral variety $X$ over a number field, the sequence (\ref{CTKSconj-eq}) is exact.
\end{conj}

Unlike the case of rational points, Conjecture \ref{Conj-zero} does not require $X$ to be rationally connected,
i.e., it is expected that the Brauer--Manin obstruction is the only one for zero-cycles.
On the other hand, 
the notion of descent obstruction and of \'etale Brauer--Manin obstruction cannot be generalized to zero-cycles directly.
Some attempts at these generalizations have been formulated  in \cite{Ba24} and \cite{Ba26},
but none of them has been proven to be essentially stronger than the Brauer--Manin obstruction.

\medskip

For a projective smooth geometrically integral variety $X$,
we denote by $\Sym^nX$ the symmetric powers of $X$ and by $\tSym^nX$ a proper resolution of singularity of $\Sym^nX$.

In the present article, we study an object between rational points $X(k)$ and Chow group $\CH_0(X) $: 
the rational point of symmetric powers $\tSym^nX(k) $, 
because symmetric powers $\Sym^nX$ parametrize effective zero-cycles.
More specifically, for any field extension $L/k$, we have the canonical map
$$[-]_L\colon \ \tSym^nX(L) \to \Sym^nX(L) \to \CH_0(X_L) .$$
This induces the canonical map 
\begin{equation}\label{def-eq-thm-1-cor}
[-] \colon \tSym^nX(\RA_k) \to \Sym^nX(\RA_k) \to \widehat{\CH_{0,\RA}(X)},
\end{equation}
which sends $\tSym^nX(k) $ to the image of $\CH_0(X)$.
By \cite[Lem. 1.8]{W12}, the map $[-]$ is continuous.
This induces a natural question:

\begin{ques}\label{quetion}
Can we use cohomological obstructions for rational points on $\tSym^nX $ to define new cohomological obstruction for zero-cycles on $X$?
\end{ques}

If $\Pic(X_{\bk})$ is torsion free, by \cite[Thm. A and B]{LLZ26}, we have an isomorphism $\frac{\Br(X)}{\Br(k)}\cong \frac{\Br(\tSym^nX)}{\Br(k) }$ 
and the Brauer--Manin obstruction of $\tSym^nX$ cannot provide new cohomological obstruction for zero-cycles on $X$.
Without this assumption, we cannot expect the isomorphism of Brauer groups even over $\bk$; see Remark \ref{rem-S4-last}. Therefore, for general smooth geometrically integral $X$, one cannot answer Question \ref{quetion}, even for Brauer--Manin obstruction, directly.

 To address this problem, we define the canonical transfer $ \tr_{X,\Br}\colon  \Br(X)\to \Br(\tSym^nX)$ as in (\ref{eq-def-transfer}).
By Lemma \ref{thm-1-cor-lem}, this homomorphism is compatible with our map $[-]$, 
i.e., for any $a\in \Br(X)$ and any $(x_v)\in \tSym^nX(\RA_k)$, we have
\begin{equation}
 \langle [(x_v)],a\rangle_{\BM} = \langle (x_v), \tr_{X,\Br}(a)\rangle_{\BM}  . 
 \end{equation}

The key result of present article is to compare the Brauer--Manin obstruction on $\tSym^nX$ 
given by the whole Brauer group and  by the subgroup coming from $X$.

\begin{thm}\label{thm-1}
Let $X$ be a quasi-projective smooth geometrically integral $k$-variety.
Then for any positive integer $n$, we have
$$\tSym^nX(\RA_k)^{\Im(\tr_{X,\Br})} =\tSym^nX(\RA_k)^{\Br}.$$
\end{thm}

As a consequence, we prove that the Brauer--Manin obstruction on $\tSym^nX$ cannot provide new cohomological obstruction for zero-cycles on $X$.
More precisely, we have the following result.

\begin{cor}\label{thm-1-cor}
Let $X$ be a projective smooth geometrically integral $k$-variety.
Let $n$ be a positive integer and $[-] $ be the map in  (\ref{def-eq-thm-1-cor}).
Then, for all $(x_v)\in \tSym^nX(\RA_k)$, we have $(x_v)\in \tSym^nX(\RA_k)^{\Br}$ if and only if
$$[(x_v)]\in \Ker \left(\widehat{\CH_{0,\RA}(X)} \to \Hom(\Br(X),\BQ/\BZ )\right) .$$
\end{cor}

\medskip

The proof of Theorem \ref{thm-1} contains two key intermediate ingredients: 
a Dold--Thom type result for Brauer groups and a precise descent formula for adelic points; see Theorem \ref{thm-2} and Proposition \ref{desc-result-1-prop}, respectively.
Both of them are related to the following algebraic version of Hurewicz map (\ref{def-epsilon-eq}).

More concretely, for any connected scheme $X$ and any finite Galois cover $X'\to X$, 
the Hochschild--Serre spectral sequence (\cite[Thm. III.2.20]{Mi80}) induces a homomorphism
$H^2(\Gal(X'/X),\mu_{\infty})\to H^2(X,\mu_{\infty}) $.
Passing to the limit, we obtain a natural homomorphism
\begin{equation}\label{def-epsilon-eq}
\epsilon_X\colon  H^2_{cont}(\pi_1(X),\mu_{\infty})\to H^2(X,\mu_{\infty}).
\end{equation}
This is an algebraic version of the \textbf{Hurewicz map}.
We denote
\begin{equation}\label{def-epsilon-eq1}
\epsilon'_X\colon H^2_{cont}(\pi_1(X),\mu_{\infty})\xrightarrow{\epsilon_X} H^2(X,\mu_{\infty}) \to \Br(X),
\end{equation}
where the second homomorphism is induced by the Kummer sequence.

The following result is an algebraic analogue of Dold--Thom theorem for Brauer group.  Moreover, the Dold--Thom theorem (Theorem \ref{thm- toplogical DT}) is an important result in algebraic topology.

\begin{thm}[cf. {Corollary \ref{cor-DT thm for deg2}}]\label{thm-2}
Let $k$ be an algebraically closed field of characteristic $0$ and $X$ be a quasi-projective smooth integral $k$-variety.
Then for any proper resolution $\tSym^nX\to \Sym^nX$, we have
\begin{enumerate}
\item 
the \'etale fundamental groups $\pi_1(\tSym^nX) \cong \pi_1(X)^{ab}$ is commutative for any $n\geq 2$;
\item 
the natural homomorphism
$$( \epsilon'_{\tSym^nX} , \tr_{X,\Br})\colon \ H^2_{cont}(\pi_1(\tSym^nX),\mu_{\infty}) \oplus \Br(X)  \twoheadrightarrow \Br(\tSym^nX) $$
is surjective for any positive integer $n$.
\end{enumerate}
\end{thm}

Now, we introduce the other ingredient---the descent formula.
The descent method is a widely used method to study Hasse principle and the weak approximation for rational points.
It was firstly established by Colliot-Th\'el\`ene and Sansuc in \cite{CTS87}.

In the last twenty years, many precise descent formulas for the adelic points orthogonal to some certain Brauer groups were established
by Demarche, Harpaz, Wittenberg, Xu and the first author 
(see \cite[Thm. 1.2]{CDX}, \cite[Prop. 5.2 and Thm. 5.9]{C1},  \cite[Thm. 2.1]{HW20}, \cite[Prop. 5.1]{C4} and \cite[Prop.~6.2]{C5}).

Here we establish a descent formula with respect to the Hurewicz map in the case the geometric fundamental group is abelian. 
This formula shows that the algebraic Brauer group provides more cohomological obstruction than the image of the Hurewicz map.

\begin{prop}\label{desc-result-1-prop}
Let $X$ be a quasi-projective smooth geometrically integral $k$-variety and $\Br_1(X)$ be its algebraic Brauer group as in (\ref{def-eq-H^i_1}).
Assume $\pi_1(X_{\bk})$ is abelian. 
Then 
\begin{enumerate}
\item we have $ X(\RA_k)^{\Br_1(X)}\subset X(\RA_k)^{\Im(\epsilon'_X)}$;
\item for any finite commutative $k$-group $A$ and any $A$-torsor $f\colon Y\to X$ with $Y$ geometrically integral,
we have 
$$ X(\RA_k)^{\Br_1(X)} \subset 
\bigcup_{\sigma\in H^1(k,A)}f_{\sigma} (Y_{\sigma}(\RA_k)^{ \Im(\epsilon'_{Y_{\sigma}})}).$$
\end{enumerate}
\end{prop}

The paper is organized as follows. In Section \ref{section2}, we recall some basic results on symmetric powers of varieties, review the Hurewicz map, the K\"unneth formula, and the Dold--Thom theorem in the classical topological setting, and compare them with their \'etale counterparts. 
In Section \ref{section3}, we study the Brauer group of symmetric powers of curves by using relative Picard scheme.
Section \ref{section4} contains the proof of Theorem~\ref{thm-2}. 
In Section \ref{section5}, we establish the descent formula (see Proposition \ref{desc-result-1-prop}) and use it to prove Theorem \ref{thm-1}. Finally, Section \ref{section6} is devoted to the proof of Corollary \ref{thm-1-cor}.

\section{Preliminaries}\label{section2}

\subsection{Notations and conventions}
Throughout the paper, without emphasis, we denote by $k$ a field of characteristic zero and we fix an algebraic closure $\bk$ of $k$. 
We denote by $\Gamma_k\coloneqq \mathrm{Gal}(\bk/k)$ the absolute Galois group of $k$.

Let $k_{\et}$ be the small \'etale site of $k$ and $F$ be a torsion sheaf over $k_{\et}$. Then $F(\bk)$ is a torsion $\Gamma_k$-module. 
Write $\mu_m\coloneqq \{a\in \bk^{\times}\ |\ a^m=1\}$ and $\mu_{\infty}\coloneqq \cup_m\mu_m\subset \bk^{\times}$.

A $k$-variety $X$ is a separated $k$-scheme of finite type. For every $k$-variety $X$, we write $X_{\bar{k}}\coloneqq X\times_{\bar{k}}k$. We denote by $H^i(X,F)$ the \'etale cohomology group of $X$ with coefficient $F$ and
$\Br(X)\coloneqq H^2(X,\BG_m)$ the cohomological Brauer group,
and we denote the algebraic part of \'etale cohomology by
\begin{equation}\label{def-eq-H^i_1}
H^i_1(X,F)\coloneqq \Ker(H^i(X,F) \to H^i(X_{\bk},F) )\ \ \ \text{and}\ \ \ \Br_1(X)\coloneqq \Ker(\Br(X)\to \Br(X_{\bk})).
\end{equation}
For any connected $k$-variety $X$, let $\pi_1(X)$ be the \'etale fundamental group of $X$. 
Here, we ignore the choice of the base point without emphasis. 
Let $\pi_1(X)^{ab}$ be the maximal abelian quotient of  $\pi_1(X)$. For any  torsion sheaf $F$ on $k_{\et}$ and any connected $k$-variety $X$, we have a canonical isomorphism (\cite[\S XI.5]{SGA1})
\begin{equation}\label{eq-fund-gps-3}
H^1(\pi_1(X),F(\bk))\iso H^1(X,F)\ \ \ \text{and hence}\ \ \ H^1(X_{\bk},F)\cong \mathrm{Hom}_{cont}(\pi_1(X_{\bk}),F(\bk)),
\end{equation}
where the action of $\pi_1(X)$ on $F(\bk)$ is induced by the action of $\Gamma_k$.

\medskip

We recall the notion of presheaves with transfers (see \cite{Voe06} for more details).

Let $X,Y$ be two smooth varieties over $k$.
If $Y$ is connected, an {\bf elementary correspondence} from $Y$ to $X$ is an integral closed subscheme $W\subset Y\times X$ which is finite and surjective over $Y$,
and the {\bf finite correspondences} $Cor(Y,X)$ is the free abelian group generated by elementary correspondences (\cite[Def. 1.1]{Voe06}).
In general, let $Cor(Y,X)\coloneqq \oplus_i Cor(Y_i,X)$ if $Y=\sqcup_i Y_i$ with $Y_i$ connected.
For example, any morphism $f\colon  Y\to X$ induces its graph $\Gamma_f\subset Y\times X$, which is a finite correspondence.

Let  $Cor_k$ be the category whose objects are smooth varieties over $k$ and morphisms are $Cor(Y,X)$.
A presheaf with transfers is a contravariant additive functor from $Cor_k$ to  the category of abelian groups.

Let $F$ be a torsion sheaf over $k_{\et}$.
Then the \'etale cohomology $H^i(-,F)$ is a presheaf with transfers; see \cite[Lem. 6.11, Lem. 6.21]{Voe06}.
Moreover, for any morphism $f\colon Y\to X$, 
the homomorphisms $ \Gamma_f^*\colon H^i(X,F)\to H^i(Y,F)$
induced by the finite correspondence $\Gamma_f$ are exactly the pullback $f^*$ of the cohomology groups.

Let $G$ be a linear algebraic group and $f\colon Y\to X$ be a (left) $G$-torsor.
For any $\sigma\in H^1(k,G)$, we denote by $G_{\sigma}$ (resp. $f_{\sigma}\colon Y_{\sigma}\to X $) the twist of $G$ (resp. $f$) by $\sigma$.
Therefore, $f_{\sigma} \colon Y_{\sigma}\to X$ is a (left) $G_{\sigma}$-torsor.
If $G$ is commutative, then $G_{\sigma}=G$ for any $\sigma\in H^1(k,G)$.

\subsection{Basic properties for symmetric powers}
In this part, we recall some basic properties for symmetric powers, which are well known to the experts.

At the beginning, we construct a canonical resolution of an open subset of the symmetric powers,
and we will prove this resolution has very nice cohomological properties in the following sections.

For a quasi-projective smooth geometrically integral $k$-variety $X$, we define:
\begin{enumerate}
    \item $\Sym^nX\coloneqq X^n/S_n$ to be the $n$-th symmetric power of $X$, where $S_n$ is the symmetric group of permutations of $n$ elements and $X^n\coloneqq X\times\cdots\times  X$ ($n$-times);
    \item $U^n(X)\hookrightarrow X^n$ to be the open subscheme of pairwise distinct points, 
    and $V^n(X)\hookrightarrow X^n$ to be the open subscheme of at most two points coincide;
    \item $U(X)$ (resp. $V(X)$) to be the image of $U^n(X)$ (resp. $V^n(X)$) under the natural projection $X^n\to \Sym^nX$;
    \item $\Bl_{V\setminus U}(V^n(X)) $ to be the blowing-up of $V^n(X)$ along $V^n(X)\setminus U^n(X)$, 
    and $\tilde{V}(X)$ to be the quotient $ \Bl_{V\setminus U}(V^n(X))/S_n$, 
    where the action of $S_n$ on $V^n(X)$ and $U^n(X)$ lifts to an action on the blowing-up;
    \item $h\colon \tSym^nX\to \Sym^nX$ to be a proper resolution of singularities that is an isomorphism over $U(X)$.
\end{enumerate}
Then $U(X)$ is an open subscheme of all $\tilde{V}(X)$, $\tSym^nX$ and $\Sym^nX$. 
By Proposition \ref{notation-sym-prop} below, $\tilde{V}(X)$ is already smooth. 
 We may assume  $\tilde{V}(X)\subset \tSym^nX$.
Therefore, we have open immersions
\begin{equation}\label{notation-sym-eq}
U^n(X)\subset V^n(X)\subset X^n,\ \ \ U(X)\subset V(X)\subset \Sym^nX \ \ \ 
\text{and}\ \ \ U(X)\subset \tilde{V}(X)\subset \tSym^nX.
\end{equation}

Moreover, if $X$ is a curve, $V^n(X)\setminus U^n(X)$ is a smooth divisor of $U^n(X)$. 
Then the blowing-up is trivial and $\tilde{V}(X)=V(X)$.

\begin{prop}\label{notation-sym-prop}
Under the above hypothesis, $\tilde{V}(X)$ is smooth 
and $\tilde{V}(X)\setminus U(X)$ is a geometrically integral divisor of $\tilde{V}(X)$.
\end{prop}

\begin{proof}
For any $1\leq j<l\leq n$, let $D_{j<l}\coloneqq\{(x_i)\in V^n(X)|\  x_j=x_l\} $.
Denote $\Delta$ the image of the diagonal $X\to X^2$.
Then $D_{1<2}$ is an open subset of $\Delta\times X^{n-2}$, hence it is smooth and geometrically integral.
The same statement holds for all $D_{j<l}$.
Since $V^n(X)\setminus U^n(X)=\coprod\limits_{1\leq j<l\leq n} D_{j<l} $,
we have 
\begin{equation*}
    (\Bl_{V\setminus U}V^n(X) ) \setminus U^n(X) =\coprod_{1\leq j<l\leq n} E_{j<l},
\end{equation*}
where $E_{j<l}$ is the inverse image of $D_{j<l}$.
By the property of blowing-up along a smooth center, $E_{j<l}$ is a divisor of $\Bl_{V\setminus U}V^n(X) $ 
and the projection $E_{j<l}\to D_{j<l}$ is a projective bundle. 
Hence $E_{j<l}$ is geometrically integral.
Since the $S_n$ transfers all $E_{j<l}$, 
the quotient $(\coprod\limits_{j<l}E_{j<l})/S_n = \tilde{V}(X)\setminus U(X)$ is geometrically integral.

Now we prove that $\tilde{V}(X)$ is smooth. 
The case $n=1$ is trivial.
If $n=2$, this follows from \cite[Example 7.3.1 (3)]{FGAexplained}.
For $n\geq3$, we need to show that $\tilde{V}(X)$ is smooth at the image of all $P\in \sqcup_{j<l}E_{j<l}$.
We may assume $P\in E_{1<2}$.
Denote $U_{1<2}\coloneqq U^n(X)\cup E_{1<2}$ and $Q\in D_{1<2}$ the image of $P$. 
We recall the notation from \cite[\S V.2]{SGA1}: the inertia subgroup $I_P$ of $P$ is the $\sigma\in S_n$
such that $\sigma(P)=P$ and the induced $\sigma^*\colon k(P)\to k(P)$ is identity.
Since $Q\in D_{1<2}$, we have 
$$I_P\subset I_Q\subset S_2=\{{\rm id},(1,2)\}\subset S_n.$$
By \cite[Prop. V.2.2]{SGA1}, $U_{1<2}/S_2 \to  \tilde{V}(X)$ is \'etale at $P$.
Hence we need only to show $U_{1<2}/S_2$ is smooth.
Since $U_{1<2} $ is an open subset of $ Bl_{V\setminus U}(X^2)\times X^{n-2}$,
the quotient $U_{1<2}/S_2 $ is an open subset of $ \Bl_{V\setminus U}(X^2)/S_2\times X^{n-2} $.
Since $ \Bl_{V\setminus U}(X^2)/S_2$ is smooth (the $n=2$ case),  $U_{1<2}/S_2 $ is smooth and the result follows.
\end{proof}

Now we define the canonical transfer homomorphism for the Brauer group of symmetric powers.

Since the universal family $W_0\subset X\times \Sym^nX$ is the subvariety induced by the following closed embedding
\begin{equation}\label{equ-universal family}
\begin{split}
      X\times \Sym^{n-1}X&\hookrightarrow X\times \Sym^nX\\
      (x,[x_1,\cdots,x_{n-1}])&\mapsto (x,[x,x_1,\cdots,x_{n-1}]),
\end{split}
\end{equation}
 the natural morphism $W_0\to \Sym^nX$ is finite of degree $n$. Note that for every resolution of singularity $h\colon \tSym^nX\to \Sym^nX$, the subscheme 
 \begin{equation}\label{equ-def-W}
    W\coloneqq (\mathrm{id}_X\times h)^{-1}(W_0)\in Cor(\tSym^nX,X) 
 \end{equation}
  corresponds to a finite correspondence $W$; then for any torsion sheaf $F$ over $k_{\et}$ and for every $i\in \mathbb{N}$, we define
\begin{equation}\label{eq-def-transfer}
\tr^i_X\coloneqq  W^* \colon H^i(X,F)\to H^i(\tSym^nX,F)\ \ \ \text{and}\ \ \ 
\tr_{X,\Br}\coloneqq  W^* \colon  \Br(X)\to \Br(\tSym^nX) .
\end{equation}
We call them the \textbf{canonical transfers}. 

\begin{lem}\label{lem-pullback-of-corr}
   Let $X$, $Y$ and $Z$ be smooth $k$-varieties with a flat $k$-morphism $f\colon X\to Y$. Then for every correspondence $W\in Cor(Y,Z)$, $\Gamma_f^*\circ W^*=(W')^*\in Cor(X,Z)$, where $W'$ is the pull-back of $W$ along the flat morphism $f\times {\rm id}:X\times Z\to Y\times Z$.
\end{lem}

\begin{proof}
   See \cite[Prop. 16.1.1]{fulton}. Note that by our assumptions, everything in the proof can be boiled down to the level of cycles, without passing to equivalence class of cycles. 
\end{proof}

\begin{prop}\label{eq-def-transfer1}
Let $k$ be a field of $\chara(k)=0$.
Let $X$ be a quasi-projective smooth $k$-variety.
Then the two compositions 
\begin{equation}\label{equ-def-transfer2}
     H^i(X,F) \xrightarrow{tr^i_X} H^i(\tSym^nX,F) \to H^i(U(X),F)\to H^i(U^n(X),F)
\end{equation}
 and
\begin{equation}\label{equ-def-transfer3}
    H^i(X,F) \xrightarrow{\sum_i\pr_i^*} H^i(X^n,F) \to H^i(U^n(X),F)
\end{equation}
are the same, where $\pr_i:X^n\to X$ is the $i$-th projection.
\end{prop}

\begin{proof}
By Lemma \ref{lem-pullback-of-corr}, the morphism (\ref{equ-def-transfer2}) is induced by the correspondence $W^n(X)$ satisfying the following commutative diagram
\begin{equation*}
    \xymatrix{W^n(X)\ar[r]\ar[d]&W|_{U(X)}\ar[r]\ar[d]&W\ar[d]\\
    U^n(X)\ar[r]^{\pi|_{U^n(X)}}&U(X)\ar@{^(->}[r]^j&\tSym^nX,
    }
\end{equation*}
where both squares are Cartesian, all horizontal arrows are flat and all vertical arrows are finite and surjective. Moreover, since $W$ is induced by the resolution of singularity $h\colon \tSym^nX\to \Sym^nX$, by universal property, the correspondence $W^n(X)$ is induced by the morphism
\begin{equation*}
    U^n(X)\xrightarrow{\pi|_{U^n(X)}}U(X)\xrightarrow{j}\tSym^nX\xrightarrow{h}\Sym^nX.
\end{equation*}
On the other hand, the above composition coincides with the following composition of morphisms
\begin{equation*}
    U^n(X)\xrightarrow{j'}X^n\xrightarrow{\pi}\Sym^nX,
\end{equation*}
where $j'$ is an open immersion. Hence, $W^n(X)$ is exactly the correspondence $(j'\times {\rm id})^*(W'')=W''|_{U^n(X)}$ by Lemma \ref{lem-pullback-of-corr}, where $W''\in Cor(X^n,X)$ is the correspondence  associated to the morphism $\pi$. Thus, 
\begin{equation*}
  W^n(X)^*=((j'\times {\rm id})^*W'')^*=j'^*\circ W''^*\colon H^i(X,F)\xrightarrow{W''^*} H^i(X^n,F)\xrightarrow{j'^*} H^i(U^n(X),F). 
\end{equation*}

Therefore, to finish this proof, it suffices to show that the morphism $W''^*\colon H^i(X,F)\to H^i(X^n,F)$ is exactly the morphism $ H^i(X,F) \xrightarrow{\sum_i\pr_i^*} H^i(X^n,F)$, i.e., $W''=\sum_i\Gamma_{\pr_i}\in Cor(X^n,X)$. In fact, note that the universal family over $\Sym^nX$ is the subvariety induced by (\ref{equ-universal family})
and we have the following Cartesian diagram
\begin{equation*}
    \xymatrix{\cup_i \Gamma_{\pr_i}\ar@{^(->}[rr]\ar[d]&&X\times X^n\ar[d]^{{\rm id}\times \pi}\\
    X\times \Sym^{n-1}X\ar@{^(->}[rr]&&X\times \Sym^nX,
    }
\end{equation*}
where the left morphism restricts to $\Gamma_{\pr_i}$ is the morphism 
$$(x,x_1,\cdots,x_i,\cdots,x_n)\mapsto (x,[x_1,\cdots,\widehat{x_i},\cdots,x_n])$$ 
with $x=x_i$. This implies that we have $W''=\sum_i\Gamma_{\pr_i}\in Cor(X^n,X)$ as cycles.
\end{proof}

\subsection{Topological Hurewicz map and K\"unneth formula}

In this section, we recall two fundamental constructions in algebraic topology: the classical Künneth formula and the classical Hurewicz map. 
We mainly focus on the homology and cohomology of degree $n=1,2$. 
Moreover, we discuss their connections with group homology and cohomology.

\begin{prop}\label{top-ex-Kunn}
Let $A$ and $B$ be two connected pointed CW complexes, and let $p_1$ $($\emph{resp.} $p_2)$ be the projection of  $A\times B$ to $A$ $($\emph{resp.} $B)$.  For any abelian group $F$, we have the following natural isomorphism
\begin{equation}\label{equ-top Kunn for H^1}
    (p_1^*,p_2^*)\colon \ H^1_{\top}(A,F)\oplus H^1_{\top}(B,F) \iso H^1_{\top}(A\times B,F) 
\end{equation}
and the following natural exact sequence
\begin{equation}\label{equ-top kun for H^2}
    0\to H^2_{\top}(A,F)\oplus H^2_{\top}(B,F)\xrightarrow{(p_1^*,p_2^*)} H^2_{\top}(A\times B,F)\xrightarrow{\lambda^{\top}_{A,B}}\Hom_\BZ(\pi_1^{\top}(A)^{ab}\otimes_{\BZ} \pi_1^{\top}(B)^{ab},F)\to 0.
\end{equation}
\end{prop}

\begin{proof}
Since $H_0^{\top}(A,\BZ)\cong \BZ\cong H_0^{\top}(B,\BZ)$, the K\"unneth formula (\cite[Thm. 3B.6]{Hatcher}) implies two canonical isomorphisms
\begin{equation*}
    H_1^{\top}(A,\BZ)\oplus H_1^{\top}(B,\BZ)\xrightarrow{\sim}H_1^{\top}(A\times B,\BZ)
\end{equation*}
and
\begin{equation*}
    H_2^{\top}(A,\BZ)\oplus H_2^\top(B,\BZ)\oplus H_1^{\top}(A,\BZ)\otimes_{\BZ}H_1^{\top}(B,\BZ)\xrightarrow{\sim} H_2^{\top}(A\times B,\BZ),
\end{equation*}
where the morphisms are given by cross products of cycles. Moreover, one could check that for a connected pointed topological space $X$ with the base point $x\in 
X$, the following composition 
\begin{equation}\label{equ-com is id}
    H_n^{\top}(Y,\BZ)\cong H_n^{\top}(Y,\BZ)\otimes H_0^{\top}(X,\BZ)\xrightarrow{\times}H_n^{\top}(Y\times X,\BZ)\xrightarrow{\mathrm{projection}} H_n^{\top}(Y,\BZ)
\end{equation}
is identity for all $n\geq 0$. Indeed, this morphism coincides with the natural push-forward along the morphism $$\text{id}\colon Y\cong Y\times x\hookrightarrow Y\times X\xrightarrow{p}Y.$$

The universal coefficient theorem (see \cite[Thm. 3.3.2]{Hatcher}) says that there exists a natural exact sequence, for every abelian group $F$ and for every topological space $X$ 
\begin{equation}\label{equ-UCT for topo coho}
    0\to \mathrm{Ext}^1(H_{n-1}^{\top}(X,\BZ),F)\to H^n_{\top}(X,F)\to \mathrm{Hom}_\BZ(H_n^{\top}(X,\BZ),F)\to 0.
\end{equation}
Let $n=1$, noting that $\mathrm{Ext}^1(H_{0}^{\top}(X,\BZ),F)=0$ as the group $H_{0}^{\top}(X,\BZ)$ is free, then we have the following commutative diagram
\begin{equation*}
    \xymatrix{H^1_{\top}(A,F)\oplus H^1_{\top}(B,F)\ar[rr]^{\sim}\ar[d]^{(p_1^*,p_2^*)}& &\mathrm{Hom}_\BZ(H_1^{\top}(A,\BZ)\oplus H_1^{\top}(B,\BZ),F)\ar[d]^{({p_1}_*,{p_2}_*)}\\
    H^1_{\top}(A\times B,F)\ar[rr]^{\sim}&& \mathrm{Hom}_\BZ(H_1^{\top}(A\times B,\BZ),F).
    }
\end{equation*}
Thus, the isomorphism in (\ref{equ-top Kunn for H^1}) is proved, since the morphism $({p_1}_*,{p_2}_*)$, which is induced by the inverse of K\"unneth formula, is an isomorphism.

Similarly, for $n=2$, we have the following commutative diagram of exact sequence
\begin{equation*}
\begin{adjustbox}{max width=\textwidth}
    \xymatrix{0\ar[r]& \mathrm{Ext}^1(H_{1}^{\top}(A,\BZ)\oplus H_{1}^{\top}(A,\BZ),F)\ar[r]\ar[d]^{\sim}&H^2_{\top}(A,F)\oplus H^2_{\top}(B,F)\ar[d]\ar[r]&\mathrm{Hom}_\BZ(H_2^{\top}(A,\BZ)\oplus H_2^{\top}(B,\BZ),F)\ar[d]\ar[r]& 0\\
    0\ar[r]& \mathrm{Ext}^1(H_{1}^{\top}(A\times B,\BZ),F)\ar[r]&H^2_{\top}(A\times B,F)\ar[r]&\mathrm{Hom}_\BZ(H_2^{\top}(A\times B,\BZ),F)\ar[r]& 0,
    }
\end{adjustbox}
\end{equation*}
where all vertical arrows are induced by projections. Hence by (\ref{equ-com is id}), the first vertical arrow is an isomorphism and the last one is an injection with cokernel $\mathrm{Hom}_{\BZ}(H_1^{\top}(A,\BZ)\otimes_{\BZ}H_1^{\top}(B,\BZ),F)$. Thus, by snake lemma, we have the canonical exact sequence
\begin{equation*}
    0\to H^2_{\top}(A,F)\oplus H^2_{\top}(B,F)\xrightarrow{(p_1^*,p_2^*)} H^2_{\top}(A\times B,F)\xrightarrow{f_{A,B}}\mathrm{Hom}_{\BZ}(H_1^{\top}(A,\BZ)\otimes_{\BZ}H_1^{\top}(B,\BZ),F)\to 0.
\end{equation*}

Finally, noting that every connected CW complex is path-connected (\emph{cf.} after \cite[Prop. 4.A.4]{Hatcher}),  
the Hurewicz Theorem (see \cite[Thm. 2A.1]{Hatcher}) says $\pi_1^{\top}(A)^{ab}\xrightarrow{\cong }H_1(A,\BZ)$ and $\pi_1^{\top}(B)^{ab}\xrightarrow{\cong }H_1(B,\BZ)$.
Therefore, to obtain the exact sequence (\ref{equ-top kun for H^2}), it is enough to define
\begin{equation*}
    \lambda^{\top}_{A,B}\colon H^2_{\top}(A\times B,F)\xrightarrow{f_{A,B}}\mathrm{Hom}_{\BZ}(H_1^{\top}(A,\BZ)\otimes_{\BZ}H_1^{\top}(B,\BZ),F)\cong \Hom_\BZ(\pi_1^{\top}(A)^{ab}\otimes_{\BZ} \pi_1^{\top}(B)^{ab},F).
\end{equation*}
\end{proof}
As a consequence, for any  smooth integral complex variety $X$ and for any abelian group $F$, we have the natural isomorphism
\begin{equation}\label{Thm-Kunneth-eq3}
\oplus \pr_i^*\colon  \bigoplus_{i=1}^n H^1_\top(X,F)\iso H^1_\top(X^n,F)
\end{equation}
and a natural exact sequence
\begin{equation}\label{Thm-Kunneth-eq4}
0\to \bigoplus_{i=1}^n H^2_\top(X,F)\xrightarrow{ \oplus \pr_i^*}  H^2_\top(X^n,F)\xrightarrow{\oplus\lambda_{i<j}  }
\bigoplus_{1\leq i<j\leq n}  \Hom_{\BZ} ( \pi_1^\top(X)^{ab}\otimes \pi_1^\top(X)^{ab} ,F) \to 0, 
\end{equation}
where $\pr_i\colon X^n\to X$ is the $i$-th projection, $\pr_{i,j}\colon X^n\to X^2$ is the projection onto the $i$-th and $j$-th components and 
$$\lambda_{i<j}\colon    H^2_\top(X^2,F) \xrightarrow{\pr^*_{i,j}} H^2_\top(X^n,F)   \xrightarrow{\lambda_{X,X}} 
\Hom_{\BZ} ( \pi_1^\top(X)^{ab}\otimes \pi_1^\top(X)^{ab} ,F).$$

Now, let us consider their connections with group (co)homology. Let $A^u\to A $ be the universal cover of a connected pointed CW complex $A$. For any abelian group $F$, the  Cartan--Leray spectral sequence (see \cite[Thm. VII.7.9]{brown2012cohomology}) implies the following two natural spectral sequences relating to group (co)homology 
\begin{equation}\label{equ-CLss for homo}
     E^2_{p,q}\coloneqq H_p(\pi_1^{\top}(A),H_q^\top(A^u,F))\Rightarrow H_{p+q}^\top(A,F)
\end{equation}
and 
\begin{equation}\label{CLss for coho}
  E_2^{p,q}\coloneqq H^p(\pi_1^{\top}(A),H^q_\top(A^u,F))\Rightarrow H^{p+q}_\top(A,F). 
\end{equation}
Actually, for any cover of $A$, we have the similar spectral sequences. These two spectral sequences induce the following edge morphisms for all $n\geq 0$
\begin{equation}\label{top-ex-Hure}
\epsilon_n^{\top}\colon H_n^{\top}(A,F) \to H_n(\pi_1^{\top}(A),F) \quad \text{ and } \quad   \epsilon^n_{\top}\colon H^n(\pi_1^{\top}(A),F)\to H^n_\top(A,F). 
\end{equation}
We call these edge morphisms as Hurewicz maps and Hurewicz maps in cohomology.

Let $G$ be an abelian group. We have natural isomorphisms
\begin{equation}\label{equ-group homology for 012}
 \BZ\xrightarrow{\sim }H_0(G,\BZ), \qquad   G\xrightarrow{\sim} H_1(G,\BZ), \qquad\text{and}\qquad \wedge^2 G\xrightarrow{\sim}H_2(G,\BZ),
\end{equation}
where $\BZ$ has trivial $G$-action; see \cite[Thm. V.6.4]{brown2012cohomology}.

  One could define the following morphism (here we omit the coefficients $\BZ$) for connected pointed CW complexes $A$ and $B$
\begin{equation}\label{equ-a lot composition}
\begin{split}
     \pi_1^{\top}(A)^{ab}\otimes_{\BZ} \pi_1^{\top}(B)^{ab}&\xrightarrow{\sim} H_1^\top(A)\otimes_\BZ H_1^\top(B)\xrightarrow{\times} H_2^\top(A\times B)\\
     &\xrightarrow{\epsilon_2^\top} H_2(\pi_1^\top(A\times B))\xrightarrow{}H_2(\pi_1^\top(A\times B)^{ab})\xrightarrow{\sim} \wedge^2\pi_1^\top(A\times B)^{ab},
\end{split}
\end{equation}
where
\begin{enumerate}
    \item the first isomorphism is given by Hurewicz Theorem (see \cite[Thm. 2A.1]{Hatcher});
    \item the second one is cross product and is injective by K\"unneth formula (see \cite[Thm. 3B.6]{Hatcher});
    \item $\epsilon_2^\top$ is the edge morphisms of Cartan--Leray spectral sequence (see (\ref{top-ex-Hure}));
    \item the next arrow is induced by the group homomorphism $\pi_1(A\times B)\twoheadrightarrow \pi_1(A\times B)^{ab}$;
    \item the last isomorphism is the inverse of the last one in (\ref{equ-group homology for 012}).
\end{enumerate}
On the other hand, since connected CW complexes are path-connected (\emph{cf.} after \cite[Prop. 4.A.4]{Hatcher}), hence regarding fundamental groups \cite[Prop. 1.1.12]{Hatcher} says that 
\begin{equation*}
    \pi_1^{\top}(A\times B)\cong \pi_1^{\top}(A)\times \pi_1^{\top}(B),
\end{equation*}
and hence one has
\begin{equation*}
 \pi_1^{\top}(A\times B)^{ab}\cong \pi_1^{\top}(A)^{ab}\oplus \pi_1^{\top}(B)^{ab}.   
\end{equation*}
Note that for abelian groups $M$ and $N$, we have the following canonical isomorphism
\begin{equation}\label{equ-iso for wedge of abel}
 \begin{split}
      \Phi\colon \wedge^2 M\oplus (M\otimes_\BZ N)\oplus \wedge^2 N
&\longrightarrow \wedge^2(M\oplus N)\\
m_1\wedge m_2&\longmapsto (m_1,0)\wedge(m_2,0),\\
m\otimes n&\longmapsto (m,0)\wedge(0,n),\\
n_1\wedge n_2&\longmapsto (0,n_1)\wedge(0,n_2).  
\end{split}
\end{equation}

\begin{prop}\label{prop-top case for inj}
    Let $A$ and $B$ be two connected pointed CW complexes. The composite in \((\ref{equ-a lot composition})\) coincides, up to sign,  with the second component of the map $\Phi$ defined in \((\ref{equ-iso for wedge of abel})\). In particular, it is an injection.
\end{prop}
\begin{proof}
    Pick $[a]\in \pi_1^\top(A)^{ab}$ (\emph{resp.} $[b]\in \pi_1^\top(B)^{ab}$), where $a\in \pi_1^\top(A)$ (\emph{resp.} $b\in \pi_1^\top(B)$). Then the loops $a$ and $b$ define continuous maps of pointed spaces $a\colon S^1\to A$ and $b\colon S^1\to B$ such that these continuous maps map the generator $1\in \pi_1^\top(S^1)\cong \BZ$ to $a\in \pi_1^\top(A)$ and $b\in \pi_1^\top(B)$, respectively. Since every morphism in (\ref{equ-a lot composition}) is naturally commutative with morphisms of connected pointed CW complexes, we get the following commutative diagram
    \begin{equation*}
        \xymatrix{1\otimes 1\ar@{|->}[d]&\pi_1^\top(S^1)^{ab}\otimes_\BZ \pi_1^\top(S^1)^{ab}\ar[rr]^{(\ref{equ-a lot composition})}\ar[d]&&\wedge^2\pi_1^\top(S^1\times S^1)^{ab}\ar[d]& (1,0)\wedge (0,1)\ar@{|->}[d]\\
       [a]\otimes [b] &\pi_1^\top(A)^{ab}\otimes_\BZ \pi_1^\top(B)^{ab}\ar[rr]^{(\ref{equ-a lot composition})}&&\wedge^2\pi_1^\top(A\times B)^{ab}& [(a,0)\wedge (0,b)].
        }
    \end{equation*}
Here we note that
\begin{itemize}
    \item $\pi_1^\top(S^1)\cong \BZ\cong \pi_1^\top(S^1)^{ab}$ and hence $\pi_1^\top(S^1)^{ab}\otimes_\BZ \pi_1^\top(S^1)^{ab}\cong \BZ$ with a generator $1\otimes 1$;
    \item $\wedge^2\pi_1^\top(S^1\times S^1)^{ab}\xrightarrow[\sim]{\Phi^{-1}} \pi_1^\top(S^1)^{ab}\otimes_\BZ \pi_1^\top(S^1)^{ab}\cong \BZ$ with a generator $(1,0)\wedge(0,1)$;
    \item the morphism $\Phi$ is canonical so that $(1,0)\wedge (0,1)$ maps to $[(a,0)\wedge (0,b)]$ via the morphism $a\times b\colon S^1\times S^1\to A\times B$.    
\end{itemize}
Therefore, it suffices to show that (\ref{equ-a lot composition}) is an isomorphism for the case $A=B=S^1$, so that $1\otimes 1$ must map to $\pm (1,0)\wedge (0,1)$. Note that the morphism (\ref{equ-a lot composition}) is also canonical, so we can only have $1\otimes 1\mapsto (1,0)\wedge (0,1)$ or $1\otimes 1\mapsto -(1,0)\wedge (0,1)$, but not both.

It is enough to check that the morphisms ``$\times$" and ``$\epsilon_2^\top$" in (\ref{equ-a lot composition}) are isomorphic for $A=B=S^1$. The first is because $H_2(S^1,\BZ)=0$ and the K\"unneth formula (\emph{cf.} \cite[Thm. 3B.6]{Hatcher}). The later one is because that $S^1\times S^1$ is a classifying space of $\BZ\oplus \BZ$, so the universal cover $(S^1\times S^1)^u$ is contractible, which means that $H_{>0}^\top((S^1\times S^1)^u,\BZ)=0$. Hence, the Cartan--Leray spectral sequence degenerates and $\epsilon_n^\top$ is an isomorphism for all $n\geq 0$. The proof is now complete.
\end{proof}

From now on, let's modify the sign of the isomorphism in (\ref{equ-group homology for 012})
such that Proposition \ref{prop-top case for inj} holds without ``up to sign''.

\subsection{Topological Dold--Thom theorem}
In this subsection, we recall the main results of the celebrated Dold--Thom theorem (see \cite{DT58}). Moreover, we give some conclusions of (co)homology of symmetric powers on low degrees. These facts are well-known for experts, and we include the details here for the convenience of the reader.

Let $A$ be a pointed connected CW complex with a fix point $e\in A$.
Let $\Sym^nA$ be the $n$-th symmetric power of $A$, i.e., 
the quotient space of the product of $n$-fold by the permutation action of $S_n$.
The inclusion 
\begin{equation*}
  A^n\subset A^{n+1}\colon  (x_1,\cdots, x_n)\mapsto (x_1,\cdots, x_n,e)  
\end{equation*}
induces an inclusion $\Sym^nA\subset \Sym^{n+1}A$, and we define $\Sym A\coloneqq \cup_n\Sym^nA$.
\begin{thm}[cf.{\cite{DT58}}]\label{thm- toplogical DT}
Let $A$ be a pointed connected CW complex. 
\begin{itemize}
\item[(a)] There exists a natural isomorphism $\pi_n(\Sym A)\to {H}_n(A,\BZ)$ for all $n\geq 1$ which is compatible with Hurewicz homomorphism $\pi_n(A)\to {H}_n(A,\BZ)$ via the morphism $A=\Sym ^1 A\to \Sym A$. 
\item[(b)] The H-space $\Sym A$  has the weak homotopy type of a product of Eilenberg--MacLane spaces with
$\pi_i(\Sym A) \cong H_i(A,\BZ)$. That is, 
\begin{equation}\label{equ-EM space}
    \Sym A\simeq_w\prod_{n\geq1} K(H_n(A,\BZ),n).
\end{equation}
\item[(c)] For $n\geq 2$, the natural morphism $\pi_1(\Sym^nA)\to \pi_1(\Sym A)\cong H_1(A,\BZ)$ induced by the inclusion $\Sym^nA\subset \Sym A$ is an isomorphism.
\item[(d)] For all $n\geq 1$ and all $i\geq 0$, the natural morphism $H_i(\Sym^{n-1}A,\BZ)\to H_i(\Sym^{n}A,\BZ)$ induced by the inclusion $\Sym^{n-1} A\subset \Sym^n A$ is a split injection.

\end{itemize}
\end{thm}
\begin{proof}
  The (a) and (b) are the main results of \cite{DT58}; and (c) can be found in \cite[12.15.Satz]{DP61}. One could find (d) in \cite[Thm. 3.3]{kallel2008symmetric}.
\end{proof}
\begin{cor}\label{cor-topo DT}
Let $A$ be a pointed connected CW complex and let $F$ be an abelian group. For all $n\geq 1$ and all $i\geq 0$, considering the closed inclusion $\Sym^{n-1} A\subset \Sym^n A$,
\begin{itemize}
    \item[(a)] the natural morphism $H_i(\Sym^{n-1}A,F)\to H_i(\Sym^{n}A,F)$ induced by this inclusion is injective;
    \item[(b)] the natural morphism $H^i(\Sym^{n}A,F)\to H^i(\Sym^{n-1}A,F)$ induced by this inclusion is surjective.
\end{itemize}
\end{cor}
\begin{proof}
    These follow from the universal coefficient theorem and Theorem \ref{thm- toplogical DT} (d); see \cite[Thm. 3.3.\& Thm. 3A.3]{Hatcher}.
\end{proof}

\subsection{\'Etale fundamental group and Hurewicz map in algebraic geometry}

In this subsection, we first recall the basic property of the \'etale fundamental groups (see \cite{Sz} for more details) and Hochschild--Serre spectral sequences. Furthermore, we compare the Cartan--Leray spectral sequences with Hochschild--Serre spectral sequences over the complex number field.

Firstly, for any normal integral noetherian scheme $X$ and any non-empty open subset $U\subset X$,
the induced homomorphism
\begin{equation}\label{eq-fund-gps-1}
\pi_1(U)\twoheadrightarrow \pi_1(X)
\end{equation}
is surjective, and it is an isomorphism if $\codim(X\setminus U,X)\geq 2$ and if $X$ is smooth.

Secondly, let $X,Y$ be two integral $k$-varieties.
Then we have the natural isomorphism
\begin{equation}\label{eq-fund-gps-2}
\pi_1(X_{\bk}\times_{\bk} Y_{\bk})\iso \pi_1(X_{\bk})\times \pi_1(Y_{\bk}).
\end{equation}

Thirdly, if $k=\BC$ is the complex field and $X$ is connected of finite type,  
then $X(\BC)$ is a complex algebraic manifold and we have the comparison theorem (\cite[Thm. 5.7.4]{Sz})
\begin{equation}\label{eq-fund-gps-5}
\widehat{\pi_1^{\top}(X(\BC))} \iso \pi_1(X).
\end{equation}

Fourthly, for any geometrically integral $k$-variety $X$,
 we have the canonical homotopy exact sequence for fundamental groups (\cite[Thm. IX.6.1]{SGA1} or \cite[Prop. 5.6.1]{Sz})
\begin{equation}\label{eq-fund-gps-4}
1\to \pi_1(X_{\bk})\to \pi_1(X)\to \Gamma_k\to 1.
\end{equation}
Moreover, this  homotopy exact sequence has the following generalization (\cite[Lem. 2]{KL81}, see also \cite[Ch. 5, Ex. 6, pp. 208]{Sz}).

\begin{prop}[Katz, Lang]\label{prop-Katz-Lang}
Let $S$ be an integral normal noetherian scheme and $\eta\in S$ the generic point.
 Let $f\colon  X\to S$ be a smooth surjective morphism of finite type such that the geometric generic fibre $X_{\bar{\eta}}$ is connected.
Then we have
\begin{enumerate}
    \item $X$ is normal and connected;
    \item for geometric point $\bar{x}$ of $X_{\bar{\eta}}$ with image $x\in X$, we have exact sequence 
$$ \pi_1(X_{\bar{\eta}},\bar{x})\to \pi_1(X,\bar{x})\to \pi_1(S,\bar{\eta})\to 1 .$$ 
\end{enumerate}
\end{prop}

Finally, an important fact is that the fundamental groups satisfy the Lefschetz Hyperplane Theorem. 
The following proposition is a general version of this result.

\begin{prop}[\emph{cf.}{\cite[Prop. 6.2 and Rem. 6.3]{EK16}}]\label{prop-Lef}
Let $X$ be a quasi-projective smooth geometrically integral variety over a field $k$ of characteristic zero. 
Then there exists a smooth geometrically integral sub-curve $C\subset X$ 
such that the homomorphism of the \'etale fundamental groups
$\pi_1(C)\to \pi_1(X)$ is surjective. If $k=\BC$, there exists a smooth geometrically integral sub-curve $C\subset X$ such that the natural morphism $\pi_1^\top(C)\to \pi^\top_1(X)$ is surjective.
\end{prop}

\begin{proof}
    For topological cases, it follows from \cite[1.2, Thm. p.153]{GR88}. For every field of characteristic zero, by \cite[Prop. 6.2 and Rem. 6.3]{EK16}, there exists a smooth geometrically integral sub-curve $C\subset X$ such that the homomorphism of tame fundamental groups
$\pi_1^{tame}(C)\to \pi_1^{tame}(X)$ is surjective. For the definition of tame fundamental groups, see \cite[Def.~5.7.15]{Sz}. 
We finish our proof by noting that the tame fundamental groups is the same as \'etale fundamental groups since $\chara(k)=0$. The proof is complete now.
\end{proof}

Let $Y\to X$ be a finite \'etale Galois cover of a $k$-scheme $X$ with Galois group $\Gal(Y/X)$. Then we have the Hochschild--Serre spectral sequence for any torsion sheaf $F$ over $k_{\et}$
\begin{equation*}
    E_2^{p,q}\coloneqq H^p(\Gal(Y/X),H^q(Y,F))\Rightarrow H^{p+q}(X,F).
\end{equation*}
Since cohomology commutes with inverse limits of schemes and with inverse limits of groups, and direct limits are exact, by passing to the direct limit over an inverse system of finite \'etale Galois covers $\{Y_i/X\}_{i\in I}$, we form a spectral sequence
\begin{equation}\label{equ-HS ss}
    E_2^{p,q}\coloneqq \varinjlim_{\{Y_i/X\}_{i\in I}}H^p(\Gal(Y_i/X),H^q(Y_i,F))\Rightarrow H^{p+q}(X,F).
\end{equation}
for any torsion sheaf $F$ over $k_{\et}$; see \cite[Rem.~III.2.21(b)]{Mi80}. On the other hand, for connected $k$-scheme $X$, one could choose the inverse system as in \cite[\S I.5, pp. 39]{Mi80}; then we have
\begin{equation*}
   \varinjlim_{Y/X}H^p(\Gal(Y/X),H^q(Y,F))\cong H^p(\varprojlim_{Y/X}\Gal(Y/X),\varinjlim_{Y/X}H^q(Y,F))\cong H^p_{cont}(\pi_1(X),\varinjlim_{Y/X}H^q(Y,F)) 
\end{equation*}
 by the definition of \'etale fundamental groups, where $Y$ runs through all finite \'etale Galois cover (with a fixed point over the fixed geometric point $\bar{x}\to X$) and we take continuous group cohomology. Therefore, for this case, the spectral sequence (\ref{equ-HS ss}) is exactly the following spectral sequence  
\begin{equation}\label{equ-HS ss2}
  E_2^{p,q}\coloneqq   H^p_{cont}(\pi_1(X),\varinjlim_{Y/X}H^q(Y,F))\Rightarrow H^{p+q}(X,F),
\end{equation}
where $Y/X$ passes all finite \'etale Galois cover.

 Similarly to (\ref{top-ex-Hure}), for any torsion sheaf $F$ over $k_{\et}$, 
the spectral sequence (\ref{equ-HS ss2}) induces the Hurewicz map for all $i\geq0$
\begin{equation}\label{equ-algbraic Hur map}
    \epsilon_{X,F}^i\colon  H^i_{cont}(\pi_1(X),F)\to H^i(X,F).
\end{equation}
For simplicity, we write
\begin{equation}\label{eq-def HS ss for F}
    \epsilon_{X,F}\colon  H^2_{cont}(\pi_1(X),F)\to H^2(X,F).
\end{equation}
As in (\ref{def-epsilon-eq}), we denote $\epsilon_X\coloneqq \epsilon_{X,\mu_{\infty}} $.

\begin{prop}\label{Prop-Hur-1}
Let $X$ be a quasi-projective smooth geometrically integral variety over $k$ and let $F$ be a torsion sheaf over $k_{\et}$.
For any $\alpha\in H^2(X,F)$, the following statements are equivalent
\begin{enumerate}
\item  $\alpha\in \Im(\epsilon_{X,F})$;
\item  $\alpha|_{X_{\bk}}$ belongs in the image of 
$\epsilon_{X_{\bk},F}\colon  H^2_{cont}(\pi_1(X_{\bk}),F)\to H^2(X_{\bk},F)$;
\item there exists a finite \'etale Galois cover $Y\to X$ such that $\alpha|_Y=0$;
\item there exists a finite \'etale Galois cover $Z\to X_{\bk}$ such that $\alpha|_Z=0$.
\end{enumerate}
\end{prop}

\begin{proof}
Consider the spectral sequence (\ref{equ-HS ss2}).
By (\ref{eq-fund-gps-3}), for any $\alpha\in H^1(Y,F)$,
there exists a  finite \'etale Galois cover $Y'/Y/X$ such that $\alpha|_{Y'}=0$.
Therefore, $E_2^{p,1}\cong 0$ for all $p$. 
Then the spectral sequence (\ref{equ-HS ss2}) induces the exact sequence
\begin{equation}\label{Prop-Hur-1-e2}
0\to H^2_{cont}(\pi_1(X),F) \xrightarrow{\epsilon_{X,F}}  H^2(X,F) \to \varinjlim_{Y/X} H^2(Y,F).
\end{equation}

Now (1)$\Leftrightarrow$(3) follows from (\ref{Prop-Hur-1-e2}).
Applying this result to $X_{\bk}$, we have (2)$\Leftrightarrow$(4).
Moreover, (3)$\Rightarrow$(4) is clearly, and (4)$\Rightarrow$(3) follows from the fact
$H^2(Z,F)= \varinjlim_{Z/Y/X} H^2(Y,F)$, where $Y/X$ passes all finite \'etale cover, 
and each cover $Y/X$ has a finite \'etale cover $Y'/Y$ which is Galois over $X$.
\end{proof}

Let us compare the Cartan--Leray spectral sequence (\ref{CLss for coho}) with the spectral sequence (\ref{equ-HS ss2}) under Artin's comparison theorem (\cite[Cor. III.3.12]{Mi80}).

\begin{prop}\label{prop-compatible of edge morphisms}
 Let $X$ be a connected smooth variety over $\BC$, and let $F$ be a torsion abelian group. Let $X^u\to X$ be the universal cover of $X$ as CW complexes. Starting from page 2, there exists a natural morphism of spectral sequences between (\ref{equ-HS ss2}) and (\ref{CLss for coho})
 \begin{equation*}
 E_2^{p,q}\coloneqq H^p_{cont}(\pi_1(X),\varinjlim_{Y/X}H^q(Y,F))\xrightarrow{f_2^{p,q}} H^p(\pi_1^\top(X),H^q_\top(X^u,F))  \eqqcolon E_2^{'p,q}.
 \end{equation*}
 Moreover, if the topological fundamental group $\pi_1^\top(X)$ is abelian, 
 then $f_2^{p,0}$ is an isomorphism for all $p$.
\end{prop}

\begin{proof}
Let $f\colon Y\to X$ be a finite \'etale Galois cover of $X$ with Galois group $G_Y\coloneqq\Gal(Y/X)$. Then the morphism $f(\BC)\colon Y(\BC)\to X(\BC)$ is a finite topological cover with the deck transformation group $G_Y$. For the category of torsion \'etale sheaves, the following diagram commutes
\begin{equation*}
\xymatrix{
  \operatorname{Sh}(X_{\acute{e}t})\ar[dd] ^{F\mapsto F^{an}} \ar[drr]^{F\mapsto F(Y)}& & & &\\
 && G\text{-}\operatorname{Mod}\ar[rr]^{M\mapsto M^G} & &\operatorname{Ab}\\
  \operatorname{Sh}(X_{\mathrm{an}}). \ar[urr]_{E\mapsto f^*E(Y^{an})}& & & &
}
\end{equation*}
Since that both that Cartan--Leray spectral sequence and Hochschild–Serre spectral sequence can be viewed as the Grothendieck's spectral sequences of composition of functors, by Artin's comparison theorem (\cite[Cor. III.3.12]{Mi80}), we can identify the Hochschild–Serre spectral sequence of the finite Galois cover $f$ coincides with the Cartan--Leray spectral sequence of $f(\BC)$ (see \cite[Thm. VII.7.9]{brown2012cohomology})
\begin{equation*}
   E_2^{p,q}(Y_{\et})\coloneqq  H^p(G_Y,H^q_{\et}(Y,F))\xrightarrow{\cong} H^p(G_Y,H^q_\top(Y,F)) \eqqcolon E_2^{p,q}(Y_\top).
\end{equation*}
One the other hand, by the universal property of topological universal covers, there exist natural morphism of spectral sequence for every finite \'etale Galois cover $Y\to X$
\begin{equation*}
E_2^{p,q}(Y_\top)\coloneqq  H^p(G_Y,H^q_\top(Y,F))\xrightarrow{g_2^{p,q}}  H^p(\pi_1^\top(X),H^q_\top(X^u,F))  \eqqcolon E_2^{'p,q}.
\end{equation*}
Let $Y$ run through all finite \'etale Galois cover of $X$. There exist natural morphism of spectral sequences $f_2^{p,q}$ such that the following commutes
\begin{equation*} 
    \xymatrix{\varinjlim_{Y/X} E_2^{p,q}(Y_{\et})\coloneqq  \varinjlim_{Y/X}  H^p(G_Y,H^q_{\et}(Y,F))\ar[r]^-{\cong}\ar[d]^{\cong}& \varinjlim_{Y/X} H^p(G_Y,H^q_\top(Y,F)) \eqqcolon \varinjlim_{Y/X} E_2^{p,q}(Y_\top)\ar[d]^{g_2^{p,q}}\\
    E_2^{p,q}\coloneqq H^p_{cont}(\pi_1(X),\varinjlim_{Y/X}H^q(Y,F))\ar@{-->}[r]^{f_2^{p,q}}& H^p(\pi_1^\top(X),H^q_\top(X^u,F))  \eqqcolon E_2^{'p,q},
    }
\end{equation*}
which is exactly we desired one.

Now we consider the $E_2^{p,0} $ part of the above diagram.
\begin{equation}  \label{prop-compatible of edge morphisms-e}
\xymatrix{  \varinjlim_{Y/X}  H^p(G_Y,H^0_{\et}(Y,F))\ar[r]^-{\cong}\ar[d]^{\cong}& \varinjlim_{Y/X} H^p(G_Y,H^0_\top(Y,F))  \ar[d]^{g_2^{p,0}}\\
H^p_{cont}(\pi_1(X),\varinjlim_{Y/X}H^0(Y,F))\ar[r]^{f_2^{p,0}}&   H^p(\pi_1^\top(X),H^0_\top(X^u,F)) 
}
\end{equation}
In our inverse system $\{Y/X\}$, the sub-system that $Y$ is connected is cofinal, hence provides the same direct limit.
Therefore, we may assume all our $Y$ are connected.
In this case, the space $Y(\BC)$ is also connected and $H^0(Y,F)\cong F \cong H^0_{\top}(Y,F) $.
Then the diagram (\ref{prop-compatible of edge morphisms-e}) is exactly the $H^p(-,F)$ (resp. $H^p_{cont}(-,F)$) of 
\begin{equation*} 
\xymatrix{  \varprojlim_{Y/X}  G_Y   &  \varprojlim_{Y/X}   G_Y   \ar[l]_=  \\
\pi_1(X)  \ar[u]^{\cong}  &   \pi_1^\top(X) \ar[u]\ar[l]_g  .
}
\end{equation*}
By (\ref{eq-fund-gps-5}), the homomorphism $f_2^{p,0}$ is exactly the morphism induced by the profinite completion
$$\pi_1^\top(X)\to \pi_1(X)\cong\widehat{\pi_1^{\top}(X)}.$$

It remains to show that $f_2^{p,0}$ is an isomorphism for all $p$ if $\pi_1^\top(X)$ is abelian. 
Since every connected complex algebraic variety has the homotopy type of a finite CW-complex, its topological fundamental group is finitely presented, hence finitely generated. Thus, there exists a unique $r\in \BN$ such that we have $\pi_1^\top(X)\cong \BZ^{\oplus r}\oplus M$, where $M$ is a finite abelian group. Therefore, $\pi_1(X)\cong \hat{\BZ}^{\oplus r}\oplus M$. 
Since $\pi_1^\top(X)$ acts on the group $F$ trivially, the morphism $g$ equals the map defined componentwise $$H^i_{cont}(\hat{\BZ}^{\oplus r},F)\oplus H^i_{cont}(M,F)\xrightarrow{(g_1,g_2)} H^i(\BZ^{\oplus r},F)\oplus H^i(M,F).$$
Since $\BZ^{\oplus r}$ is a finitely generated free group so it is a Serre's good group (see \cite[\S 2.2.6]{serre2013galois}), $g_1$ is an isomorphism;  the morphism $g_2$ is also isomorphic as $M$ is finite. Therefore, $g$ is an isomorphism. The proof is complete.
\end{proof}

\section{The fibrations from symmetric powers of curves to relative Picard schemes}\label{section3}

Let $C$ be a quasi-projective smooth geometrically integral curve over $k$.
We define
\begin{equation}\label{eq-def-Nc}
  N_C\coloneqq \left\{\begin{array}{ll}
      2\cdot g_{\Bar{C}}+\mathrm{deg}(\Bar{C}\setminus C)-2 & \text{if } g_{\Bar{C}}>0 \\
      \mathrm{deg}(\Bar{C}\setminus C) & \text{if } g_{\Bar{C}}=0,
  \end{array} \right.
\end{equation}
where $C\subset \bar{C}$ is the smooth compactification of $C$ and $g_{\bar{C}}$ is the geometric genus of $\bar{C}$. 

\subsection{Projective cases}
Let us first recall some background of Picard functor for a smooth projective curve $C$; see \cite[Chap. 8]{Neron} for details. We define the \textbf{Picard functor}, denoted by ${Pic}_C$, as the $fppf$-sheaf associated to the functor
\begin{equation*}
    P_C\colon (\mathrm{Sch}/k)^{\rm op}\to (\mathrm{Sets}), \quad T\mapsto \mathrm{Pic}(C\times T).
\end{equation*}
For any integer $n\in \BZ$, we define the subfunctor $Pic^{n}_C$ of $Pic_C$ induced by the line bundles of degree $n$. Then the functor $Pic^{n}_C$ is represented by a projective $k$-scheme, denoted by $\Pic^n_C$. Moreover, the functor ${Pic}_C$ is represented by the disjoint union of all $\{\Pic_C^{n}\}_{n\in \BZ}$, and we denote it $\Pic_C$; see \cite[\S 8.2, Thm. 5 \& \S 8.4,  Thm.~3]{Neron}. In particular, the identity component $\Pic^0_C$ of $\Pic_C$ is an abelian variety over~$k$, and $\Pic_C^n$ is a torsor under $\Pic_C^0$ for all $n\in \BZ$; see \cite[\S 9.2, Prop. 3]{Neron}. Applying \cite[\S 8.2, Prop.~7]{Neron} to the morphism $C\to \mathrm{Spec}(k)$, one can get the following proposition. 

\begin{prop}[projective cases]\label{prop-proj case}
Let $C$ be a projective smooth geometrically integral curve over $k$. The morphism 
\begin{equation*}
 \Sym^nC\cong \mathrm{Hilb}^nC\cong \mathrm{Div}^n_C\to  \Pic^n_C    
\end{equation*}
is a projective bundle for all $n> N_C$, where 
\begin{enumerate}
    \item $\mathrm{Hilb}^nC$ is the Hilbert scheme of $C$ with degree $n$ and the first isomorphism is given by \cite[\S 9.3, Prop. 3]{Neron};
    \item $\mathrm{Div}^n_C$ represents the functor $Div^n_C\colon (\mathrm{Sch}/k)^{\rm op}\to (\mathrm{Sets})$ which associates to relative effective Cartier divisors of degree $n$;
    \item the morphism $\mathrm{Div}^n_C\to  \Pic^n_C$ is induced by the canonical morphism of functors $Div_C\to Pic_C$; see \cite[pp. 212]{Neron}.    
\end{enumerate}
\end{prop}

\begin{proof}
If $g_{C}=0$, then $C\cong \BP^1$, $\Sym^nC\cong \BP^n$ and $\Pic^n_C $ is a point. Thus this proposition holds obviously. We may assume that $g_C\geq 1$.

Let $T\coloneqq \Pic^n_C$ and let $\CL$ be the universal line bundle on $C_T=C\times \Pic^n_C$. Then $\CL$ induces the identity morphism $T\xrightarrow{\rm id} \Pic_C^n$ by representability of $Pic_C^n$. By \cite[\S 8.2, Prop.7]{Neron}, the functor $Div^n_C\times _{\Pic^n_C}T\cong Div_C^n$ is represented by a projective $T$-scheme $\BP_T(\CF)$ for some $\CO_T$-module $\CF$. Moreover, there exists a canonical choice of $\CF$. By \cite[\S 8.2, Prop.7]{Neron} again, it suffices to show that for $n> N_C$, the canonical $\CL$ is cohomologically flat in dimension zero (see \cite[pp. 206]{Neron}). Then we get that $\CF$ is a locally free $\CO_T$-module, and hence $\BP_T(\CF)\cong \Sym^n C$ is a projective bundle over $T=\Pic_C^n$.

Now, we are going to show that $\CL$ is cohomologically flat in dimension zero. Let $\CO_C(1)$ be an ample line bundle on $C$, and let $\CO_{C_T}(1)$ be the relative ample line over $T$ associated to the pullback $C\times T\to C$. Let $p\colon C\times \Pic_C^n\to \Pic^n_C$ be the projection and for a point $t\in T$, let $p_t\colon C_t\to {\rm Spec}(\kappa(t))$ be its restriction. If $n>N_C$, by Riemann--Roch theorem (\emph{cf.} \cite[\S 9.1, Thm. 1]{Neron}) and noting $\deg (\CL_t(m))\geq n>N_C$ for all $m\geq 0$, we get $H^1(C_t,\CL_t(m))=0$ for all $m\geq 0$, where $\CL_t$ is the restriction of $\CL$ to $C_t$; see \cite[Chap. IV, Ex. 1.3.4]{Hartshorne}. This implies that the function $t\mapsto \dim(H^i(C_t,\CL_t(m)))$ is the zero function on $T$ for all $i\geq 1$ and for all $m\geq 0$, since $H^{\geq 2}(C_t,\CL_t(m))\equiv0$ by Grothendieck vanishing theorem. Thus $p_*\CL$ is locally free on $T$ by \cite[Chap. III, Cor. 12.9]{Hartshorne} and $\CL$ is cohomologically flat in dimension zero since ${\rm rank}(p_*\CL)=\dim (H^0(C_t,\CL_t))=n-g+1>0$ for any $t\in T$; \emph{cf.} \cite[pp.222]{Neron}. Therefore, we can take $\CF$ to be the dual of $p_*\CL$, which is a locally free $\CO_T$-module.
\end{proof}

\subsection{Affine cases}
Let $C$ be an affine smooth geometrically integral curve over a field $k$ of characteristic zero (not necessarily algebraically closed), and let $C\subset \bar{C}$ be the smooth compactification of $C$ with boundary $Z$. First, let us recall the definition of \textbf{relative Picard functor}. The relative Picard functor, denoted by $Pic(\Bar{C},Z)$, is the $fpqc$-sheaf associated to the functor
\begin{equation*}
    P_{\Bar{C},Z}\colon (\mathrm{Sch}/k)^{\rm op}\to (\mathrm{Sets}), \quad T\mapsto \mathrm{Pic}(\Bar{C}\times T,Z\times T),
\end{equation*}
where $$\mathrm{Pic}(\Bar{C}\times T,Z\times T)\coloneqq H^1_{{\rm Zar/et}}(\Bar{C}\times T,\ker\{\BG_{m,\Bar{C}\times T}\to i_*(\BG_{m,Z\times T})\})$$ with inclusion $i\colon Z\times T\hookrightarrow \Bar{C}\times T$. 
Actually, the group $\mathrm{Pic}(\Bar{C}\times T,Z\times T)$ is the group of isomorphism classes of pairs $(\CL,\varphi)$, 
where $\CL$ is a line bundle on $\Bar{C}\times T$ and $\varphi\colon \CO_{Z\times T}\cong \CL|_{Z\times T}$ is a trivialization on $Z\times T$. Here we say $(\CL,\varphi)\sim (\CL',\varphi')$ if there exists an isomorphism $\phi\colon \CL\xrightarrow{\cong} \CL'$ such that 
$$\varphi'=\phi|_{Z\times T}\circ\varphi \colon  \CO_{Z\times T}\cong \CL'|_{Z\times T}.$$
Then the pre-sheaf $P_{\Bar{C},Z}$ is already an $fpqc$-sheaf (i.e., $P_{\Bar{C},Z}=Pic(\Bar{C},Z)$) and the functor $Pic(\Bar{C},Z)$ is represented by a $k$-group scheme $\Pic(\Bar{C},Z)$ which is locally of finite type over $k$; see \cite[Lem. 2.1.1 \& \S A.2]{BVS01}.
Also the natural transformation $Pic(\Bar{C},Z)\to Pic(\Bar{C})$ induced by $\Pic(\Bar{C}\times T,Z\times T)\to \Pic(\Bar{C}\times T)$ implies the morphism of $k$-group schemes $\Pic(\Bar{C},Z)\to \Pic(\Bar{C})$, 
and we write $\Pic^n(\Bar{C},Z)$ as the preimage of $\Pic^n(\Bar{C})$ under this morphism. Then $\Pic^0(\Bar{C},Z)$ is the identity component of $\Pic(\Bar{C},Z)$ and is a semi-abelian $k$-group scheme by \cite[Prop. 2.1.2]{BVS01}.

For every $k$-scheme $T$, the Hilbert functor $Hilb^nC$ maps $T$ to the set 
\begin{equation*}
    \left\{  W\mid \text{closed subscheme $W\subset C\times T$ such that the projection $W\to T$ is finite flat of degree $n$ }\right\}.
\end{equation*}
Thus, for any $W\in (Hilb^nC)(T)$, the scheme $W\subset C\times T\subset \Bar{C}\times T$ is closed, and $W$ induces an effective relative Cartier divisor on $\Bar{C}\times T$ over $T$ by \cite[\S 8.2, Lem. 6]{Neron}. Moreover, one can show that the line bundle associated with this effective relative Cartier divisor is exactly given by $\CO(W)\coloneqq \CI^{-1}$, where $\CI$ is the ideal sheaf of $W\subset \Bar{C}\times T$ (we sometimes write the ideal sheaf as $\CO(-W)$ for this case), and there exists a canonical inclusion $\psi\colon \CO_{\Bar{C}\times T}\hookrightarrow \CO(W)$ induced by the embedding $\CI\hookrightarrow \CO_{\Bar{C}\times T}$. Since $W\cap (Z\times T)=\emptyset$, the restriction $\varphi(W)\coloneqq\psi|_{Z\times T}$ is a trivialization on $Z\times T$. Therefore, the pair  $ (\CO(W),\varphi(W))\in P_{\Bar{C},Z}(T)=Pic(\Bar{C},Z)(T)$, and one can verify that the arrangement $W\mapsto (\CO(W),\varphi(W))$ gives a natural transformation $Hilb^nC\to Pic^n(\Bar{C},Z)$. This natural transformation induces a morphism of schemes  $\mathrm{Hilb}^nC\to   \Pic^n_{\Bar{C},Z}$ by representability of functors, and further, we have the morphism
\begin{equation}\label{equ-Hilbert to relative Pic}
 \theta\colon \Sym^nC\cong \mathrm{Hilb}^nC\to   \Pic^n_{\Bar{C},Z}.    
\end{equation}

\begin{lem}[\emph{cf.} {\cite[\S 8.3, Prop. 7]{Neron}}]\label{lem-affine case}
Let $C$ be an affine smooth geometrically integral curve over $k$, and let $C\subset \bar{C}$ be the smooth compactification of $C$ with boundary $Z$. Assume $T$ is an affine $k$-scheme with a $k$-morphism $f\colon T\to \Pic^n_{\Bar{C},Z}$ and $(\CL,\varphi)$ is the line bundle on $\Bar{C}\times T$ associated to the morphism~$f$. If $n>N_C$, then there exists a locally free $\CO_T$-module $\CF$ such that the fiber product $Hilb^nC\times _{\Pic^n_{\Bar{C},Z}}T$ induced by the pull-back of the morphism $\theta$ is represented by the vector bundle ${\rm Spec}(S(\CF^{\vee}))$.
\end{lem}

\begin{proof}
The representability of $Hilb^nC\times _{\Pic^n_{\Bar{C},Z}}T$ is given by the representabilities of the functors $Hilb^nC$ and $Pic^n_{\Bar{C},Z}$. Indeed, the functor $Hilb^nC\times _{\Pic^n_{\Bar{C},Z}}T$ is represented by ${\rm Hilb}^nC\times _{\Pic^n_{\Bar{C},Z}}T$.

 Moreover, the fiber product $Hilb^nC\times _{\Pic^n_{\Bar{C},Z}}T$ is isomorphic to the functor $D_{(\CL,\varphi)} \colon (\mathrm{Sch}/T)^{\rm op}\to (\mathrm{Sets})$, which associates a $T$-scheme $T'$ to the set of all closed subsets $W'\subset C\times T'$ such that 
 \begin{itemize}
     \item[(a)] the projection  $W'\to T'$ \text{ is flat and finite of degree } $n$;
     \item[(b)] $(\CO(W'),\varphi(W'))$ and $(\CL',\varphi')$ { give rise to the same element in } $Pic^n_{\Bar{C},Z}(T')$.
 \end{itemize}
 Here $\CL'$ (resp. $\varphi'$) denotes the pull-back of $\CL$ (resp. $\varphi$) along with the morphism $\Bar{C}\times T'\to \Bar{C}\times T$. 

\bigskip

\begin{claim}
Let $p\colon \Bar{C}\times T\to T$ be the natural projection and let $p'$ be its pull-back along $T'\to T$. Then there is a bijection $$\Ker \Big(\Gamma(T',p'_*\CL')\xrightarrow{\rm restriction}  \Gamma(T',p'_*(\CL'|_{Z\times T'})) \Big)\xrightarrow{\sim} D_{(\CL,\varphi)}(T').$$
\end{claim}
\begin{proof}[Proof of the Claim]
Indeed, the relation $(\CO(W'),\varphi(W'))\sim (\CL',\varphi')$ means that there exists a morphism $\phi\colon \CO(W')\xrightarrow{\cong} \CL'$ such that 
\begin{equation}\label{eq-lemma-affine case}
    \varphi'= \phi|_{Z\times T'}\circ \varphi(W')\colon \CO_{Z\times T'}\cong \CL'|_{Z\times T'}
\end{equation}
Note that the morphism $\varphi(W')$ comes from a canonical inclusion $\psi(W')\colon \CO_{\Bar{C}\times T'}\hookrightarrow \CO(W')$, i.e., we have  $\varphi(W')=\psi|_{Z\times T'}$. Therefore, we have $$\phi\circ \psi(W')\in \mathrm{Hom}_{\CO_{\Bar{C}\times T'}}(\CO_{\Bar{C}\times T'},\CL')=\Gamma(\Bar{C}\times T',\CL')=\Gamma(T',p'_*\CL'),$$ 
and the relation (\ref{eq-lemma-affine case}) means that
\begin{equation*}
    \varphi'=(\phi\circ \psi(W'))|_{Z\times T'}\in \mathrm{Hom}_{\CO_{Z\times T'}}(\CO_{Z\times T'},\CL'|_{Z\times T'})=\Gamma(Z\times T',\CL'|_{Z\times T'})=\Gamma(T',p'_*(\CL'|_{Z\times T'})).
\end{equation*}
Moreover, we may assume that $\phi={\rm id}$ up to a rescaling, and then there is no $\lambda\in \Gamma(T',p'_*(\CO_{\Bar{C}\times T'}^*))$ such that both $\psi(W')$ and $\lambda\cdot \psi(W')$ map to $\varphi'$ under the restriction map.  Therefore, the set $D_{(\CL,\varphi)}(T')$ is one-to-one correspond to the set of $\psi(W')\in \Gamma(T',p'_*\CL')$ such that the image of $\psi(W')$ under the linear morphism $\Gamma(T',p'_*\CL')\xrightarrow{\rm restriction}  \Gamma(T',p'_*(\CL'|_{Z\times T'})) $ is exactly $\varphi'$, i.e., the inverse image of $\varphi'$. 

To finish the proof of this claim, it suffices to show that $\varphi'$ lives in the image of the linear morphism $\Gamma(T',p'_*\CL')\to  \Gamma(T',p'_*(\CL'|_{Z\times T'})) $ for every $T$-scheme $T'$. Consider following commutative diagram
\begin{equation*}
    \xymatrix{\Gamma(T,p_*\CL)\ar[r]\ar[d]&\Gamma(T,p_*(\CL|_{Z\times T}))\ar[d] & \varphi \ar@{|->}[d]\\
    \Gamma(T',p_*\CL')\ar[r]&\Gamma(T,p_*(\CL'|_{Z\times T'}))& \varphi',
    }
\end{equation*}
where vertical arrows are pull-backs, horizontal arrows are restrictions and all morphisms are linear. Hence, it suffices to show that the morphism  $\Gamma(T,p_*\CL)\to \Gamma(T,p_*(\CL|_{Z\times T}))$ is surjective and further the surjectivity of the morphism $p_*\CL\to p_*(\CL|_{Z\times T'})$ for all $n>N_C$, since the global section functor on affine schemes is exact.

By Nakayama's lemma, it is enough to show that for any field extension $K/k$ with a morphism ${\rm Spec}(K)\to T\xrightarrow{f}\Pic^n_{\Bar{C},Z}$, we always have the surjectivity of $H^0(\Bar{C}_K,\CL_K)\to H^0(Z_K,\CL_K|_{Z_K})$, where $\CL_K$ is the pull-back of $\CL$ along ${\rm Spec}(K)\to T$. Write $i\colon Z_K\hookrightarrow \Bar{C}$ and consider the following short exact sequence of sheaves of $\CO_{\Bar{C}_K}$-modules
\begin{equation*}
    0\to \CL_K(-Z_K)\to \CL_K\to i_*(\CL_K|_{Z_K})\to0.
\end{equation*}
Taking Zariski cohomology groups, we get the following exact sequence
\begin{equation*}
     0\to H^0(\Bar{C}_K, \CL_K(-Z_K))\to H^0(\Bar{C}_K,\CL_K)\to H^0(Z_K,\CL_K|_{Z_K}) \to H^1(\Bar{C}_K, \CL_K(-Z_K))\to\cdots.
\end{equation*}
Then for every $n>N_C$, since $\deg_K(Z_K)=\deg (\Bar{C}\setminus C)$ and $\deg_K(\CL_K)=n$, we have 
\begin{equation*}
    \deg_K(\CL_K(-Z_K))>\max\{0,2\cdot g_{\Bar{C}}-2\},
\end{equation*}
which implies that $H^1(\Bar{C}_K, \CL_K(-Z_K))=0$ by Riemann--Roch theorem (\emph{cf.} \cite[Thm. IV.1.1.3]{Hartshorne}). Therefore, for all $n>N_C$, the morphisms $p_*\CL\to p_*(\CL|_{Z\times T'})$ are surjective. The proof of this claim is now complete.
\end{proof}
 
\bigskip

Note that we have the following isomorphism
\begin{equation*}
   \Ker \Big(\Gamma(T',p'_*\CL')\xrightarrow{\rm restriction}  \Gamma(T',p'_*(\CL'|_{Z\times T'})) \Big)\cong \Gamma(T',p'_*(\CL'(-Z\times T'))),
\end{equation*}
by applying the left exact functors $p_*$ and $\Gamma(T',-)$ successively to the exact sequence
\begin{equation*}
    0\to \CL'(-Z\times T')\to \CL'\to \iota_*(\CL'|_{Z\times T'})\to0
\end{equation*}
with $\iota\colon Z\times T'\to \Bar{C}\times T'$. Thus we have
\begin{equation}
  {\rm Hom}_{\CO_{T'}}(\CO_{T'},p'_*(\CL'(-Z\times T'))) =\Gamma(T',p'_*(\CL'(-Z\times T')))\xrightarrow[\sim]{1-1}D_{(\CL,\varphi)}(T').
\end{equation}
Let $\CO_{\Bar{C}}(1)$ be an ample line bundle on $\Bar{C}$ and let $\CO_{\Bar{C}_T}(1)$ be the associated relative ample line bundle to $p\colon \Bar{C}\times T\to T$.  For any point $t\in T$, let $p_t\colon \Bar{C}_t\to {\rm Spec}(\kappa(t))$ be its restriction. Since for every $m\geq 0$, $$\deg(\CL_t(m)\otimes \CO_{\Bar{C}_t}(-Z_t))\geq \deg(\CL_t(-Z_t))>\max\{0,2\cdot g_{\Bar{C}}-2\},$$
by Riemann--Roch theorem again, $H^1(\Bar{C}_t,\CL_t(m)\otimes \CO_{\Bar{C}_t}(-Z_t))=0$. This implies that the coherent shaves $p_*(\CL(m)\otimes \CO_{\Bar{C}\times T}(-Z\times T))$ are locally free for all $m\geq 0$, and hence $\CL(-Z\times T)$ is flat over $T$; see \cite[Prop 4.37 \& Lem. 5.5]{FGAexplained}. Then the functor $ Hilb^nC\times _{\Pic^n_{\Bar{C},Z}}T\colon (\mathrm{Sch}/T)^{\rm op}\to (\mathrm{Sets})$ is equal to the functor
\begin{equation*}
\begin{split}
     T'\mapsto &{\rm Hom}_{\CO_{T'}}(\CO_{T'},p'_*(\CL'(-Z\times T')))\\
   \cong  &{\rm Hom}_{\CO_{\Bar{C}\times T'}}(\CO_{\Bar{C}\times T'},\CL'(-Z\times T'))\\
   \cong  &{\rm Hom}_{\CO_{\Bar{C}\times T'}}(\eta^*(\CO_{\Bar{C}\times T}),\eta^*(\CL(-Z\times T))),
\end{split}
\end{equation*}
where $\eta\colon \Bar{C}\times T'\to \Bar{C}\times T$. Therefore, by \cite[Thm. 5.8]{FGAexplained} and the flatness of $\CL(-Z\times T)$,  the functor $ Hilb^nC\times _{\Pic^n_{\Bar{C},Z}}T$ is represented by an affine vector bundle ${\rm Spec}(S(\CF^{\vee}))$ with locally free coherent sheaf $\CF=p_*(\CL(-Z\times T))$. 
\end{proof}

Applying Lemma \ref{lem-affine case} to every affine subscheme $f\colon T\hookrightarrow \Pic^n_{\Bar{C},Z}$, one gets that the morphism $\theta$ is an affine bundle for all $n>N_C$. 

\begin{prop}[affine cases]\label{prop-affine case}
 Let $C$ be an affine smooth geometrically integral curve over $k$. Let $C\subset \bar{C}$ be the smooth compactification of $C$ with boundary $Z$.  The morphism $\theta$ defined in (\ref{equ-Hilbert to relative Pic}) is an affine bundle for all $n> N_C$.
\end{prop}

\begin{rem}
 From Proposition \ref{prop-affine case}, $\Pic^n_{\Bar{C},Z}$ is geometrically integral for $n>N_C$, because this holds for $\Sym^nC$.
Since $\Pic_{\Bar{C},Z}$ is a group scheme, all $\Pic^i_{\Bar{C},Z}$ are isomorphic over $\bk$, 
hence are all geometrically integral.
Moreover, $\Pic^0_{\Bar{C},Z}$ is the identity component of $\Pic(\Bar{C},Z)$
and is a connected semi-abelian $k$-group scheme by \cite[Prop. 2.1.2]{BVS01}.   
\end{rem}

\subsection{Cohomology of symmetric powers of curves}

In this part, we apply our fibration to calculate the cohomology of symmetric powers of curves. 
We first consider the topological case, which need the following results about topological properties of abelian Lie groups.

\begin{lem}\label{exam-Hure-top-lie-gp}
Let $G$ be a connected finite dimensional abelian Lie group. Then 
\begin{itemize}
    \item[(i)] the universal cover of $G$ is contractible and $\pi_1^{\top}(G)\cong \BZ^d$ for some $d\geq 0$;
    \item[(ii)] the edge morphism $\epsilon_n^{\top}$ and $\epsilon^n_{\top}$ are isomorphism for all $n\geq 0$;
    \item[(iii)] $H_n(G,\BZ)\cong H_n(\pi_1^\top(G),\BZ)\cong \wedge^n \BZ^d$ is a finitely generated free abelian group for all $n\geq 0$.
\end{itemize}
\end{lem}

\begin{proof}
Actually, let $\mathfrak{g}$ be the corresponding Lie algebra of $G$. Then, for any $a,b\in \mathfrak{g}$, the Lie bracket $[a,b]=0$, and hence the Dynkin product $a\ast b$ in \cite[Def. 5.7.3]{Donne25} is equal to the sum $a+b$.
By \cite[Thm.~9.4.22]{Donne25}, the exponential map $\mathrm{exp}\colon  \mathfrak{g}\to G$ is a universal cover of $G$ with kernel $\Ker(\mathrm{exp})\subset \mathfrak{g}$ a discrete subgroup. Then (i) follows from the fact that $\mathfrak{g}$ is a finite dimensional real (or complex) vector space and hence the topological space $\mathfrak{g}$ is contractible, and the discrete subgroup $\Ker(\mathrm{exp})\cong \BZ^d$ for an integer $d$. Moreover, (ii) is induced by (i), and then (iii) is obtained by \cite[Thm. V.6.4]{brown2012cohomology}.
\end{proof}

\begin{prop}\label{lem-DT-1}
Let $C$ be an affine smooth integral curve over the complex number field $\BC$. For every abelian group $A$ and for any $n> N_C$, we have
\begin{itemize}
    \item[(i)]  the topological fundamental group $\pi_1^\top(\Sym^n C)$ is a finitely generated free abelian group;
    \item[(ii)] the Hurewicz map $\epsilon_\top^i\colon  H^i(\pi_1^\top(\Sym^nC),A)\to H^i_\top(\Sym^nC,A)$ is an isomorphism 
    for all $i\geq 0$;
    \item[(iii)]  for any fixed closed point $P$ of $C$, the morphism $$l_P\colon   C\xrightarrow{id_C\times P\times \cdots \times P} C^n\xrightarrow{p} \Sym^nC $$ induces a surjective homomorphism $l_{P,*}\colon \pi_1^\top(C)\to \pi_1^\top(\Sym^nC) $;
    \item[(iv)]  the canonical morphism $p_*\colon \pi_1^\top(C^n)\to \pi_1^\top(\Sym^nC)$ is split surjective and the natural morphism $p^*\colon H^i_\top(\Sym^nC,A)\to H^i_\top (C^n,A)$ is split injective for all $i\geq 0$.
\end{itemize}
\end{prop}
\begin{proof}
Let $\Pic_{\bar{C},Z}$ be the relative Picard scheme as above.
Let 
$$\theta\colon  \Sym^nC\to \Pic^n_{\bar{C},Z} $$
be the morphism in (\ref{equ-Hilbert to relative Pic}).
From Proposition \ref{prop-affine case}, 
$\theta$ is an affine bundle.

Fix a closed point $P\in C(\BC)$ and we denote by $Q\coloneqq \theta( P,\cdots,P)$.
Since $\Pic_{\bar{C},Z}$ is a group scheme,
the morphism $\Pic^0_{\bar{C},Z} \to \Pic^n_{\bar{C},Z}\colon  t\mapsto t+ Q$ is an isomorphism.
Hence $\Pic^n_{\bar{C},Z}$ is also a connected semi-abelian variety with unit $\theta(Q)$. In particular, $\Pic^n_{\bar{C},Z}(\BC)$ is a connected finite dimensional abelian Lie group.

For (i), since $\theta$ is an affine bundle, $\theta_*\colon \pi_1^{\top}(\Sym^nC)\xrightarrow{\sim }\pi_1^\top (\Pic^n_{\bar{C},Z})$; the second group is a finitely generated free abelian group by Lemma \ref{exam-Hure-top-lie-gp} (i).

For (ii), 
the morphism $\theta$ induces a commutative diagram
$$\xymatrix{H^i(\pi_1^\top(\Pic^n_{\bar{C},Z} ),A) \ar[rr]^{\epsilon_\top^i}  \ar[d]^{H^i(\theta_{*} ,A)} &&
H^i_\top(\Pic^n_{\bar{C},Z},A) \ar[d]^{\theta^*} \\
 H^i(\pi_1^\top(\Sym^nC),A)\ar[rr]^{\epsilon_\top^i} && H^i_\top(\Sym^nC,A).
}$$
Since $\theta$ is an affine bundle, 
$\theta^*\colon  H^i_\top(\Pic^n_{\bar{C},Z},A)\to H^i_\top(\Sym^nC,A)$ is an isomorphism for any $i\geq 0$; and since $\theta_*\colon \pi_1^{\top}(\Sym^nC)\xrightarrow{\sim }\pi_1^\top (\Pic^n_{\bar{C},Z})$, the morphism $H^i(\theta_{*} ,A)$ is an isomorphism. Now the result follows from the case of connected finite dimensional abelian Lie groups, which is given by Lemma \ref{exam-Hure-top-lie-gp}~(ii). Actually, one can also show that the universal cover of $\Sym^n C$ is contractible, since it is an affine bundle over some complex vector space (the universal cover of $\Pic^n_{\bar{C},Z}$).

For (iii), since $\pi_1^{\top}(\Sym^nC)\iso \pi_1^\top(\Pic^n_{\bar{C},Z})$ is finitely generated and free abelian, the image of $l_{P,*}$ is a normal subgroup. 
Suppose that $l_{P,*}$ is not surjective. Then one could find an abelian subgroup $H$ of $\pi_1^\top(\Pic^n_{\bar{C},Z})$ with finite index such that $\Im(l_{P,*})\subset H\subsetneq \pi_1^\top(\Pic^n_{\bar{C},Z})$, and then the group $H$ gives a covering space $G\to  \Pic^n_{\bar{C},Z}$  of connected abelian Lie groups with $\pi_1^\top(G)\cong H$ such that the morphism $C\xrightarrow{l_P} \Sym^nC \to \Pic^n_{\bar{C},Z} $ factors through $G$, i.e.,  we have a morphism $\psi\colon  C\to G$. Moreover, the morphism $\psi$ induces a morphism of complex manifolds 
$$\psi^n\colon  C^n\to G \colon   (c_1,\cdots,c_n)\mapsto \sum_i \psi(c_i).$$
Since $\psi^n$ is symmetric on $C$, it induces $\Sym^nC \to G$. That means that the morphism $\theta$ factors through $G$,  this contradicts to the surjectivity of $\pi_1(\Sym^nC)\iso \pi_1(\Pic^n_{\bar{C},Z})$.

It remains to show (iv).
Suppose that $\pi_1^\top(\Sym^nC)\cong \BZ^{\oplus r}$ for some $r\geq 0$ by (i). Since the morphism $l_{P,*}\colon \pi_1^\top(C)\to \pi_1^\top(\Sym^nC) $ is surjective , the Hurewicz theorem \cite[Thm. 2A.1]{Hatcher} implies that $$r=\rank(H_1^\top(\Sym^nC,\BZ))\leq \rank(H_1^\top(C,\BZ))\leq n$$ as $n>N_C$ and $N_C$ are defined as (\ref{eq-def-Nc}). Indeed, for affine curve $C$, one could show that 
$$\rank(H_1^\top(C,\BZ))=2\cdot g_{\Bar{C}}+\mathrm{deg}(\Bar{C}\setminus C)-1.$$
Let $e_1,\cdots,e_r$ be a basis of $\BZ^{\oplus r}\cong \pi_1^\top(\Sym^nC)$. One could lift each $e_j$ to $\alpha_j\in \pi_1^\top(C)$ via the surjective morphism $l_{P,*}\colon \pi_1^\top(C)\to \pi_1^\top(\Sym^nC) $. Then for every $1\leq j\leq r$, under the natural isomorphism $ \pi_1^\top(C^n)\cong \pi_1^\top(C)\times \cdots \times \pi_1^\top(C)$, the element $(1,\cdots,\alpha_j,\cdots,1)$ maps to $e_j$ under the morphism $p_*$, where $\alpha_j$ occurs in the $j$-th component. Indeed, the morphism $l_P$ is equal to the morphism $p\circ (P\times \cdots\times \mathrm{id}_C\times\cdots\times P)$ and the morphism $P\times \cdots\times \mathrm{id}_C\times P$ maps $\alpha_j$ to $(1,\cdots,\alpha_j,\cdots,1)$, where $\mathrm{id}_C$ and $\alpha_j$ both occurs in the $j$-th component. Define the homomorphism of groups
\begin{equation*}
    \begin{split}
      \varphi\colon \BZ^{\oplus r}&\to \pi_1^\top(C^n)\cong \pi_1^\top(C)\times \cdots \times \pi_1^\top(C) \qquad \text{($n$-times)}  \\
      e_j&\mapsto (1,\cdots, \alpha_j,\cdots,1),
    \end{split}
\end{equation*}
where $\alpha_j$ occurs in the $j$-th component. Then the group homomorphism
\begin{equation*}
    p_*\circ \varphi \colon \BZ^{\oplus r}\to \pi_1^\top(C^n)\to \pi_1^\top(\Sym^nC)\cong \BZ^{\oplus r}
\end{equation*}
is exactly the identity. In other words, morphism $p_*\colon \pi_1^\top(C^n)\to \pi_1^\top(\Sym^nC)$ is split surjective.

 By Riemann mapping theorem \cite[\S 27]{For81}, since $C$ is affine, the universal cover of $C$ is  contractible, so is the universal cover of $C^n$. Thus the edge morphism of Cartan--Leray spectral sequence (\ref{CLss for coho}) gives an isomorphism $H^i_\top(C^n,A)\cong  H^i(\pi_1^\top(C^n),A)$ for all $i\geq 0$. By the naturality of Cartan--Leray spectral sequences, it suffices to show that the morphism $${p}^*\colon H^i(\pi_1^\top(\Sym^nC),A)\to H^i(\pi_1^\top(C^n),A)$$ is split injective, since we also have $  H^i_\top(\Sym^nC,A)\cong H^i(\pi_1^\top(\Sym^nC),A)$. 
 It follows from the fact that $p_*\circ \varphi=\mathrm{id}\colon \pi_1^\top(\Sym^nC)\to \pi_1^\top(\Sym^nC)$. The proof is complete. 
\end{proof}
In fact, for $i=0,1,2$, the similar results of Proposition \ref{lem-DT-1} already hold from $n=2$, by applying the Dold--Thom theorem (see Theorem \ref{thm- toplogical DT}) and Proposition \ref{lem-DT-1}.

\begin{prop}\label{prop-n>1}
Let $C$ be an affine smooth integral complex curve. For every abelian group $A$ and any $n\geq 2$, we have
\begin{itemize}
    \item[(i)]  the topological fundamental group $\pi_1^\top(\Sym^n C)$ is a finitely generated free abelian group;
    \item[(ii)] the Hurewicz map $\epsilon_\top^i\colon  H^i(\pi_1^\top(\Sym^nC),A)\to H^i_\top(\Sym^nC,A)$ is surjective for all $i\geq 0$ 
    and is an isomorphism for $0\leq i\leq 2$;
    \item[(iii)]  for any fixed closed point $P$ of $C$, the morphism $$l_P\colon   C\xrightarrow{id_C\times P\times \cdots \times P} C^n\xrightarrow{p} \Sym^nC $$ induces a surjective homomorphism $l_{P,*}\colon \pi_1^\top(C)\to \pi_1^\top(\Sym^nC) $. In particular, the canonical morphism $p_*\colon \pi_1^\top(C^n)\to \pi_1^\top(\Sym^nC)$ is surjective;
    \item[(iv)] the natural morphism $p^*\colon H^i_\top(\Sym^nC,A)\to H^i_\top (C^n,A)$ is a split injection.
\end{itemize}
\end{prop}

\begin{proof}
The (i) follows from Theorem \ref{thm- toplogical DT} (iii) and Proposition \ref{lem-DT-1}. For any $n\geq 2$, there exists an integer $N\geq n$ with $N>N_C$. For (ii), consider the following commutative diagram induced by the inclusion $\Sym^nC\subset \Sym^N C$
\begin{equation*}
    \xymatrix{H^i(\pi_1^\top(\Sym^NC),A)\ar[rr]^{\cong}\ar[d]^{\cong}&& H^i(\Sym^NC,A)\ar@{->>}[d]\\
    H^i(\pi_1^\top(\Sym^nC),A)\ar[rr]^{\epsilon_\top^i}& &H^i(\Sym^nC,A),
    }
\end{equation*}
where 
\begin{itemize}
\item the top and the left are isomorphic by Proposition \ref{lem-DT-1} (ii) and  Theorem \ref{thm- toplogical DT} (c), respectively;
\item if $i\leq 2$, the bottom arrow is injective by Cartan--Leray spectral sequence (\ref{CLss for coho});
\item the right arrow is surjective by Corollary \ref{cor-topo DT}.
\end{itemize}
This implies (ii).

Consider the following commutative diagram
\begin{equation*}
\xymatrix{\pi_1^\top(C)\ar@/^1pc/[rrrr]^{l'_{P,*}}\ar[rr]\ar@{=}[d]&&\pi_1^\top(C^N)\ar[rr]_{p'_*}&&\pi_1^\top (\Sym^NC)\\
\pi_1^\top(C)\ar@/_1pc/[rrrr]_{l_{P,*}}\ar[rr]&&\pi_1^\top(C^n)\ar[u]\ar[rr]^{p_*}&&\pi_1^\top (\Sym^nC),\ar[u]^{\cong}
}
\end{equation*}
where $l_P\colon   C\xrightarrow{id_C\times P\times \cdots \times P} C^n\xrightarrow{p} \Sym^nC $ and $l'_P\colon   C\xrightarrow{id_C\times P\times \cdots \times P} C^N\xrightarrow{p'} \Sym^NC$. Then (iii) follows from the surjectivity of $l'_{P,*}\colon \pi_1^\top(C)\to \pi_1^\top(\Sym^NC)$; see Proposition \ref{lem-DT-1} (iii).

It remains to show (iv). It is clear that for $i=0,1$, the morphism $p^*$ is a split injection by universal coefficient theorem \cite[Thm. 3.3.2]{Hatcher}, Theorem \ref{thm- toplogical DT} (c) and K\"unneth formula (\ref{equ-top Kunn for H^1}). For $i=2$, since $H_j^\top(\Sym^nC,\BZ)$ and $H_j^\top(C,\BZ)$ are free abelian for every $j\geq 0$ (see Lemma \ref{exam-Hure-top-lie-gp}), applying universal coefficient theorem \cite[Thm. 3.3.2]{Hatcher} and K\"unneth formula \cite[Thm. 3B.6]{Hatcher} again, we have that 
\begin{equation*}
    H^2_\top(X,A)\cong \Hom_{\BZ}(H_2^\top(X,\BZ),A)
\end{equation*}
for $X=\Sym^nC$ or $C^n$. Therefore, it suffices to show that the morphism $p_*\colon H_2^\top(C^n,\BZ)\to H_2^\top(\Sym^n C,\BZ)$
is surjective and hence it is split as the group $H_2^\top(\Sym^n C,\BZ)$ is a finitely generated free abelian group by Lemma \ref{exam-Hure-top-lie-gp}. Similar to the proof of Proposition \ref{lem-DT-1} (iv), it suffices to show the surjectivity of the morphism
\begin{equation*}
    p_*\colon H_2(\pi_1^\top(C^n),\BZ)\to H_2(\pi_1^\top(\Sym^nC),\BZ).
\end{equation*}
Note that $H_2(\pi_1^\top(\Sym^nC),\BZ)\cong \wedge^2\pi_1^\top(\Sym^nC)$ since $\pi_1^\top(\Sym^nC)$ is abelian; see \cite[Thm. V.6.4]{brown2012cohomology}. If $\alpha\wedge\beta\in \wedge^2\pi_1^\top(\Sym^nC)$ with $\alpha,\beta\in \pi_1^\top(\Sym^nC)\cong H_1(C,\BZ)\cong \pi_1^\top(C)^{ab}$, the $\alpha\wedge\beta$ lies in the image of $H_2(\BZ^{\oplus2},\BZ)\to H_2(\pi_1^\top(\Sym^nC),\BZ)$ via the morphism
\begin{equation}\label{equ-Z2toSymnC}
        \BZ^{\oplus 2}\to \pi_1^\top(\Sym^nC)\cong \pi_1^\top(C)^{ab}\quad (a,b)\mapsto a\cdot\alpha+b\cdot \beta,
\end{equation}
because $H_2(\BZ^{\oplus2},\BZ)\cong \wedge^2(\BZ^{\oplus2})\cong \BZ$ with a generator $(1,0)\wedge (0,1)$.  Lift $\alpha$ and $\beta$ to $\pi_1^\top(C)^{ab}$, say $\alpha'$ and $\beta'$, respectively. Following the same arguments of the proof of Proposition \ref{lem-DT-1} (iv), the morphism (\ref{equ-Z2toSymnC}) factors through the morphism (note that $n\geq 2$)
\begin{equation*}
    \varphi\colon\BZ^{\oplus 2}\to \pi_1^\top(C^n)\cong \pi_1^\top(C)\times \cdots\times \pi_1^\top(C) \quad \text{($n$-times) } \quad (a,b)\mapsto (a\cdot\alpha',b\cdot \beta,0,\cdots,0)
\end{equation*}
via the natural morphism $p_*\colon \pi_1^\top(C^n)\to \pi_1^\top(\Sym^nC)$. Therefore, under the morphism $$p_*\colon H_2(\pi_1^\top(C^n),\BZ)\to H_2(\pi_1^\top(\Sym^nC),\BZ),$$ the element $\varphi_*(1,0)\wedge (0,1))\in H_2(\pi_1^\top(C^n),\BZ)$ maps to $\alpha\wedge\beta$. We finish our proof by noting that such elements $\alpha\wedge\beta$ generate the group $H_2(\pi_1^\top(\Sym^nC),\BZ)$. The proof is complete.
\end{proof}

Similar to topological cases, one could get the following corollary under the \'etale topology.
\begin{prop}\label{prop-curves}
Let $k$ be an algebraically closed field of characteristic zero. Let $C$ be an affine smooth integral curve over $k$. Then for any torsion sheaf $F$ over $k_{\et}$, we have
\begin{itemize}
    \item[(i)]  the \'etale fundamental group $\pi_1(\Sym^nC)$ is commutative for all $n\geq 2$;
    \item[(ii)] the Hurewicz map $\epsilon^i\colon  H^i_{cont}(\pi_1(\Sym^nC),F)\to H^i(\Sym^nC,F)$ is an isomorphism if either $n>N_C$ and $i\geq 0$, or $n\geq 2$ and $i\leq2$;
    \item[(iii)]  for any fixed closed point $P\in{C}(k)$, the morphism \begin{equation*}
      l_P\colon  {C}\xrightarrow{\mathrm{id}\times P\times \cdots\times P}{C}^n\xrightarrow{\rm proj.} \Sym^n C
\end{equation*}
induces a surjective homomorphism $l_{P,*}\colon \pi_1({C})\to \pi_1(\Sym^nC) $ for all $n\geq 2$;
    \item[(iv)]  the canonical morphism $p_*\colon \pi_1^\top(C^n)\to \pi_1^\top(\Sym^nC)$ is surjective for $n\geq 2$ which is split for $n>N_C$. The natural morphism $p^*\colon H^i(\Sym^nC,F)\to H^i ({C}^n,F)$ is split injective if either $n>N_C$ and $i\geq 0$, or $n\geq 2$ and $i\leq2$.
\end{itemize}
\end{prop}

\begin{proof}
By passing to the limit, there exists a subfield $\BQ\subset L\subset k$ and an affine  smooth integral curve $C'$ over $L$ such that the field $L$ is algebraically closed with an inclusion  $L\subset \BC$ and $C'\times_Lk\cong C$. Also, we may assume that the closed point $P$ is defined over the field $L$.

By the rigidity of fundamental groups and the comparison theorem (\ref{eq-fund-gps-5}), we have isomorphisms
\begin{equation*}
  \pi_1(\Sym^nC)\cong \pi_1(\Sym^nC') \cong \pi_1(\Sym^nC'_{\BC}) \cong \widehat{ \pi_1^{\top}(\Sym^nC'(\BC))},  \qquad  
\end{equation*}
and 
\begin{equation*}
   \pi_1({C}^n)\cong \pi_1({C}'^n) \cong \pi_1({C}'^n_{\BC}) \cong \widehat{ \pi_1^{\top}({C}'^n(\BC))}  .
\end{equation*}
These isomorphisms are compatible with the natural morphisms $l_{P,*}$ and $p_*$. Therefore, (i) and (iii) follow from Proposition \ref{prop-n>1} (i) and (iii), respectively; since the profinite completion of an abelian group is abelian, and profinite completion preserves surjections of arbitrary groups.

By the rigidity (\cite[Cor. VI.4.3]{Mi80}) and Artin comparison theorem (\cite[Cor. III.3.12]{Mi80}) of \'etale cohomology, we have natural isomorphisms
\begin{equation}
     H^i(\Sym^nC,F )\cong H^i(\Sym^nC',F) \cong H^i(\Sym^nC'_{\BC},F)\cong H^i_{\top}(\Sym^nC(\BC),F),
\end{equation}
and these isomorphism also hold after replacing $\Sym^nC$ and $\Sym^nC'$ with ${C}^n$ and ${C}'^n$, respectively.
Thus, (iv) is given by Proposition \ref{lem-DT-1} (iv) and Proposition \ref{prop-n>1} (iv); (ii) follows from Proposition \ref{prop-compatible of edge morphisms}, Proposition \ref{lem-DT-1} (ii) and Proposition \ref{prop-n>1} (ii). The proof is complete.
\end{proof}

\section{Dold--Thom theorem for Brauer groups}\label{section4}
Let $k$ be a field of characteristic zero and $F$ be a torsion sheaf over $k_{\et}$.
Let $X$ be a quasi-projective smooth geometrically integral $k$-variety.

By  Proposition \ref{prop-Lef}, there exists a smooth geometrically integral sub-curve $C\subset X$ (not necessary closed) 
such that $\pi_1(C)\to \pi_1(X)$ is surjective.
We may assume $C$ is affine, 
because, if $C$ is projective, we may replace $C$ by its non-trivial open subset. Actually, by Lemma \ref{lem-DT-1-2} below, 
the surjectivity of $\pi_1(C)\to \pi_1(X)$ is equivalent to the surjectivity of $\pi_1(C_{\bk})\to \pi_1(X_{\bk})$.

Let $\Sym^nX$ be the symmetric power of $X$ and let $\tSym^nX$ be a resolution of singularity of $\Sym^nX$.
We recall the notation $U^n(X),U(X),V(X),\tilde{V}(X) $ as in (\ref{notation-sym-eq}). 
Let $\pr_i\colon X^n\to X$  be the $i$-th projection.

\subsection{Dold--Thom of degree 1}

Since the morphism $\tilde{V}(X)\to V(X)$ is birational proper, 
the projection morphism $ V^n(X)\to V(X)$ can be lifted to a rational map $p\colon  X^n \dashrightarrow \tilde{V}(X)$,
which is well defined on all codimension one points of $X^n$.
Hence there exists an open subset $U_X\subset X^n$ such that $p$ is well defined on $U_X$ and $\codim(X^n\setminus U_X,X^n)\geq 2$. 
By the purity of \'etale fundamental groups (see (\ref{eq-fund-gps-1})) and of \'etale cohomology (see \cite[Cor. VI.5.3]{Mi80}), 
we have $\pi_1(U_X)\cong \pi_1(X^n) $ and $H^i(U_X,F)\cong H^i(X^n,F)$ for $i=1,2$,
and hence the push-forward and the pullback
\begin{equation}\label{eq-def-p}
p_*\colon \pi_1(X^n)\cong \pi_1(U_X)\to \pi_1(\tilde{V}(X)) \ \ \ \text{and}\ \ \   p^*\colon  H^i(\tilde{V}(X),F)\to H^i(U_X,F)\cong H^i(X^n,F)
\end{equation}
are well defined for $i=1,2$.

\begin{thm}\label{thm-DT-1}
Let $k$ be algebraically closed of characteristic zero.
Let $X$ be a  quasi-projective smooth geometrically integral $k$-variety, 
and let $F$ be a torsion sheaf over $k_{\et}$.
For any $n\geq 2$, we have
\begin{itemize}
    \item[(i)] the fundamental groups $\pi_1(\tilde{V}(X))\cong \pi_1(\tSym^nX)\cong \pi_1(X)^{ab}$ are commutative;
    \item[(ii)] the $p^*$ in (\ref{eq-def-p}) induces isomorphisms 
    $H^1(\tSym^nX,F)\iso H^1(\tilde{V}(X),F)\iso   H^1(X^n,F)^{S_n} $;
    \item[(iii)] the canonical transfer $\tr_X^1\colon  H^1(X,F)\to H^1(\tSym^nX,F)$ in (\ref{eq-def-transfer}) is an isomorphism;
    \item[(iv)] there exists an isomorphism $\pi_1(X)^{ab} \iso \pi_1(\tilde{V}(X))$ such that
   the $p_*$ in (\ref{eq-def-p}) is equal to the composition
    $$\pi_1(X^n)\to \pi_1(X^n)^{ab}\xrightarrow{\sum_i \pr_{i,*}}  \pi_1(X)^{ab} \iso \pi_1(\tilde{V}(X))  .$$
\end{itemize} 
\end{thm}

Before proving Theorem \ref{thm-DT-1}, 
we consider an immersion $f\colon Y\to X$ and prove some reduction results in Proposition \ref{lem-DT-1-4}.

\begin{lem}\label{lem-DT-1-3}
Let $X$, $Y$ be two quasi-projective smooth geometrically integral $k$-varieties of dimension~$>0$ and let $f\colon Y\to X$ be an immersion.
Then $f$ induces a canonical morphism $\tilde{f}_V\colon  \tilde{V}(Y)\to \tilde{V}(X) $ and the commutative diagram
\begin{equation}\label{lem-DT-1-3-eq}
\xymatrix{Y^n \ar[d]^{f^n}  &  U^n(Y)\ar[d]^{f_U^n} \ar[r]^{\psi_Y}\ar@{_(->}[l]_{i_Y} &
U(Y)\ar[d]^{f_U} \ar@{^(->}[r]^{j_Y}  & \tilde{V}(Y) \ar[d]^{\tilde{f}_V}  \\
X^n  &  U^n(X)\ar[r]^{\psi_X} \ar@{_(->}[l]_{i_X}  &  U(X) \ar@{^(->}[r]^{j_X}  & \tilde{V}(X)  .
}
\end{equation}
\end{lem}

\begin{proof}
Since $f$ is an immersion, it induces morphisms $f^n$, $f_U^n$, $f_U$ and $f_V^n\colon V^n(Y)\to V^n(X)$.
By \cite[Cor. II.7.15]{Hartshorne}, $f^n_V$ induces a canonical morphism 
$\Bl_{V\setminus U}V^n(Y)\to \Bl_{V\setminus U}V^n(X) $ (see the notation above (\ref{notation-sym-eq})),
and hence a morphism
$$\tilde{f}_V\colon  \tilde{V}(Y)\coloneqq  \Bl_{V\setminus U}V^n(Y)/S_n \to  \Bl_{V\setminus U}V^n(X)/S_n \eqqcolon \tilde{V}(X) .$$
In addition, the commutative diagram can be checked directly.
\end{proof}

\begin{lem}\label{lem-DT-1-2}
Let $X,Y$ be smooth geometrically integral $k$-varieties and let $f\colon Y\to X$ be a morphism.
Then the morphism $f_*\colon \pi_1(Y)\to \pi_1(X)$ is surjective (\emph{resp.} injective, isomorphism) 
if and only if the morphism  $f_{\bk,*}\colon \pi_1(Y_{\bk})\to \pi_1(X_{\bk})$ is surjective (\emph{resp.} injective, isomorphism).

Moreover, in this case, the morphisms $f^n_*\colon\pi_1(Y^n)\to \pi_1(X^n)$ and $f^n_{\bk,*}\colon\pi_1(Y^n_{\bk})\to \pi_1(X^n_{\bk})$ are surjective (\emph{resp.} injective, isomorphism).
\end{lem}

\begin{proof}
The (\ref{eq-fund-gps-4}) induces a commutative diagram
$$\xymatrix{1\ar[r] & \pi_1(Y_{\bk})\ar[r]\ar[d]^{f_{\bk,*}}  &\pi_1(Y)\ar[r]\ar[d]^{f_*} &\Gal(\bk/k)\ar[d]^=\ar[r] &1\\
1\ar[r] & \pi_1(X_{\bk})\ar[r]  &\pi_1(X)\ar[r] &\Gal(\bk/k)\ar[r] &1.
}$$
Then the first statement follows from a diagram chasing.

For the ``moreover'' part, since $f_{\bk,*}$ is surjective  (injective, isomorphism),
(\ref{eq-fund-gps-2}) implies that $f^n_{\bk,*}$ satisfies the same property,
and the first statement implies that $f^n_*$ also satisfies the same property.
\end{proof}

\begin{prop}\label{lem-DT-1-4}
Under the hypothesis of Lemma \ref{lem-DT-1-3}, assume that $k$ is algebraically closed.
If the morphism $f_*\colon\pi_1(Y)\to \pi_1(X)$ is surjective, then the morphism $\tilde{f}_{V,*}\colon \pi_1(\tilde{V}(Y))\to \pi_1(\tilde{V}(X))$ is surjective.
\end{prop}

\begin{proof}
If $\dim(X)=1$, then $\dim(Y)=1$ and $f$ is an open immersion.
In this case, $\tilde{V}(X)=V(X)$, $\tilde{V}(Y)=V(Y)$ and $\tilde{f}_V$ also is an open immersion.
Hence $\tilde{f}_{V,*}$ is surjective by (\ref{eq-fund-gps-1}).

We assume $\dim(X)\geq 2$.
Applying $\pi_1(-)$ to the diagram (\ref{lem-DT-1-3-eq}), we get the commutative diagram
\begin{equation}\label{eq-inclu-CX-dia1}
\xymatrix{\pi_1(Y^n) \ar[d]^{f^n_{*}} &  \pi_1(U^n(Y))\ar[d]^{f^n_{U,*}} \ar[r]^{\psi_{Y,*}}\ar@{->>}[l]_{i_{Y,*}} &
\pi_1(U(Y))\ar[d]^{f_{U,*}} \ar@{->>}[r]^{j_{Y,*}}  & \pi_1(\tilde{V}(Y)) \ar[d]^{\tilde{f}_{V,*}} \\
\pi_1(X^n)  &  \pi_1(U^n(X))\ar[r]^{\psi_{X,*}} \ar@{->>}[l]_{i_{X,*}}  & \pi_1(U(X)) \ar@{->>}[r]^{j_{X,*}}  
& \pi_1(\tilde{V}(X)) ,
}
\end{equation}
where $i_{Y,*},i_{X,*}, j_{Y,*},j_{X,*}$ are surjective by (\ref{eq-fund-gps-1}).
Since $\dim(X)\geq 2$, we have $\codim(X^n\setminus U^n(X),X^n)\geq 2 $ and $i_{X,*}$ is an isomorphism
by the purity of \'etale fundamental groups.
By Lemma \ref{lem-DT-1-2}, $f^n_*$ is surjective.
Hence $f^n_{U,*}$ is surjective.
Since $\psi_Y,\psi_X$ are $S_n$-torsor, we have commutative diagram with exact rows
$$\xymatrix{1\ar[r] & \pi_1(U^n(Y))\ar[d]^{f^n_{U,*}}\ar[r]  & \pi_1(U(Y))\ar[r]\ar[d]^{f_{U,*}} &S_n\ar[r]\ar[d]^= &1\\
1\ar[r] & \pi_1( U^n(X))\ar[r]  & \pi_1(U(X)) \ar[r] &S_n\ar[r] &1.
}$$
Hence $f_{U,*}$ is surjective and therefore $\tilde{f}_{V,*}$ is also surjective.
\end{proof}

\begin{proof}[Proof of Theorem \ref{thm-DT-1}]
Let $f\colon C\to X$ be an affine sub-curve such that the morphism $\pi_1(C)\to \pi_1(X)$ is surjective. We denote by $\tilde{f}_V\colon \tilde{V}(C)\to \tilde{V}(X)$ the morphism induced in Lemma \ref{lem-DT-1-3},
and denote by $l\colon  \tilde{V}(X)\hookrightarrow \tSym^nX$ the open immersion.
By Proposition \ref{lem-DT-1-4}, all morphisms in
$$ \pi_1(\tilde{V}(C))\xrightarrow{\tilde{f}_{V,*}} \pi_1(\tilde{V}(X))\xrightarrow{l_*}  \pi_1(\tSym^nX )$$
are surjective.

For (i), by Proposition \ref{prop-curves} (i), $\pi_1(\Sym^nC)$ is commutative.
Since $C$ is a curve, we have $\tilde{V}(C)=V(C)$ and $\codim(\Sym^nC\setminus V(C),V(C))\geq 2$.
Hence $\pi_1(\tilde{V}(C))\cong \pi_1(\Sym^nC)$ are commutative, by the purity of \'etale fundamental groups.
Since $\tilde{f}_{V,*},$ and $l_*$ are surjective, 
the groups $\pi_1(\tilde{V}(X))$ and $\pi_1(\tSym^nX)$ are also commutative. It remains to show that $l_*$ is an isomorphism.

For (ii) and (iii), by (\ref{eq-fund-gps-1}) and Proposition \ref{prop-curves} (iv), $\pi_1(U^n(C))\to \pi_1(C^n)\to \pi_1(\Sym^nC)\cong \pi_1(\tilde{V}(C))$ is surjective.
Applying the diagram (\ref{eq-inclu-CX-dia1}) to the case $Y=C$ and using the surjectivity of $\tilde{f}_{V,*}$, 
the composition $j_{X,*}\circ \psi_{X,*} $ is surjective.
By (\ref{eq-fund-gps-3}), $\psi_X^*\circ j_X^*\colon  H^1(\tilde{V}(X),F)\to H^1(U^n(X),F) $ is injective,
and hence $p^*$ is injective.
Then we have inclusions
$$ H^1(\tSym^nX  ,F)\xhookrightarrow{l^*}   H^1(\tilde{V}(X),F) \xhookrightarrow{p^*}  H^1(X^n,F) ^{S_n}.$$

Since $H^1(X^n,F)\hookrightarrow H^1(U^n(X),F)$ is injective, by Proposition \ref{eq-def-transfer1}, we have
\begin{equation}\label{Thm-DT-1-eq-1}
 p^*\circ l^*\circ  \tr_X^1=\sum_i\pr_i^* \colon  H^1(X,F)\to H^1(X^n,F). 
 \end{equation}
By the K\"unneth formula (\ref{Thm-Kunneth-eq3}), 
$\oplus \pr_i^*\colon \oplus_{i=1}^n H^1(X,F)\to H^1(X^n,F) $ is an isomorphism.
Hence the morphism $\sum_i\pr_i^*$ is injective and we have
$$H^1(X^n,F) ^{S_n}= \Im(\sum_i\pr_i^* ) =\Im (p^*\circ l^*\circ  \tr_X^1).$$
This implies (ii) and (iii) and also $l^*$ is an isomorphism.

 For any smooth integral $k$-variety $Y$, the Pontryagin dual of (\ref{eq-fund-gps-3}) provides
 \begin{equation}\label{equ-Pontryagin dual}
     \pi_1(Y)^{ab}\cong  \Hom_{cont}(H^1(Y,\BQ/\BZ),\BQ/\BZ).
 \end{equation}
Apply (iii) to the sheaf $\BQ/\BZ$, since $\pi_1(\tilde{V}(X))$ is commutative,
then the Pontryagin dual of $l^*\circ tr^1_X$ induces an isomorphism 
$\pi_1(\tilde{V}(X))\cong \pi_1(\tilde{V}(X))^{ab}\iso \pi_1(X)^{ab}$.
Since $p_*$ is the Pontryagin dual of $p^*$, 
the Pontryagin dual of (\ref{Thm-DT-1-eq-1}) implies (iv).

Finally, applying (\ref{equ-Pontryagin dual}) and (ii), we have $\pi_1(\tilde{V}(X))\cong \pi_1(\tSym^nX)$. The proof is complete.
\end{proof}

\begin{rem}\label{rem-DT1-top}
   One could also get that $\pi_1(\tSym^nX)\cong\pi_1(\Sym^nX)\cong \pi_1(X)^{ab}$ directly, after applying \cite[Thm. (7.5.2)]{Kollar1993Shafarevich} and the rigidity theorem of fundamental groups (\ref{eq-fund-gps-5}). 
   Moreover, if $k=\BC$ and we take topological fundamental groups and singular cohomology groups, one could get the topological version of Theorem \ref{thm-DT-1}. Because we still could find an affine sub-curve $C\subset X$ such that the morphism $\pi_1^\top(C)\to \pi_1^\top(X)$ is surjective by Proposition \ref{prop-Lef}.
\end{rem}

\begin{cor}\label{Thm-DT-1-cor1}
Let $X$ be a  quasi-projective smooth geometrically integral $k$-variety, 
and $F$ a torsion sheaf over $k_{\et}$.
Then, for any $n\geq 2$, the morphism $p^*$ in (\ref{eq-def-p}) induces an isomorphism
$$ p^*\colon  H^1(\tSym^nX,F)\cong H^1(X^n,F)^{S_n}.$$
\end{cor}

\begin{proof}
    Let us consider the following commutative diagram
    \begin{equation*}
        \xymatrix{0\ar[r] &  H^1(k,F)\ar[r]\ar@{=}[d]& H^1(\tSym^nX,F)\ar[r]\ar[d]^{p^*}  &
H^1(\tSym^nX_{\bk},F)^{\Gal(\bk/k)}\ar[r]\ar[d]^{p_{\bk}^*}_{\sim}  & H^2(k,F)\ar@{=}[d]    \\
0\ar[r] &  H^1(k,F)\ar[r]\ar@{=}[d] & H^1(X^n,F)^{S_n}\ar[r]^-{res}\ar@{^(->}[d]  &  (H^1(X^n_{\bk},F)^{S_n})^{\Gal(\bk/k)} \ar[r]\ar@{^(->}[d] &   H^2(k,F)\ar@{=}[d] \\
        0\ar[r] &  H^1(k,F)\ar[r] & H^1(X^n,F)\ar[r]^-{res}  &  H^1(X^n_{\bk},F)^{\Gal(\bk/k)} \ar[r] &   H^2(k,F),
        }
    \end{equation*}
    where the first and the third rows are exact by the Hochschild-Serre spectral sequence (\cite[Thm. III.2.20]{Mi80}).
    The middle row is exact at $H^1(k,F) $ and $ H^1(X^n,F)^{S_n}$ by the fact that taking $S_n$-invariant is left exact.
    From Theorem \ref{thm-DT-1} (ii), the morphism $p_{\bk}^*$ is an isomorphism.
    The diagram chasing shows that the middle row is exact.
    By the five lemma, the morphism $p^*$ is an isomorphism.
\end{proof}

\subsection{Dold--Thom of degree 2}

In this subsection, we calculate $H^2(\tilde{V}(X),F)$ over algebraically closed fields of characteristic zero
by using the Hurewicz map ${\epsilon}_{\tilde{V}(X)}$ in (\ref{eq-def HS ss for F}) 
and the canonical transfer $\tr_X^2$ in (\ref{eq-def-transfer}); see Theorem \ref{thm-DT-2}.

If there is no ambiguity, we also denote by $\tr_X^2$ for the composition of original $\tr_X^2$ with the natural restriction 
$ H^2(\tSym^nX,F) \to H^2(\tilde{V}(X),F)$. Let $D_X\coloneqq \tilde{V}(X)\setminus U(X)$. 
By Proposition \ref{notation-sym-prop}, $D_X$ is smooth geometrically integral.
The Gysin sequence (cf.\cite[Cor. VI.5.3]{Mi80}) induces the following exact sequences for any admissible $(k,F)$
\begin{equation*}
    F(-1)\xrightarrow{i_{D_X}}  H^2(\tilde{V}(X),F)\to H^2(U(X),F) 
\end{equation*}
and 
\begin{equation*}
    F(-1)^{\oplus m}\xrightarrow{i_*} H^2(\tSym^nX,F)\to H^2(U(X),F),
\end{equation*}
where $m$ is the number of the irreducible complements of $\tSym^n X\setminus U(X)$ with codimension $1$.

\begin{thm}\label{thm-DT-2}
Let $k$ be an algebraically closed field,
and let $X$ be a quasi-projective smooth integral $k$-variety.
Then, for any $n\geq 2$ and any torsion \'etale sheaf $F$, the following homomorphisms are surjective
\begin{equation}\label{equ-thm4.61}
    (\epsilon_{\tSym^nX},\tr_X^2,i_*)\colon   H^2_{cont}(\pi_1(\tSym^nX),F)\oplus H^2(X,F) \oplus 
F(-1)^{\oplus m}\to H^2(\tSym^nX,F)
\end{equation}
and 
\begin{equation}\label{equ-thm4.62}
   (\epsilon_{\tilde{V}(X)},\tr_X^2,i_{D_X})\colon   H^2_{cont}(\pi_1(\tilde{V}(X)),F)\oplus H^2(X,F) \oplus 
F(-1)\to H^2(\tilde{V}(X),F). 
\end{equation}
 
\end{thm}
\begin{rem}
Suppose $k=\BC$ and we take singular cohomology. In fact, we first show that for any $n\geq 2$ and any abelian group $F$, the following homomorphisms are surjective
\begin{equation}\label{equ-thm4.63}
    (\epsilon^\top_{\tSym^nX},\tr_X^2,i_*)\colon   H^2(\pi^\top_1(\tSym^nX),F)\oplus H^2_\top(X,F) \oplus 
F(-1)^{\oplus m}\to H^2_\top(\tSym^nX,F)
\end{equation}
and 
\begin{equation}\label{equ-thm4.64}
   (\epsilon^\top_{\tilde{V}(X)},\tr_X^2,i_{D_X})\colon   H^2(\pi^\top_1(\tilde{V}(X)),F)\oplus H^2_\top(X,F) \oplus 
F(-1)\to H^2_\top(\tilde{V}(X),F). 
\end{equation}
Then we could prove Theorem \ref{thm-DT-2} after applying Artin's comparison theorem and Proposition \ref{prop-compatible of edge morphisms}.
\end{rem}

We only prove the surjectivity of (\ref{equ-thm4.64}) (and hence the surjectivity of (\ref{equ-thm4.62})), and the surjectivity of (\ref{equ-thm4.63}) and (\ref{equ-thm4.61}) follows from the surjectivity of (\ref{equ-thm4.62}) and Gysin exact sequence; we leave to the reader.

Recall that the K\"unneth formula (\ref{equ-top kun for H^2}) induces a canonical projection for any abelian group $F$
$$\lambda_X\coloneqq \bigoplus \lambda_{i< j}\colon  H^2_\top(X^n,F)\to 
\bigoplus\limits_{1\leq i<j\leq n}  \Hom_{\BZ} ( \pi_1^\top(X)^{ab}\otimes_\BZ \pi_1^\top(X)^{ab} ,F) . $$
We recall the $p^*\colon  H^2_\top(\tilde{V}(X) ,F)\to H^2_\top(X^n,F)$ in  (\ref{eq-def-p}). We first suppose $\dim(X)\geq 2$, the by purity of \'etale cohomology, this morphism is exactly the composition of 
\begin{equation*}
    p^*\colon H^2_\top(\tilde{V}(X),F)\xrightarrow{\rm restr.} H^2_\top(U(X),F)\xrightarrow{\rm pullback} H^2_\top(U^n(X),F)\xleftarrow{\sim} H^2_\top(X^n,F).
\end{equation*}
Hence, the morphism $p^*$ can be factored through as 
\begin{equation}
   H^2_\top(\tilde{V}(X),F)\to H^2_\top(U^n(X),F)^{S_n}\cong H^2_\top(X^n,F)^{S_n}\hookrightarrow H^2_\top(X^n,F).
\end{equation}
If $\dim(X)=1$, the morphism $p^*$ is the composition of the following, by purity,
\begin{equation*}
    p^*\colon H^2_\top(\tilde{V}(X),F)\xleftarrow{\sim} H^2_\top(\Sym^n X,F)\xrightarrow{\rm pullback} H^2_\top(X^n,F)^{S_n}\subset H^2_\top (X^n,F).
\end{equation*}
Also, the image of the morphism  $H^2_{\top}(\tSym^nX,F)\xrightarrow{\rm restr.}H^2_{\top}(\tilde{V}(X),F)\to H^2_\top(X^n,F)$ is contained in $H^2_\top(X^n,F)^{S_n}$. Therefore, for simplicity, we sometimes write $p^*\colon H^2_{\top}(\tSym^nX,F)\to H^2_\top(X^n,F)^{S_n}$ and $p^*\colon H^2_{\top}(\tilde{V}(X),F)\to H^2_\top(X^n,F)^{S_n}$ for the pullbacks.

Moreover, since $k$ is algebraically closed, we ignore the Tate twisted $(-1)$ and then we write $F(-1)$ as $F$ for simplicity. 

\begin{lem}\label{lem-DT-2-prop1}
Let $X$ be a quasi-projective smooth integral complex variety. For any $n\geq 2$, we have $\Im(\lambda_X\circ p^* )=\Im(\lambda_X\circ p^* \circ {\epsilon}_{\tilde{V}(X)})$. Here all groups related to these morphisms are singular cohomology with an abelian group $F$.
\end{lem}

\begin{proof}
Consider the following commutative diagram

\begin{equation*}
    \xymatrix{H^2(\pi_1^\top(\tilde{V}(X)))\ar@{->>}[r]\ar[d]^{p^*}\ar@/^1pc/[rrr]^{\epsilon_{\tilde{V}(X)}}&\Hom_\BZ(\wedge^2\pi_1^\top(\tilde{V}(X)),F)\ar[r]_{\epsilon^{\tilde{V}(X)}}\ar[d]^{(p_*)^\vee}&\Hom_\BZ(H_2^\top(\tilde{V}(X),\BZ),F)\ar[dd]^{(p_*)^\vee}& H^2_\top(\tilde{V}(X))\ar@{->>}[l]\ar[dd]^{p^*}\\
H^2(\pi_1^\top(X^n)^{ab})\ar@{->>}[r]\ar[d]&\Hom_\BZ(\wedge^2\pi_1(X^n)^{ab},F)\ar[d]&&\\
H^2(\pi_1^\top(X^n))\ar@{->>}[r]\ar@/_1pc/[rrr]_{\epsilon_{X^n}}&\Hom_\BZ(\wedge^2\pi_1^\top(X^n),F)\ar[r]^{\epsilon^{X^n}}&\Hom_\BZ(H_2^\top(X^n,\BZ),F)\ar@{->>}[d]^{(\lambda^X)^\vee}& H^2_\top(X^n)\ar@{->>}[l]\ar@{->>}[ld]^{\lambda_X}\\
&&\bigoplus\limits_{1\leq i<j\leq n}\Hom_\BZ((\pi_1(X)^{ab})^{\otimes 2},F)&,
        }
\end{equation*}
where we omit the coefficients $F$ for all cohomology groups and the first two vertical arrows $p^*$ and $(p_*)^\vee$ exist as $\pi_1^\top(\tilde{V}(X))$ is abelian (see Remark \ref{rem-DT1-top}). Therefore, we have
\begin{equation*}
   \Im(\lambda_X\circ p^* )=\Im((\lambda^X)^\vee\circ (p_*)^\vee)\quad \text{and}\quad  \Im(\lambda_X\circ p^* \circ {\epsilon}_{\tilde{V}(X)})=\Im((\lambda^X)^\vee\circ (p_*)^\vee \circ {\epsilon}^{\tilde{V}(X)});
\end{equation*}
and hence it suffices to show that we have
\begin{equation}\label{equ-prop 4.7(i)}
    \Im((\lambda^X)^\vee\circ (p_*)^\vee)=\Im((\lambda^X)^\vee\circ (p_*)^\vee \circ {\epsilon}^{\tilde{V}(X)}).
\end{equation}

We calculate the composition $(\lambda^X)^\vee\circ (p_*)^\vee \circ {\epsilon}^{\tilde{V}(X)}$ first. Note that we have the commutative diagram again
\begin{equation*}
 \xymatrix{\Hom_{\BZ}(\wedge^2 \pi_1^\top(\tilde{V}(X))^{ab},F)  \ar[d]^{\Hom_\BZ(\wedge^2p_*,F)} 
\ar[r]^-{{\epsilon}^{\tilde{V}(X)}} & \Hom_\BZ(H_2^\top(\tilde{V}(X),\BZ),F)  \ar[d]^{(p_*)^{\vee}} &   \\
\Hom_{\BZ}(\wedge^2 \pi_1^\top(X^n)^{ab}, F)  \ar[r]^-{{\epsilon}^{X^n}} & \Hom_\BZ(H_2^\top(X^n,\BZ),F)\ar[r]^-{(\lambda^X)^\vee}
&\bigoplus\limits_{1\leq i<j\leq n}  \Hom_{\BZ} ( \pi_1^\top(X)^{ab}\otimes \pi_1^\top(X)^{ab} ,F) .
}   
\end{equation*}

Note that 
$(\lambda^X)^\vee\circ {\epsilon}^{X^n}$ is the $\Hom_\BZ(-,F)$ of the homomorphism of fundamental groups (\ref{equ-a lot composition}).
By Theorem~\ref{thm-DT-1} (iv), 
the composition $(\lambda^X)^\vee\circ (p_*)^\vee \circ {\epsilon}^{\tilde{V}(X)}$  is the $\Hom_\BZ(-,F)$ of the surjection
\begin{equation}\label{equ-Hom from tensor to wedge}
    \begin{split}
        q_X\colon \bigoplus_{1\leq i<j\leq n}  \pi_1^\top(X)^{ab}\otimes \pi_1^\top(X)^{ab} &\to 
\wedge^2 \pi^\top_1(\tilde{V}(X))^{ab}\cong \wedge^2 \pi_1^\top(X)^{ab} \\
(a_{ij}\otimes b_{ij})_{1\leq i<j\leq n}&\mapsto  \sum_{1\leq i<j\leq n}(a_{ij}\wedge b_{ij}).
    \end{split}
\end{equation}

Let $C\subset X$ be an affine sub-curve such that the morphism $\pi_1^\top(C)\to \pi_1^\top(X)$ is surjective; see Proposition \ref{prop-Lef}. 
This induces a commutative diagram
\begin{equation*}
  \xymatrix{ \Hom_{\BZ}(\wedge^2 \pi_1^\top(\tilde{V}(X))^{ab} ,F)  \ar@{^(->}[d]^{\psi_1}  \ar[r]^-{{\epsilon}^{\tilde{V}(X)}} &
\Hom_\BZ(H_2^\top(\tilde{V}(X),\BZ),F)  \ar[rr]^-{(\lambda^X)^\vee\circ (p_{X,*})^\vee} \ar[d] & &  
\bigoplus\limits_{1\leq i<j\leq n}  \Hom_{\BZ} ( (\pi_1^\top(X)^{ab})^{\otimes 2} ,F) \ar@{^(->}[d]^{\psi_2} \\
 \Hom_{\BZ}(\wedge^2 \pi_1^\top(V(C))^{ab} ,F) \ar[r]^-{{\epsilon}^{V(C)}}_-{\sim} &
   \Hom_\BZ(H_2^\top({V}(C),\BZ),F)   \ar[rr]^-{(\lambda^C)^\vee\circ (p_{C,*})^\vee}  &&
   \bigoplus\limits_{1\leq i<j\leq n}  \Hom_{\BZ} ( (\pi_1^\top(C)^{ab})^{\otimes 2},F),
   }
\end{equation*}
where 
\begin{itemize}
    \item  $\psi_1$ and $\psi_2$ are injective since $\pi_1^\top(C)\to \pi_1^\top(X)$ is surjective;
    \item ${\epsilon}^{V(C)}$ is an isomorphism as $\codim(\Sym^nC\setminus V(C),\Sym^nC)=2$ and by Proposition \ref{prop-n>1} (ii).
\end{itemize}
Thus, $\Im((\lambda^C)^\vee\circ (p_{C,*})^\vee)=\Im((\lambda^C)^\vee\circ (p_{C,*})^\vee \circ {\epsilon}^{{V}(C)})$.

On the other hand, proving (\ref{equ-prop 4.7(i)}) is equivalent to show that for every $$\beta\in \Im((\lambda^X)^\vee\circ (p_*)^\vee )\subset \bigoplus_{1\leq i<j\leq n}  \Hom_{\BZ} ( \pi_1^\top(X)^{ab}\otimes \pi_1^\top(X)^{ab} ,F), $$ we have $\beta (\Ker(q_X))=0$; because we have the following exact sequence 
\begin{equation*}
    0\to \Hom_{\BZ}(\wedge^2 \pi_1^\top(\tilde{V}(X))^{ab},F) \xrightarrow{q_X^*} \bigoplus_{1\leq i<j\leq n}  \Hom_{\BZ} ( \pi_1^\top(X)^{ab}\otimes \pi_1^\top(X)^{ab} ,F)\to \Hom_{\BZ}(\Ker(q_X),F).
\end{equation*}
Note that $\psi_2(\beta)\in \Im((\lambda^X)^\vee\circ (p_*)^\vee \circ {\epsilon}^{{V}(C)})$, and so we have $\psi_2(\beta)(\Ker(q_C))=0$. Hence, it is enough to show that for any $a\in \Ker(q_X)$, there exists $b\in \Ker(q_C)$ satisfying $\beta(a)=\psi_2(\beta)(b)$.

Let $K\coloneqq\Ker(\pi_1^\top(C)^{ab}\twoheadrightarrow \pi_1^\top(X)^{ab})$. Then one could obtain the following commutative diagram of exact sequences
\begin{equation*}
\begin{adjustbox}{max width=\textwidth}
    \xymatrix{\bigoplus\limits_{1\leq i<j\leq n}(K\otimes \pi_1^\top(C)^{ab}\oplus \pi_1^\top(C)^{ab}\otimes K)\ar[r]\ar@{->>}[d]&\bigoplus\limits_{1\leq i<j\leq n}  \pi_1(C)^{ab}\otimes \pi_1^\top(C)^{ab}\ar@{->>}[rr]^{\oplus (i_*\otimes i_*)}\ar@{->>}[d]^{q_C}&&\bigoplus\limits_{1\leq i<j\leq n}  \pi_1^\top(X)^{ab}\otimes \pi_1^\top(X)^{ab}\ar@{->>}[d]^{q_X}\\
    K\otimes \pi_1^\top(C)^{ab}\ar[r]&\wedge^2 \pi_1^\top(C)^{ab}\ar@{->>}[rr]^{\wedge^2i_*}&&\wedge^2 \pi_1^\top(X)^{ab}.
    }
\end{adjustbox}
\end{equation*} 
Thus by diagram chasing, one has the surjection $\Ker(q_C)\twoheadrightarrow \Ker(q_X)$. Now, suppose that $a\in \Ker(q_X)$. One could lift $a$ to $b \in\Ker( q_C)$, and hence, $\beta(a)=\psi_2(\beta)(b)=0$. The proof is complete.
\end{proof}

\begin{prop}\label{thm-DT-2-prop1}
Let $k$ be an algebraically closed field and $X$ be a quasi-projective smooth integral $k$-variety. 
For any $n\geq 2$ and any torsion abelian group $F$, we have
\begin{itemize}
    \item[(1)] $\Im(p^*\circ (\epsilon_{\tSym^nX},\tr_X^2))=\Im(p^*)\subset H^2(X^n,F)^{S_n}$;
    \item[(2)] $\Im(p^*\circ (\epsilon_{\tilde{V}(X)},\tr_X^2))=\Im(p^*)\subset H^2(X^n,F)^{S_n}$.
\end{itemize}
 The same holds for the topological case if $k=\BC$ and $F$ is an arbitrary abelian group.
\end{prop}

\begin{proof}
    We first prove these under topological setting. Since $\pi_1^\top(\tSym^nX)\cong \pi_1^\top(\tilde{V}(X))\cong \pi_1^\top(X)^{ab}$ (see Theorem \ref{thm-DT-1} and Remark \ref{rem-DT1-top}), the statement (1) is equivalent to the statement (2) and we only prove (2). 

    By K\"unneth formula, we have the short exact sequence
    (\ref{Thm-Kunneth-eq4}), and taking $S_n$-invariance, we get the following exact sequence
\begin{equation}\label{equ-exact sequence for S_n}
0\to (\bigoplus_{i=1}^n H^2_\top(X,F))^{S_n}\xrightarrow{ \oplus \pr_i^*}  H^2_\top(X^n,F)^{S_n}\xrightarrow{\lambda_X }
(\bigoplus_{1\leq i<j\leq n}  \Hom_{\BZ} ( \pi_1^\top(X)^{ab}\otimes \pi_1^\top(X)^{ab} ,F) )^{S_n}. 
\end{equation}
Since $S_n$ acts transitively on the $n$ factors of $X^n$, we have $(\bigoplus_{i=1}^n H^2_\top(X,F))^{S_n}\cong H^2_\top(X,F)$ and under this isomorphism, $H^2_\top(X,F)\subset \bigoplus_{i=1}^n H^2_\top(X,F)$ is the diagonal embedding. On the other hand, by Proposition~\ref{eq-def-transfer1}, the morphism $p^*\circ \tr_X^2\colon H^2_\top(X,F)\to H^2_\top(X^n,F)^{S_n}$ is equal to the morphism $H^2_\top(X,F)\cong (\bigoplus_{i=1}^n H^2_\top(X,F))^{S_n} \xrightarrow{ \oplus \pr_i^*}  H^2_\top(X^n,F)^{S_n}$.

Suppose that $\alpha\in \Im(p^*\colon H^2_\top(\tilde{V}(X),F)\to H^2_\top(X^n,F))\subset H^2_\top(X^n,F)^{S_n}$, applying Lemma~\ref{lem-DT-2-prop1}, there exists $\beta\in H^2(\pi_1^\top(\tilde{V}(X)),F)$ with $\lambda_X\circ p^* \circ {\epsilon}_{\tilde{V}(X)}(\beta)=\lambda_X(\alpha)$. Hence, by (\ref{equ-exact sequence for S_n}), we have $$\alpha-p^* \circ {\epsilon}_{\tilde{V}(X)}(\beta)\in \Ker(\lambda_X)=\Im(\bigoplus_{i=1}^n H^2_\top(X,F))^{S_n}\xrightarrow{ \oplus \pr_i^*}  H^2_\top(X^n,F)^{S_n})=\Im(p^*\circ \tr_X^2).$$
This leads to the equation (2) for the topological cases.

For a general algebraically closed field $k$ of characteristic zero, passing to the limit, one may assume that all varieties appearing in Theorem \ref{thm-DT-2} are defined over an algebraically closed subfield $L\subset k$ with $L\subset \BC$. Then by the rigidity of fundamental groups and of \'etale cohomology (\cite[Cor. VI.4.3]{Mi80}), and by the comparison theorems (see (\ref{eq-fund-gps-5}) and \cite[Cor. III.3.12]{Mi80}), we may assume $k=L$. Then we have isomorphisms
\begin{equation*}
  \pi_1(\tilde{V}(X))\cong  \pi_1(\tilde{V}(X)) \cong \widehat{ \pi_1^{\top}(\tilde{V}(X)(\BC))},    
\end{equation*}
\begin{equation*}
   H^i(X,F )\cong H^i(X_{\BC},F) \cong H^i_{\top}(X(\BC),F) \quad\text{ and } \quad H^i(\tilde{V}(X),F )\cong H^i(\tilde{V}(X)_{\BC},F) \cong H^i_{\top}(\tilde{V}(X)(\BC),F) .
\end{equation*}
By applying Proposition \ref{prop-compatible of edge morphisms}, one could reduce to the case $k=\BC$, noting that $\pi_1^\top(\tilde{V}(X)(\BC))$ is abelian by Theorem \ref{thm-DT-1} (i). Hence (2) and also (1) hold. The proof is complete.
\end{proof}

If the group $F$ is injective, e.g., $F=\BQ/\BZ$, the morphisms $p^*\circ (\epsilon_{\tSym^nX},\tr_X^2)$ and  $p^*\circ (\epsilon_{\tilde{V}(X)},\tr_X^2)$ are actually injective; see Corollary \ref{cor-iso for Q/Z}.

\begin{cor}\label{cor-iso for Q/Z}
Let $k$ be an algebraically closed field and let $X$ be a quasi-projective smooth integral $k$-variety. For any $n\geq 2$ and any injective torsion abelian group $F$, the morphisms $p^*\circ (\epsilon_{\tSym^nX},\tr_X^2)$ and  $p^*\circ (\epsilon_{\tilde{V}(X)},\tr_X^2)$ are injective. Moreover, the same holds for the topological case, either if $F$ is injective or if $H_1^\top(X,\BZ)$ is torsion free.
\end{cor}

\begin{proof}
 Similar to the proof of Proposition \ref{thm-DT-2-prop1}, we only consider the topological cases and the injectivity of $p^*\circ (\epsilon_{\tilde{V}(X)},\tr_X^2)$. Suppose either $F$ is injective or $H_1^\top(X,\BZ)$ is torsion free, by universal coefficient theorem for group cohomology (see \cite[Prop. I.0.8 \& Thm. V.6.4]{brown2012cohomology}), we have 
\begin{equation*}
    H^2(\pi_1^{\top}(\tilde{V}(X)),F)\cong \Hom_\BZ(H_2(\pi_1^{\top}(\tilde{V}(X))),F)\cong \Hom_\BZ(\wedge^2\pi_1^{\top}(\tilde{V}(X))),F).
\end{equation*}
For every $(\alpha,\beta)\in \Ker(p^*\circ (\epsilon_{\tilde{V}(X)},\tr_X^2))\subset H^2(\pi^\top_1(\tilde{V}(X)),F)\oplus H^2_\top(X,F)$, then we have 
\begin{equation*}
    0=\lambda_X\circ p^*\circ (\epsilon_{\tilde{V}(X)},\tr_X^2)(\alpha,\beta)=\lambda_X\circ p^*\circ \epsilon_{\tilde{V}(X)}(\alpha).
\end{equation*}
Since $F$ is injective, the morphism $\lambda_X\circ p^*\circ \epsilon_{\tilde{V}(X)}$ is exactly the composition $(\lambda^X)^\vee\circ (p_*)^\vee \circ {\epsilon}^{\tilde{V}(X)}$, which  is the $\Hom_\BZ(-,F)$ of the surjection (\ref{equ-Hom from tensor to wedge}). Hence, $\lambda_X\circ p^*\circ \epsilon_{\tilde{V}(X)}$ is injective and then one has $\alpha=0\in H^2(\pi^\top_1(\tilde{V}(X)),F)$. In addition, we have $\beta=0\in H^2_\top(X,F)$, because $p^*\circ \tr_X^2$ is exactly the morphism $(\bigoplus_{i=1}^n H^2_\top(X,F))^{S_n}\xrightarrow{ \oplus \pr_i^*}  H^2_\top(X^n,F)^{S_n}$ which is an injection. The proof is complete.
\end{proof}

\begin{rem}
 The morphisms $p^*\circ (\epsilon_{\tSym^nX},\tr_X^2)$ and  $p^*\circ (\epsilon_{\tilde{V}(X)},\tr_X^2)$ are generally not surjective. 
 For instance, let $G\coloneqq \BZ/2$ act on an affine space $\BA^N$ with sufficiently large $N$, and let $X=U/G$, where $U\subset \BA^N$ is an open subset such that $G$ acts on $U$ freely. Then $H^2_\top(X,F)=H^2(G,F)$ by Cartan--Leray spectral sequence. In particular, if $F$ is injective, then $H^2_\top(X,F)\cong \Hom_\BZ(\wedge^2G,F)=0$, and hence the K\"unneth formula says that $H^2(X^n,F)^{S_n}\xrightarrow{\lambda_X} (\bigoplus_{1\leq i<j\leq n}  \Hom_{\BZ} ( \pi_1^\top(X)^{ab}\otimes \pi_1^\top(X)^{ab} ,F) )^{S_n}$. Since the group $S_n$ acts transitively on the set of unordered pairs
$\bigl\{\{i,j\}\mid 1\leq i<j\leq n\bigr\}$, it follows that an $S_n$-invariant element
$$(\varphi_{ij})_{i<j}\in\bigoplus_{1\leq i<j\leq n} \Hom_{\BZ} ( \pi_1^\top(X)^{ab}\otimes \pi_1^\top(X)^{ab} ,F) $$
is uniquely determined by its $(1,2)$-component $\varphi_{12}$, provided that this component is invariant under the stabilizer of $\{1,2\}$. Note that the only nontrivial condition on $\varphi_{12}$ comes from the transposition $(12)\in S_n$. Since the two copies of $\pi_1^\top(X)^{ab}\cong H_1(X,\BZ)\cong \BZ/2$ arise from the degree-one K\"unneth factors corresponding to the first and second factors of $X^n$, interchanging these two factors introduces the Koszul sign $(-1)^{1\cdot 1}=-1$. Therefore, the element $(12)\in S_n$ acts on $\Hom_{\BZ} ( \pi_1^\top(X)^{ab}\otimes \pi_1^\top(X)^{ab} ,F)$ by $((12)\cdot\varphi)(a\otimes b)=-\varphi(b\otimes a)$ for $a,b\in \pi_1^\top(X)^{ab}\cong H_1(X,\BZ)$.
Thus, $\varphi_{12}\in \Hom_{\BZ} ( \pi_1^\top(X)^{ab}\otimes \pi_1^\top(X)^{ab} ,F)$ is invariant under $(12)\in S_n$ if and only if $\varphi_{12}(a\otimes b)=-\varphi_{12}(b\otimes a)$ for all $a,b\in \pi_1^\top(X)^{ab}$,i.e., $\varphi_{12}(a\otimes b+b\otimes a)=0$ for all $a,b\in \pi_1^\top(X)^{ab}$. Write
\begin{equation*}
{\rm Alt}(\pi_1^\top(X)^{ab})\coloneqq \pi_1^\top(X)^{ab}\otimes_{\BZ}\pi_1^\top(X)^{ab}/<a\otimes b+b\otimes a\mid a,b\in \pi_1^\top(X)^{ab}>\cong \BZ/2.    
\end{equation*}
Then we have
\begin{equation*}
    (\bigoplus_{1\leq i<j\leq n}  \Hom_{\BZ} ( \pi_1^\top(X)^{ab}\otimes \pi_1^\top(X)^{ab} ,F) )^{S_n}\cong \Hom_\BZ({\rm Alt}(\pi_1^\top(X)^{ab}),F)\cong \Hom_\BZ(\BZ/2,F).
\end{equation*}

However, for injective $F$, the morphism $\lambda_X\circ p^*\circ \epsilon_{\tilde{V}(X)}$ is zero as $$\wedge^2\pi_1^\top(X)^{ab}=\pi_1^\top(X)^{ab}\otimes_{\BZ}\pi_1^\top(X)^{ab}/<a\otimes a\mid a\in \pi_1^\top(X)^{ab}>=0,$$
and hence the morphism $p^*\circ \epsilon_{\tilde{V}(X)}$ is zero by noting that $\lambda_X$ is an isomorphism. Therefore, the morphisms $p^*\circ (\epsilon_{\tSym^nX},\tr_X^2)$ and  $p^*\circ (\epsilon_{\tilde{V}(X)},\tr_X^2)$ are zero but the target $H^2(X^n,F)^{S_n}\cong \Hom_\BZ(\BZ/2,F)$ is non-zero if $F[2]\neq 0$. 
\end{rem}

\begin{lem}\label{thm-DT-2-lem-1}
Let $k$ be an algebraically closed field. Let $X$ and $Y$ be two quasi-projective smooth geometrically integral $k$-varieties of dimension $\geq 2$ with an immersion morphism $f\colon Y\to X$. 
For any $n\geq 2$ and any torsion \'etale sheaf $F$, we consider the induced commutative diagram of varieties
\begin{equation*}
  \xymatrix{U^n(Y)\ar[r]^{f_U^n}\ar[d]^{\psi_Y} &  U^n(X)\ar[d]^{\psi_X}\\
U(Y)\ar[r]^{f_U}& U(X) ,
}  
\end{equation*}
and the induced commutative diagram of cohomology groups
\begin{equation*}
 \xymatrix{H^2(U(X),F)\ar[r]^{f_U^*} \ar[d]^{\psi_X^*}  &  H^2(U(Y),F) \ar[d]^{\psi_Y^*}\\ 
H^2(U^n(X),F)\ar[r]^{f_U^{n,*}} & H^2(U^n(Y),F).
}   
\end{equation*}
Let us consider the induced homomorphism $f_U^*|_{\Ker}\colon  \Ker(\psi_X^*)\to \Ker(\psi_Y^*)$.
\begin{itemize}
    \item[(i)]  If $f_*\colon  \pi_1(Y)\to \pi_1(X)$ is surjective, then $f_U^*|_{\Ker}$ is injective.
    \item[(ii)] If $f_*\colon  \pi_1(Y)\to \pi_1(X)$ is an isomorphism, then $f_U^*|_{\Ker}$ is an isomorphism.
\end{itemize} 
These also hold for the topological case if $k=\BC$ and $F$ is an abelian group.
\end{lem}

\begin{proof}
We only prove (i) for \'etale cohomology, and (ii) follows from the same argument. For the topological version, the same proof applies without any changes.

Since $f_*\colon  \pi_1(Y)\to \pi_1(X)$ is surjective, 
by (\ref{eq-fund-gps-3}), $f^*|_{H^1}\colon  H^1(X,F)\to H^1(Y,F)$ is injective.
By the {K\"unneth formula} (see \cite[Thm. 2.1]{C4c}), we have an exact sequence of $S_n$-modules
$$ 0\to H^1(X^n,F)\to H^1(Y^n,F) \to \Coker(f^*|_{H^1})^{\oplus n}  \to 0 ,$$
where the action of $S_n$ on $\Coker(f^*|_{H^1})^{\oplus n} $ is permutation.
Hence the morphism
$$ H^1(Y^n,F)^{S_n}\cong (H^1(Y,F)^{\oplus n})^{S_n} \to (\Coker(f^*|_{H^1})^{\oplus n})^{S_n}   $$
is surjective. Therefore, $\phi_1\colon  H^1(S_n, H^1(X^n,F))\to H^1(S_n,H^1(Y^n,F)) $ is injective.

Since $\dim(Y)\geq 2$, one has $\codim(Y^n\setminus U^n(Y),Y^n)\geq 2$ and $H^1(Y^n,F)\cong H^1(U^n(Y),F)$.
The same statement holds for $X$.
By the purity of \'etale cohomology,
the homomorphism
$$\bar{\phi_1}\colon   H^1(S_n, H^1(U^n(X),F))\cong  H^1(S_n, H^1(X^n,F))   \xrightarrow{\phi_1} 
 H^1(S_n, H^1(Y^n,F))\cong  H^1(S_n, H^1(U^n(Y),F))  $$
is injective.

On the other hand, the composition 
$$H^1(\tilde{V}(Y),F)\to H^1(U(Y),F) \xrightarrow{\phi_2} H^1(U^n(Y),F)^{S_n}\cong H^1(Y^n,F)^{S_n}$$
is exactly the $p^*$ in (\ref{eq-def-p}), and it is an isomorphism by Theorem \ref{thm-DT-1}.
Therefore, $\phi_2$ is surjective.

We consider the Hochschild--Serre spectral sequence (\cite[Thm. III.2.20]{Mi80}) for the Galois covers $\psi_X$ and $\psi_Y$,
and we get the commutative diagram with exact rows
$$\xymatrix{ & H^2(S_n,F)\ar[r]\ar[d]^=  & 
\Ker(\psi_X^*) \ar[r] \ar[d]^{f_U^*}   &   H^1(S_n, H^1(U^n(X),F))\ar[r]  \ar[d]^{\bar{\phi_1}}    &  H^3(S_n,F)\ar[d]^=   \\
0 \ar[r]  &  H^2(S_n,F)\ar[r]^{\phi_3}  & 
\Ker(\psi_Y^*) \ar[r]   &   H^1(S_n, H^1(U^n(Y),F))\ar[r]   &  H^3(S_n,F) ,
}$$
where $\phi_3$ is injective because $\phi_2$ is surjective.
Since $\bar{\phi_1} $ is injective,
a diagram chasing implies (i).
\end{proof}

\begin{prop}\label{thm-DT-2-prop2}
Let $k$ be an algebraically closed field and $X$ be a quasi-projective smooth integral $k$-variety with $\dim(X)\geq 2$. We have $\Ker(p^*)=\Im(i_{D_X}) $ for any $n\geq 2$, where $i_{D_X}\colon F\to H^2(\tilde{V}(X),F)$ is the Gysin morphism. 
The same holds for the topological case if $k=\BC$ and $F$ is an abelian group.
\end{prop}

\begin{proof}
We prove this for \'etale cohomology and for topological version, the same proof applies without any changes. Take $U\coloneqq U^n(\BA^2)\subset \BA^{2n}$ with free action of $S_n$.
In this case, by Leray spectral sequence and noting that $H^i(U,F)\cong H^i(\BA^{2n},F)=0$ for $0<i\leq 2$ and $\pi_1(U)\cong \pi_1(\BA^{2n})=1$,
for any smooth variety $Y$, the pullback induces isomorphism
$H^i(Y,F)\iso H^i(Y\times U,F) $ for $0\leq i\leq 2$, and
if $S_n$ acts on $Y$, then the quotient $(Y\times U)/S_n$ is smooth.

Let $C\subset X$ be an affine sub-curve such that the morphism $\pi_1^\top(C)\to \pi_1^\top(X)$ is surjective; see Proposition \ref{prop-Lef}.  We consider the following commutative diagram
\begin{equation*}
 \xymatrix{  C^n \ar[d] & C^n\times U\ar[r]\ar[d]^{q_C} \ar[l]& X^n\times U \ar[r]\ar[d]^{q_X} &  X^n\ar[d]  \ar@{-->}@/^3pc/[dd] \\
\Sym^n C  \ar@{-->}[rrrd] & (C^n\times U)/S_n \ar[r] \ar[l]& ( X^n\times U )/S_n\ar[r] \ar@{-->}[rd]  & \Sym^n(X)  \\
       & & &                              \tilde{V}(X) \ar[u] ,
}   
\end{equation*}
where all rational maps are well defined over codimension $1$ points. This induces the following commutative diagram of cohomology groups by purity of \'etale cohomology
\begin{equation}\label{equ-dia1}
 \xymatrix{H^2(\tilde{V}(X),F)\ar[r]_-{\psi}\ar[d]^{p^*}  \ar@/^1pc/[rrr]& H^2((X^n\times U)/S_n,F) \ar[r]\ar[d]^{q_X^*} &
H^2((C^n\times U)/S_n,F) \ar[d]^{q_C^*}  &   H^2(\Sym^nC,F)\ar[l] \ar[d]^{p_C^*}\\
H^2(X^n,F)\ar[r]^{\cong} &  H^2(X^n\times U,F)  \ar[r]  & H^2(C^n\times U,F) & H^2(C^n,F).\ar[l]_{\cong}
}   
\end{equation}
\begin{claim}
  The morphism $\Ker(q_X^*)\to \Ker(q_C^*)$  induced by the commutative diagram (\ref{equ-dia1}) is injective.
\end{claim}
\begin{proof}[Proof of the claim]
By definition, $C^n\times U \cong C^n\times U^n(\BA^2)$ is an open subset of $U^n(C\times \BA^2)$ with complement of codimension $\geq 2$.
Thus $(C^n\times U)/S_n $ is an open subset of smooth variety $U(C\times \BA^2)$ with complement of codimension $\geq 2$.
By the purity of \'etale cohomology again,
one gets the following natural isomorphism
    \begin{equation*}
        H^2(U(C\times \BA^2),F)\xrightarrow[\sim]{\rm restr.}H^2((C^n\times U)/S_n,F).
    \end{equation*}
  Similarly, we obtain the same isomorphism if we replace $C$ by $X$.
    In this case, we can identify the middle part of diagram (\ref{equ-dia1}) with the following commutative diagram
    \begin{equation*}
        \xymatrix{H^2(U(X\times \BA^2),F)\ar[r]\ar[d]&H^2(U(C\times \BA^2),F)\ar[d]\\
        H^2(U^n(X\times \BA^2),F)\ar[r] & H^2(U^n(C\times \BA^2),F),
        }
    \end{equation*}
    where the horizontal arrows are induced by the immersion $f\colon C\times \BA^2\to X\times \BA^2$  and the vertical arrows are induced by $S_n$-quotient. The proof is complete after applying the surjectivity of $\pi_1(C)\to \pi_1(X)$ and Lemma \ref{thm-DT-2-lem-1} (i).
\end{proof}

The curve case implies that $p_C^*$ is injective; see Proposition \ref{lem-DT-1} (iv) and Proposition \ref{prop-curves} (iv).
The diagram chasing shows that $\Ker(p^*)=\Ker(\psi)$.
Since the projection morphism $$(U^n(X)\times \BA^N)/S_n\to U^n(X)/S_n=U(X)$$ is a $\BA^N$-bundle 
and one has $\codim(\BA^N\setminus U,\BA^N)\geq 2$, 
then the natural pull-back $$H^2(U(X),F)\xrightarrow{\cong} H^2((U^n(X)\times \BA^N)/S_n,F)\cong H^2((U^n(X)\times U)/S_n,F) $$ is an isomorphism. 
Therefore, $\Ker(p^*)=\Ker(\psi)=\Ker(H^2(\tilde{V}(X),F)\xrightarrow{\rm restr.} H^2(U(X),F))=\Im(i_{D_X})$. The proof is complete.
\end{proof}

\begin{proof}[Proof of Theorem \ref{thm-DT-2}]
If $\dim(X)=1$, noting that $\codim(\Sym^nX\setminus V(X),\Sym^nX)=2$, to prove Theorem~\ref{thm-DT-2}, it suffices to show the following morphism is surjective (and then it is an isomorphism by Corollary~\ref{cor-iso for Q/Z})
\begin{equation}\label{equ-proof for curve case}
    (\epsilon_{\Sym^nX},\tr_X^2)\colon   H^2_{cont}(\pi_1(\Sym^nX),F)\oplus H^2(X,F) \to H^2(\Sym^nX,F).
\end{equation}
If $X$ is an affine curve, the surjectivity of (\ref{equ-proof for curve case}) is implied by Proposition \ref{prop-curves} (ii). Suppose that $X$ is a projective curve. Pick a closed point $x\in X$ with the affine open complement $U\coloneqq X\setminus \{x\}$. Consider the following commutative diagram of exact sequences (i.e., Gysin sequences)
\begin{equation*}
\xymatrix{H^0(x,F)\ar[r]^{\sim}\ar[d]^{\tr}_{\cong}& H^2(X,F)\ar[r]\ar[d]^{\tr^2_X}&H^2(U,F)=0\ar[d]^{\tr^2_U}\\
H^0(Z,F)\ar@{^(->}[r]&H^2(\Sym^nX,F)\ar[r]&H^2(\Sym^nU,F),
}
\end{equation*}
where $Z\coloneqq \Sym^nX\setminus \Sym^nU$ is isomorphic to $\Sym^{n-1}X$, and the morphism ``$\tr$" is a correspondence action induced by the correspondence $W_0|_x\cong x\times \Sym^{n-1}X\cong x\times Z$; see (\ref{equ-universal family}). 
Thus, the morphism $\tr$ is an isomorphism as $Z$ is smooth and connected. Here the first square is commutative, because both $Z$ and $\Sym^nX$ are smooth and projective so that the Gysin map can be identified with the pushout of homology groups, and then we could apply \cite[Prop. 16.1.1]{fulton}. Therefore, we have 
\begin{equation}\label{equ-proof of thm 4.7}
    \Ker(H^2(\Sym^nX,F)\xrightarrow{\rm restr.} H^2(\Sym^nU,F))=\Im(H^2(X,F)\xrightarrow{\tr^2_X}H^2(\Sym^nX,F)).
\end{equation}
On the other hand, we have another commutative diagram induced by the open immersion $U\subset X$
\begin{equation*}
\xymatrix{    H^2_{cont}(\pi_1(\Sym^nX),F)\ar[rrr]^{\epsilon_{\Sym^nX}}\ar[d]^\cong&&&H^2(\Sym^nX,F)\ar[d]^{\rm restr.}\\
    H^2_{cont}(\pi_1(\Sym^nU),F)\ar[rrr]^{\cong}_{\epsilon_{\Sym^nU}}&&&H^2(\Sym^nU,F),
}
\end{equation*}
where the left arrow is isomorphism since $\pi_1(\Sym^nX)\cong \pi_1(X)^{ab}\cong \pi_1(U)^{ab}\cong \pi_1(\Sym^nU)$. By diagram chasing and (\ref{equ-proof of thm 4.7}), one could get the surjectivity of (\ref{equ-proof for curve case}) for the projective curve $X$.

We may assume that $\dim(X)\geq 2$. We first suppose that $k=\BC$ and hence we work on topological singular cohomology. For every $\alpha\in H^2_\top(\tilde{V}(X),F)$, by applying Proposition \ref{thm-DT-2-prop1} (2), one could find $(\beta,\gamma)\in H^2(\pi_1^\top(\tilde{V}(X)),F)\oplus H^2_\top(X,F)$ such that 
		\begin{equation*}
			p^*(\alpha)=p^*\circ (\epsilon_{\tilde{V}(X)},\tr_X^2)(\beta,\gamma)\in H^2_\top(X^n,F).
		\end{equation*}
	This means that $\alpha-(\epsilon_{\tilde{V}(X)},\tr_X^2)(\beta,\gamma)\in \Ker(p^*)=\Im(i_{D_X})\in H^2_\top(\tilde{V}(X),F)$ by Proposition \ref{thm-DT-2-prop2}. Hence, there exists $\omega\in F$ such that 
	\begin{equation*}
		(\epsilon_{\tilde{V}(X)},\tr_X^2,i_{D_X})(\beta,\gamma,\omega)=\alpha,
	\end{equation*}
which means that Theorem \ref{thm-DT-2} holds for the case $k=\BC$.

For a general algebraically closed field $k$ of characteristic zero, passing to the limit, one may assume that all varieties appearing in Theorem \ref{thm-DT-2} are defined over an algebraically closed subfield $L\subset k$ with $L\subset \BC$. Then by the rigidity of fundamental groups and of \'etale cohomology (\cite[Cor. VI.4.3]{Mi80}), and by the comparison theorems (see (\ref{eq-fund-gps-5}) and \cite[Cor. III.3.12]{Mi80}), we may assume $k=L$. Then we have isomorphisms
\begin{equation*}
  \pi_1(\tilde{V}(X))\cong  \pi_1(\tilde{V}(X)) \cong \widehat{ \pi_1^{\top}(\tilde{V}(X)(\BC))},    
\end{equation*}
\begin{equation*}
   H^i(X,F )\cong H^i(X_{\BC},F) \cong H^i_{\top}(X(\BC),F) \quad\text{ and } \quad H^i(\tilde{V}(X),F )\cong H^i(\tilde{V}(X)_{\BC},F) \cong H^i_{\top}(\tilde{V}(X)(\BC),F) .
\end{equation*}
By applying Proposition \ref{prop-compatible of edge morphisms}, one could reduce to the case $k=\BC$, noting that $\pi_1^\top(\tilde{V}(X)(\BC))$ is abelian by Theorem \ref{thm-DT-1} (i). The proof is complete.
\end{proof}
The following corollary implies Theorem \ref{thm-2}.
\begin{cor}\label{cor-DT thm for deg2}
Let $k$ be algebraically closed and $X$ be a quasi-projective smooth integral $k$-variety. Let $\tSym^nX$ be a proper resolution of $\Sym^nX$. Then, for any $n\geq 2$,
\begin{itemize}
    \item[(1)] the \'etale fundamental groups $\pi_1(\tSym^nX) \cong \pi_1(X)^{ab}$ is commutative;
    \item[(2)] the $p^*$ in (\ref{eq-def-p}) induces isomorphisms 
    $H^1(\tSym^nX,F)\iso   H^1(X^n,F)^{S_n} $;
    \item[(3)] the canonical transfer $\tr_X^1\colon  H^1(X,F)\to H^1(\tSym^nX,F)$ in (\ref{eq-def-transfer}) is an isomorphism;
    \item[(4)] any torsion \'etale sheaf $F$, the following homomorphism is surjective 
\begin{equation}\label{equ-unramified coh}
    ( \epsilon'_{\tSym^nX} , \tr_{nr}): H^2_{cont}(\pi_1(\tSym^nX),F) \oplus H^2_{nr}(X,F)  \twoheadrightarrow H^2_{nr}(\tSym^nX,F),
\end{equation}
where
\begin{itemize}
\item $H^2_{nr}(\dagger,F)$ is the second unramified cohomology of $\dagger$ with $F$-coefficients (cf. \cite[\S 4]{colliot1995birational});
\item the morphism $\epsilon'_{\tSym^nX}\colon H^2_{cont}(\pi_1(\tSym^nX),F)\to H^2_{nr}(\tSym^nX,F)$ is the composition of 
\begin{equation*}
    H^2_{cont}(\pi_1(\tSym^nX),F)\xrightarrow{\epsilon_{\tSym^nX}} H^2(\tSym^nX,F)\xrightarrow{\rm restr.} H^2_{nr}(\tSym^nX,F);
\end{equation*}
\item the morphism $H^2_{nr}(X,F)  \xrightarrow{\tr_{nr}} H^2_{nr}(\tSym^nX,F)$ is induced by the correspondence action $W$ (\ref{equ-def-W}).
\end{itemize}
In particular, the morphism
\begin{equation}\label{equ-surj for Br}
    ( \epsilon'_{\tSym^nX} , \tr_{X,\Br})\colon H^2_{cont}(\pi_1(\tSym^nX),\mu_{\infty}) \oplus \Br(X)  \twoheadrightarrow \Br(\tSym^nX)
\end{equation}
is surjective.
\end{itemize}

The same holds for the case $k=\BC$ with topological cohomology groups.
\end{cor}

\begin{proof}
We only prove these under \'etale topological setting, and for classical topological case, the proof works words by words. 

For (1)-(3), since that there groups are all proper birational invariants for smooth varieties (see \cite{colliot1995birational}), it suffices to prove these for some particular resolution $\tSym^nX$. Therefore, one could get (1)-(3) from Theorem \ref{thm-DT-1} (i)-(iii).

For (4), we only prove (\ref{equ-unramified coh}) as we have $\Br(X)\cong H^2_{nr}(X,\mu_{\infty})$ for any smooth quasi-projective $X$. 
For the topological case, the proof works words by words. 

We first recall some basic facts about unramified cohomology. For every smooth $X$, note that the restriction morphism $H^2(X,F)\to H^2_{nr}(X,F)$ is surjective, and for any irreducible divisor $D\subset X$, the composition 
\begin{equation*}
    F\cdot[D]\xrightarrow{\rm Gysin} H^2(X,F)\xrightarrow{\rm restr.} H^2_{nr}(X,F)
\end{equation*}
is zero by Bloch--Ogus theory; see \cite[Prop. 3.9 \& Thm. 4.2]{BO74}. Moreover, if $X\xrightarrow{f} Y$ is a proper birational morphism between smooth varieties, then we have the isomorphism $f\colon H^2_{nr}(Y,F)\xrightarrow{\sim}H^2_{nr}(X,F)$; 
cf. \cite[Prop. 2.1.8 \& Thm. 4.1.1]{colliot1995birational}. In fact, \cite{colliot1995birational} only proves that it is a birational invariance for smooth proper varieties. But the arguments hold for proper morphism of smooth quasi-projective varieties by applying codimension one purity theorem \cite[Thm. 3.8.2]{colliot1995birational} and the valuative criterion of properness \cite[Thm. II.4.7]{Hartshorne}.

For every $\alpha\in H^2_{nr}(\tSym^nX,F)$, one could lift it to $\alpha'\in H^2(\tSym^nX,F)$. By Theorem \ref{thm-DT-2}, there exists $\beta\in \ H^2_{cont}(\pi_1(\tSym^nX),\mu_{\infty})$, $\gamma\in H^2(X,F)$ and $\omega\in F^{\oplus m}$, such that $\alpha'=(\epsilon_{\tSym^nX},\tr_X^2,i_*)(\beta,\gamma,\omega)$. Since transfer morphisms are compatible with restriction, then we have 
\begin{equation*}
    \alpha={\rm restr.}(\alpha')={\rm restr.}\cdot(\epsilon_{\tSym^nX},\tr_X^2,i_*)(\beta,\gamma,\omega)=(\epsilon'_{\tSym^nX},\tr_{nr})(\beta,{\rm restr.}(\gamma))\in H^2_{nr}(\tSym^nX,F).
\end{equation*}
Because we have ${\rm restr.}\cdot i_*(\omega)=0\in H^2_{nr}(\tSym^nX,F)$. The proof is complete.
\end{proof}

\begin{rem}\label{rem-S4-last}
    The surjection (\ref{equ-surj for Br}) is not injective in general. For instance, if $X$ is an elliptic curve, then for any $n\geq 2$, $\Sym^nX\to J(X)\cong X$ is a projective bundle by Riemann--Roch theorem, and hence $$H^2_{nr}(\Sym^nX,F)\cong H^2_{nr}(X,F)=0;$$ see \cite[Lem. 4.5]{schreieder2021unramified}.  However, we have $\pi_1(\Sym^nX)\cong \pi_1(X)\cong \hat{\BZ}\oplus \hat{\BZ}$, which implies that $$H^2_{cont}(\pi_1(\Sym^nX),F)\cong F\neq 0.$$ 

    On the other hand, we cannot expect that we always have an isomorphism $$\frac{\Br(X)}{\Br(k)}\cong \frac{\Br(\tSym^nX)}{\Br(k) }$$ for general $X$, even over algebraically closed fields. 
    In fact, if $X$ is a smooth projective complex curve with $g(C)\geq 2$, then $\Br(\Sym^nX)\cong \Br(J(X))$ for sufficiently large $n$, by Riemann--Roch theorem (\emph{cf.} \cite[Thm. IV.1.1.3]{Hartshorne}) and \cite[Lem. 4.5]{schreieder2021unramified}. Since abelian varieties always have torsion-free integral cohomology groups, one has that $\Br(J(X))\cong (\BQ/\BZ)^{b_2-\rho}$, where $b_2$ is the dimension of the $\BQ$-vector space $H^2(J(X),\BQ)$ and $\rho$ is the Picard number of $J(X)$; see \cite[Prop. 4.2.3]{colliot1995birational}. Note that $J(X)$ is a complex abelian variety of dimension $g$. Then we have $b_2=\binom{2g}{2}=g(2g-1)$, while $$\rho\leq h^{1,1}(J(X))=b_2-2h^{2,0}(J(X))<b_2$$ as $g>1$ and hence $h^{2,0}(J(X))=\binom{g}{2}>0$. Therefore, we get that $\Br(\Sym^nX)\cong \Br(J(X))\neq 0$ while $\Br(X)=0$, since $\dim(X)=1$.
\end{rem}

\section{Descent formula and proof of Theorem \ref{thm-1}}\label{section5}

In this section, let $k$ be a number field. Let $X$ be a quasi-projective smooth geometrically integral variety.
We consider the canonical restriction homomorphism
$$\Phi_{r,X}\colon  H^2(X,\mu_{\infty})\to H^2(X_{\bk},\mu_{\infty}) ,$$
the exact sequence induced by the Kummer exact sequence
\begin{equation}\label{S5-e-Ku}
 0\to \Pic(X)\otimes \BQ/\BZ\to  H^2(X,\mu_{\infty})\xrightarrow{\Phi_{Ku,X}}  \Br(X) \to 0.
\end{equation} 
and the Hurewicz map over $X_{\bk}$ as in (\ref{def-epsilon-eq1})
$$\bar{\epsilon}_X\colon  H^2(\pi_1(X_{\bk}),\mu_{\infty})  \to  H^2(X_{\bk},\mu_{\infty}).$$
For any subgroup $H\subset H^2(X_{\bk},\mu_{\infty}) $, we denote
$$\Phi_X(H)\coloneqq  \Phi_{Ku,X}(\Phi_{r,X}^{-1}( H ) ) .$$

Now, we define the algebraic part $H^2_1(X,\mu_{\infty})\coloneqq \Phi_{r,X}^{-1}(0) $ as in (\ref{def-eq-H^i_1}),
and the natural subgroup $\Br_{2/3}(X)\coloneqq \Phi_X(0) $ of $\Br(X)$ as in  \cite[(2.2)]{C5}.

In this section, for any morphism $f\colon  Y\to X$, we denote the pullbacks $f^*,\bar{f}^*,f^*_{\Br} $ as in
\begin{equation}\label{S5-e-Ku-1}
\xymatrixcolsep{3pc}\xymatrix{
H^2(X_{\bk},\mu_{\infty})\ar[d]^{\bar{f}^*}  &   H^2(X,\mu_{\infty})\ar[d]^{f^*}\ar[l]_-{ \Phi_{r,X}}\ar[r]^-{\Phi_{Ku,X}}&
\Br(X) \ar[d]^{f^*_{\Br}}  \\
H^2(Y_{\bk},\mu_{\infty})   & H^2(Y,\mu_{\infty}) \ar[l]_-{ \Phi_{r,Y}}\ar[r]^-{\Phi_{Ku,Y}} &  \Br(Y) .
}\end{equation}

\subsection{Descent formula with respect to $\Br_{2/3}(X)$}

As we mentioned in the introduction, 
the descent method is a widely used method to study Hasse principle and the weak approximation for rational points.
It was firstly established by Colliot-Th\'el\`ene and Sansuc in \cite{CTS87}.
In this part, we establish a new precise descent formula for the adelic points orthogonal to our natural subgroup $\Br_{2/3}$.

\begin{prop}\label{desc-prop2}
Let $X$ be a quasi-projective smooth geometrically integral $k$-variety.
Let $A$ be a finite commutative $k$-group and let $f\colon Y\to X$ be an $A$-torsor with $Y$ geometrically integral.
Assume $\pi_1(X_{\bk})$ is abelian. 
Then, for any subgroup $H\subset H^2(X_{\bk},\mu_{\infty})^{\Gamma_k}$, we have
$$X(\RA_k)^{\Phi_X(H)}= \bigcup_{\sigma\in H^1(k,A)}f_{\sigma}(Y_{\sigma}(\RA_k)^{\Phi_{Y_{\sigma}}(\bar{f}_{\sigma}^*(H) ) }).$$

Moreover, the special case $H=0$ implies that
$$X(\RA_k)^{\Br_{2/3}(X)}= \bigcup_{\sigma\in H^1(k,A)}f_{\sigma}(Y_{\sigma}(\RA_k)^{\Br_{2/3}(Y_{\sigma})}).$$
\end{prop}

By the commutative diagram (\ref{S5-e-Ku-1}), it is clear that $f_{\sigma,\Br}^*(\Phi_X(H))\subset \Phi_{Y_{\sigma}}(\bar{f}_{\sigma}^*(H) )  $.
This result is strictly  stronger than the form
$$ X(\RA_k)^{\Phi_X(H)}= \bigcup_{\sigma\in H^1(k,A)}f_{\sigma}(Y_{\sigma}(\RA_k)^{ f_{\sigma,\Br}^*(\Phi_X(H))},$$
because it is possible that an element $h\in H$ does not belong to $\Im(\Phi_{r,X})$, 
but $\bar{f}_{\sigma}^*(h) $ belongs to $\Im(\Phi_{r,Y_{\sigma}})$.
This strong form is essential in the proof of Theorem \ref{thm-1}.

To prove this result, we recall the notions of invariant Brauer subgroup defined in \cite[Def. 1.3]{C4} and in \cite[Def. 1.1]{C5}.
Let $X$ be a smooth geometrically integral variety with an action of a connected algebraic group $G$, 
we denote by $p_1\colon  G\times X\to G, p_2\colon  G\times X\to X $ two projections and $\rho\colon  G\times X\to X$ the action.
we define the invariant Brauer subgroup as the natural subgroup
$$\Br_G(X)\coloneqq \{b\in \Br(X)  \ |\  (\rho^*_{\Br}(b)-p_{2,\Br}^*(b)) \in p_{1,\Br}^*\Br(G) \} \subset \Br(X)$$
and its analogue
$$H^2_G(X,\mu_{\infty})\coloneqq \{a\in H^2(X,\mu_{\infty})\ |\ (\rho^*(a)-p_2^*(a)) \in p_1^* H^2(G,\mu_{\infty})   \}\subset   H^2(X,\mu_{\infty})  .$$
By definition, $\Phi_{Ku,X}(H^2_G(X,\mu_{\infty})) \subset \Br_G(X)$.

\begin{lem}\label{desc-prop2-lem1-lem}
We have $H^2_G(X,\mu_{\infty}) = \Phi_{r,X}^{-1} ( H^2_{G_{\bk}}(X_{\bk},\mu_{\infty})) $.
\end{lem}

\begin{proof}
By definition,  $\Phi_{r,X}( H^2_G(X,\mu_{\infty})) \subset  H^2_{G_{\bk}}(X_{\bk},\mu_{\infty})$.
For the converse, let $a\in H^2(X,\mu_{\infty})$ be an element with $\bar{a}\coloneqq \Phi_{r,X}(a)\in H^2_{G_{\bk}}(X_{\bk},\mu_{\infty})$.
By \cite[Prop. 3.5 (3)]{C5}, $(\bar{\rho}^*-\bar{p}_2^*)(\bar{a})=0 $.
Hence $ \rho^*(a)-p_2^*(a)$ belongs to the algebraic part 
\begin{equation*}
 H^2_1(G\times X,\mu_{\infty})\cong H^2_{1,e}(G,\mu_\infty)\oplus H^2(X,\mu_\infty) ;  
\end{equation*}
see (\ref{def-eq-H^i_1}) and \cite[(2.6)]{C5}. Since $ (\rho^*(a)-p_2^*(a))|_{e_G\times X}=0 $, we have 
\begin{equation*}
 \rho^*(a)-p_2^*(a)\in p_1^*H^2_{1,e}(G,\mu_{\infty})\subset  p_1^*H^2(G,\mu_{\infty})   
\end{equation*}
and the result follows. The proof is complete now.
\end{proof}

\begin{rem}
The Brauer group version of Lemma \ref{desc-prop2-lem1-lem} does not hold for general $G$.
But if $G$ is linear, we still have $\Br_G(X)=\Phi_{r,\Br}^{-1}(\Br_{G_{\bk}}(X_{\bk})) $,
where $\Phi_{r,\Br}\colon  \Br(X)\to \Br(X_{\bk}) $.
This follows from  \cite[Lem.~6.6]{San81} and \cite[Prop.~3.7]{C1} by the same proof.
\end{rem}

\begin{lem}\label{desc-prop2-lem2}
Consider the following commutative diagram of exact sequences of abelian groups
\begin{equation*}
    \xymatrix{A\ar[r]\ar[d]^e&B\ar[r]^{\varphi}\ar[d]^f&C\ar[r]\ar[d]^g&D\ar[d]^h\\
    A'\ar[r]&B'\ar[r]^{\varrho}&C'\ar[r]&D',
    }
\end{equation*}
where $e$ is surjective and $h$ is injective.
\begin{itemize}
    \item[(1)] For any subgroup $H$ of $C$, we have 
    \begin{equation*}
        \varrho^{-1}(g(H))=f(\varphi^{-1}(H)).
    \end{equation*}
    \item[(2)] For any subgroup $H$ of $B'$, we have
    \begin{equation*}
        g^{-1}(\varrho(H))=\varphi(f^{-1}(H)).
    \end{equation*}
\end{itemize}
\end{lem}

\begin{proof}
    (1) It suffices to show that $\varrho^{-1}(g(H))\subset f(\varphi^{-1}(H))\subset B'$. For every $\alpha\in  \varrho^{-1}(g(H))$, there exists $\beta \in H$ such that $g(\beta)=\varrho(\alpha)\in C'$. Since $h$ is injective, there exists $\gamma\in B$ such that $\varphi(\gamma)=\beta$. Thus, $f(\gamma)-\alpha\in A'$, and since $e$ is surjective, we could lift it to $A$, say $\theta$. Then we have $f(\gamma+\theta)=\alpha$ with $\varphi(\gamma+\theta)=\varphi(\gamma)=\beta\in H$, i.e., $\alpha\in  f(\varphi^{-1}(H))$.

    (2) It suffices to show that $g^{-1}(\varrho(H))\subset\varphi(f^{-1}(H))\subset C$. The proof is dual to that of (1), and we leave the details to the reader.
\end{proof}

\begin{lem}\label{desc-prop2-lem1}
Let $X$ be a quasi-projective smooth geometrically integral $k$-variety and let $g\colon Z\to X$ be a torsor under a quasi-trivial torus $T$.
Assume $[Z]\in H^1(X,T)$ is a torsion element.
Then, for any subgroup $H\subset H^2(X_{\bk},\mu_{\infty})$, we have that
\begin{itemize}
    \item[(a)] the inclusion $\Phi_Z(\bar{g}^*(H)) \subset \Br_T(Z) $;
    \item[(b)] the descent formula $X(\RA_k)^{\Phi_X(H)}=g(Z(\RA_k)^{\Phi_{Z}(\bar{g}^*(H))})$.
\end{itemize}
\end{lem}

\begin{proof}
We denote by $\rho\colon  T\times Z \to Z$ the action and $p_2\colon  T\times Z \to Z$ the projection.
Since $g\circ \rho=g\circ p_2$, we have $(\bar{\rho}^*-\bar{p}_2^*)\circ \bar{g}^*=0$
and hence we have $\bar{g}^*(H)\subset H^2_{T_{\bk}}(Z_{\bk},\mu_{\infty}) $.
Lemma \ref{desc-prop2-lem1-lem} implies the inclusions $\Phi_{r,Z}^{-1} (\bar{g}^*(H) )\subset H^2_T(Z,\mu_{\infty})$ and $\Phi_Z(\bar{g}^*(H)) \subset \Br_T(Z)$.
This proves (a).

Now we prove (b). By  \cite[(5-2)]{C4}, we have the descent formula
$$ X(\RA_k)^{(g_{\Br}^*)^{-1}B}=g(Z(\RA_k)^{B}) $$
 for any subgroup $B\subset \Br_{T}(Z)$.
 Taking $B\coloneqq  \Phi_Z(\bar{g}^*(H)) $,
therefore, to prove (b), we only need to show
\begin{equation*}
  (g^*_{\Br})^{-1}(\Phi_Z(\bar{g}^*(H)) ) = \Phi_X(H).  
\end{equation*}

The Kummer exact sequence (\ref{S5-e-Ku}) induces the following commutative diagram of exact sequences
$$\xymatrix{0\ar[r] & \Pic(X)\otimes \BQ/\BZ  \ar[r]\ar@{->>}[d]^{g^*_{\Pic}\otimes \BQ/\BZ}  &  H^2(X,\mu_{\infty})\ar[rr]^{\Phi_{Ku,X}} \ar[d]^{g^*} && \Br(X) \ar[r]\ar[d]^{g^*_{\Br}} &0\\
0\ar[r] & \Pic(Z)\otimes \BQ/\BZ  \ar[r]  &  H^2(Z,\mu_{\infty})\ar[rr]^{\Phi_{Ku,Z}} & &\Br(Z) \ar[r] &0.
} $$
Since $\Pic(T)\cong H^1(k,\hat{T})$=0, the pullback $g^*_{\Pic}\colon  \Pic(X)\to \Pic(Z)$ is surjective 
by Sansuc's exact sequence (\cite[Prop. 6.10]{San81}).
Therefore, $ g^*_{\Pic}\otimes \BQ/\BZ$ is surjective and further $\Coker(g^*)\cong \Coker(g^*_{\Br})$.
As a consequence, to prove (b), we only need to show $(g^*)^{-1}(\Phi_{r,Z}^{-1}(\bar{g}^*(H)) ) = \Phi_{r,X}^{-1}(H)$; in fact, by Lemma~\ref{desc-prop2-lem2}~(2), we then have 
\begin{equation*}
    \Phi_X(H)={\Phi_{Ku,X}}((g^*)^{-1}(\Phi_{r,Z}^{-1}(\bar{g}^*(H)) ) )=(g_{\Br}^*)^{-1}(\Phi_{Ku,Z}(\Phi_{r,Z}^{-1}(\bar{g}^*(H))))=(g^*_{\Br})^{-1}(\Phi_Z(\bar{g}^*(H)) ). 
\end{equation*}
By the commutative diagram (\ref{S5-e-Ku-1}) (with replacing $f\colon Y\to X$ by $g\colon  Z\to X$), 
it suffices to prove that $ (\bar{g}^*)^{-1}(\bar{g}^*(H)) =H$.
For this, we only need to show $\bar{g}^* $ is injective.

Since $[Z]\in H^1(X,T)$ is torsion, there exists an integer $N$ such that $N\cdot [Z]=0\in H^1(X,T)$.
Hence $\psi_{N,*}([Z])=0=[T\times X]\in H^1(X,T)$, where $\psi_N\colon T\to T\colon  t\mapsto t^N$.
Applying the $\mu_{\infty}$-version of Sansuc's exact sequence (\cite[Thm. 4.1]{C5}) 
to the $T$-torsor $g$ and  to the trivial $T$-torsor $p_2\colon  T\times X\to X$ respectively, it implies a commutative diagram with exact rows
$$\xymatrix{H^1(T_{\bk},\mu_{\infty}) \ar[r]^{\partial_1} \ar[d]^{\psi_N^*|_{H^1}}  &  H^2(X_{\bk},\mu_{\infty})\ar[r]^-{\bar{p}_2^*} \ar[d]^=
& H^2((T\times X)_{\bk},\mu_{\infty})\ar[d] \\
H^1(T_{\bk},\mu_{\infty}) \ar[r]^{\partial} & H^2(X_{\bk},\mu_{\infty})\ar[r]^{\bar{g}^*} 
&H^2(Z_{\bk},\mu_{\infty}) .
}$$
Since, over algebraically closed fields, the projection $T\times_{\bar{k}} X\to X_{\bar{k}}$ always has a section, $\bar{p}_2^* $ is injective and hence $\partial_1=0 $.
The Kummer sequence implies the short exact sequence
\begin{equation*}
    0\to \bar{k}[T]^{\times}\otimes_{\BZ} \BQ/\BZ\to H^1(T_{\bar{k}},\mu_\infty)\to \mathrm{Pic}(T_{\bar{k}})_{\rm tors}\to 0,
\end{equation*}
hence 
\begin{equation}\label{equ-H1(T)}
    H^1(T_{\bk},\mu_{\infty})\cong\bar{k}[T]^{\times}\otimes_{\BZ} \BQ/\BZ\cong (\bar{k}^*\oplus \hat{T}_{\bk})\otimes_{\BZ} \BQ/\BZ\cong \hat{T}_{\bk}\otimes_{\BZ} \BQ/\BZ
\end{equation}
by Rosenlicht's lemma (cf. \cite[Lem. 6.5]{San81}).
Therefore, $\psi_N^*|_{H^1}$ is surjective and hence $\partial=0$.
This implies that $ \bar{g}^*$ is injective and therefore (b) holds. The proof is complete.
\end{proof}

\begin{proof}[Proof of Proposition \ref{desc-prop2}]
We follow the proof of \cite[Prop. 5.1]{C4}.

By Colliot-Th\'el\`ene and Sansuc's flasque resolution (\cite[Prop. 1.3]{CTS87b}), we have an exact sequence
\begin{equation}\label{eq-flasquedescent}
0\to A\xrightarrow{\phi_1} T\xrightarrow{\phi_2} T_1\to 0, 
\end{equation}
where $T$ is a quasi-trivial torus and $T_1$ is a coflasque torus, i.e $H^1(k',\hat{T_1})=0$ for every field extension $k'/k$.
By \cite[Lem. 5.4]{C1}, $H^3(k,\hat{T}_1)=0$.
Moreover, this exact sequence induces the connecting map $\partial\colon  T_1(k)\to H^1(k,A)$.
For any $t\in T_1(k)$, we have the twist $f_{\partial(t)}\colon  Y_{\partial(t)}\to X $.

Let $Z\coloneqq Y\times^AT$ the contract product and $g\colon Z\to X$ the corresponding morphism.
Then we have $[Z]=\phi_{1,*}([Y])\in H^1(X,T)$ and hence $\phi_{2,*}([Z])=0$.
Therefore, the isomorphism $Z\times^{T}T_1\cong X\times T_1$ induces a morphism $h\colon  Z\to T_1$ which is compatible with $T$-action.
Therefore, we have the commutative diagram
\begin{equation}\label{desc-prop2-eq-diag}
\xymatrix{Y\ar[rd]_-f  \ar[r] &  Z:=Y\times^AT \ar[r]^-h \ar[d]^g &  T_1 \\
&X&.
}
\end{equation}
By a classical argument (see \cite[Thm. 5.9]{C1} for the proof), for any $t\in T_1(k)$, the fibre $$Z_t\coloneqq h^{-1}(t)\cong Y_{\partial(t)}$$
and the composition $g_t\colon   Z_t\hookrightarrow Z\xrightarrow{g}X$ is exactly 
the twist $f_{\partial(t)}\colon  Y_{\partial(t)}\to X$.

We denote by $i_t\colon  Z_t\to Z$ the inclusion and the pullback $i_t^*, \bar{i}_t^*, i^*_{t,\Br} $ as in (\ref{S5-e-Ku-1}).

\begin{claim}
We have $i^*_{t,\Br}(\Phi_Z(\bar{g}^*(H)))=\Phi_{Z_t}(\bar{g}_t^* (H)) $ 
for any $t\in T_1(k)$.
\end{claim}

\begin{proof}[Proof of the Claim]
Applying the commutative diagram (\ref{S5-e-Ku-1}) to $i_t\colon Z_t\to Z $, it suffices to show
\begin{equation}\label{desc-prop2-e4}
i^*_t(\Phi_{r,Z}(\bar{g}^*(H)))=\Phi_{r,Z_t}(\bar{g}_t^* (H)) .
\end{equation}

Since $H^3(k,\mu_{\infty})=0$, 
the Hochschild--Serre spectral sequence (\cite[Thm. III.2.20]{Mi80}) for $Z$ and for $Z_t$ induce two commutative diagrams of exact sequences
\begin{equation}\label{desc-prop2-e3}
 \xymatrix{0\ar[r] & \Ker(\Phi_{r,Z})\ar[r]\ar[d]^{\Phi_1}  &  H^2(Z,\mu_{\infty}) \ar[r]^-{\Phi_{r,Z}}\ar[d]^{i_t^*} & 
H^2(Z_{\bk},\mu_{\infty})^{\Gamma_k} \ar[r]\ar[d]^{\bar{i}_t^* } & H^2(k, H^1(Z_{\bk},\mu_{\infty}))\ar[d]^{\Phi_2} \\
0\ar[r] & \Ker(\Phi_{r,Z_t})\ar[r]  & H^2(Z_t,\mu_{\infty}) \ar[r]^-{\Phi_{r,Z_t}} & H^2(Z_{t,\bk},\mu_{\infty})^{\Gamma_k} \ar[r] &
H^2(k, H^1(Z_{t,\bk},\mu_{\infty}))
}   
\end{equation}
and 
$$\xymatrix{H^2(k,\mu_{\infty}) \ar[r]\ar[d]^= &   \Ker(\Phi_{r,Z}) \ar[r] \ar[d]^{\Phi_1} & H^1(k, H^1(Z_{\bk},\mu_{\infty})) 
\ar[d]^{\Phi_3} \ar[r] & 0\\
H^2(k,\mu_{\infty}) \ar[r]   &\Ker(\Phi_{r,Z_t})   \ar[r]  & H^1(k, H^1(Z_{t,\bk},\mu_{\infty}) ) \ar[r] &  0 .
}  $$ 

By Proposition \ref{prop-Katz-Lang} (2), we have an exact sequence
\begin{equation}\label{desc-prop2-e1}
 \pi_1(Z_{t,\bk})^{ab}\xrightarrow{\psi} \pi_1(Z_{\bk})^{ab}\to \pi_1(T_{1,\bk})^{ab}\to 0.
\end{equation}
Since $\pi_1(X_{\bk}) $ is abelian, $Z_{t,\bk}$, as a finite \'etale cover of $X_{\bk}$, 
satisfies that  $\pi_1(Z_{t,\bk})$ is abelian and the morphism $\pi_1(Z_{t,\bk})\to \pi_1(X_{\bk})$ is injective.
Hence $\psi$ is injective since the morphism $\pi_1(Z_{t,\bk})\to \pi_1(X_{\bk})$ factors through $\psi$. 
 Applying $\Hom(-,\mu_{\infty})$ to (\ref{desc-prop2-e1}), we get an exact sequence
\begin{equation}\label{desc-prop2-e2}
    0\to H^1(T_{1,\bk},\mu_{\infty})\to H^1(Z_{\bk},\mu_{\infty})\to H^1(Z_{t,\bk},\mu_{\infty})\to 0.
\end{equation}
Since $H^1(T_{1,\bk},\mu_{\infty})=\hat{T}_1\otimes_{\BZ}\BQ/\BZ$ (see (\ref{equ-H1(T)})) and $T_1$ is coflasque, one has $$H^2(k,H^1(T_{1,\bk},\mu_{\infty}))\cong H^3(k,\hat{T}_1)=0$$ by \cite[Lem. 5.4]{C1} again.
Applying Galois cohomology to (\ref{desc-prop2-e2}), the homomorphism $\Phi_2$ is injective and $\Phi_3$ is surjective.
By the five lemma, the homomorphism $\Phi_1$ is also surjective.

Since $H\subset H^2(X_{\bk},\mu_{\infty})^{\Gamma_k}$, 
we have $\bar{g}^*(H)\subset H^2(Z_{\bk},\mu_{\infty})^{\Gamma_k} $ and 
$\bar{g}_t^* (H)=\bar{i}_t^* (\bar{g}^*(H))\in  H^2(Z_{t,\bk},\mu_{\infty})^{\Gamma_k} $.
We apply Lemma \ref{desc-prop2-lem2}~(1) to diagram (\ref{desc-prop2-e3}), and we obtain (\ref{desc-prop2-e4}), which implies our claim.
\end{proof}

Let us return to our proof of Proposition \ref{desc-prop2}. Since $A$ is finite, there exists an integer $N$ such that $N\cdot [Y]=0\in H^1(X,A)$. 
Hence $N\cdot [Z]=0$ and the hypothesis of Lemma \ref{desc-prop2-lem1} is satisfied.
By Lemma \ref{desc-prop2-lem1}, for any $x\in X(\RA_k)^{\Phi_X(H)} $, 
there exists $z\in Z(\RA_k)^{\Phi_{Z}(\bar{g}^*(H))}$ such that $g(z)=x$.
Since $\Br_{2/3}(T_1)=\Br_1(T_1)$ (see \cite[Cor. 3.7]{C5}), we have $h(z)\in T_1(\RA_k)^{\Br_1(T_1)} $.
By \cite[Thm. 5.1]{CDX} (see also \cite[(5-3)]{C4}), $$T_1(\RA_k)^{\Br_1(T_1)}=\phi_2(T(\RA_k)^{\Br_1(T)})\cdot T_1(k),$$
 and there exists a $t_0\in T(\RA_k)^{\Br_1(T)}$ such that $t\coloneqq h(t_0\cdot z)\in T_1(k)$.
Since $\Phi_{Z}(\bar{g}^*(H))\subset \Br_{T}(Z)$ (see Lemma~\ref{desc-prop2-lem1}~(a)), by \cite[Prop. 3.9(1)]{C1} and the above claim,
we find that
$$t_0\cdot z\in Z(\RA_k)^{\Phi_{Z}(\bar{g}^*(H))}\cap Z_t(\RA_k)=Z_t(\RA_k)^{\Phi_{Z_t}(\bar{g}_t^* (H))}.$$
Therefore, $x=g(z)=g(t_0\cdot z)\in f_{\partial(t)}(Y_{\partial(t)}(\RA_k)^{\Phi_{Y_{\partial(t)}}(\bar{f}_{\partial(t)}^* (H) )  })$ as $Z_t\cong Y_{\partial(t)}$,
which proves our result. The proof is complete.
\end{proof}

\subsection{Descent formula with respect to Hurewicz map}

For any positive integer $m$, let $A_{m}$ be a finite commutative $k$-group such that its \emph{Cartier dual} $\hat{A}_{m}\coloneqq \Hom_{\bk}(A_m,\BG_m)$ is isomorphic to the group $ H^1(X_{\bk},\mu_m)$ as Galois modules.
As in \cite[Def. 2.1]{C4}, 
a \emph{universal torsor of $m$-torsion} is a certain $A_{m}$-torsor $\CT_{m}\to X$.
The corresponding finite \'etale abelian cover $\CT_{m,\bk}\to X_{\bk}$ induces a homomorphism $\tau_{m}\colon  \pi_1(X_{\bk})\to A_{m}(\bk)$. 
By \cite[(2-6)]{C4}, the key property of universal torsor of $m$-torsion is that
there exists an isomorphism  $\pi_1(X_{\bk})^{ab}/m\cong A_{m}(\bk)$ such that
$\tau_{m}$ is exactly the composition
$$\pi_1(X_{\bk})\to \pi_1(X_{\bk})^{ab}\to \pi_1(X_{\bk})^{ab}/m\cong A_{m}(\bk).$$
Therefore, $\CT_{m}$ is geometrically integral and the morphism $\CT_{m}\to X$ induces an exact sequence
\begin{equation}\label{desc-result-1-e-pi}
1\to \pi_1(\CT_{m,\bk} )\to \pi_1(X_{\bk})\to \pi_1(X_{\bk})^{ab}/m \to 1.
\end{equation}

\begin{lem}\label{desc-result-1-lem}
Let $X$ be a quasi-projective smooth geometrically integral $k$-variety.
If $X(\RA_k)^{\Br_{2/3}(X)}\neq \emptyset$, then, for any integer $m>0$, there exists a universal torsor of $m$-torsion.
\end{lem}

\begin{proof}
Let $p\colon X\to \Spec\ k$ be the structure morphism of $X$. In the derived category of Galois modules of $k$, consider the object $Rp_*\BG_m$  and the following mapping cone defined in \cite[\S 8.1]{HS13}
$$KD'(X)\coloneqq \mathrm{Cone} (\BG_m\to \tau_{\leq 1}Rp_*\BG_m)[1],$$
where $\tau_{\leq 1}Rp_*\BG_m$ is the good truncated of $Rp_*\BG_m$ at degree $1$ and the functor $[1]$ means the shift by $1$. In other words, we have the following exact triangle
\begin{equation}\label{equ-exact triangle}
    \BG_m\to \tau_{\leq 1}Rp_*\BG_m\to KD'(X)[-1]\xrightarrow{+1}.
\end{equation}
 Applying the functor $\Hom_k(\hat{A}_m,\cdot)$ to the exact triangle (\ref{equ-exact triangle}), one could get the following exact sequence (see \cite[Prop. 8.1]{HS13})
\begin{equation}\label{desc-type-eq}
H^1(k,A_m)\to H^1(X,A_m) \xrightarrow{\chi} \Hom_k(\hat{A}_m, KD'(X)) \xrightarrow{\partial} H^2(k,A_m)\to H^2(X,A_m).
\end{equation}
For every $A_m$-torsor $Y\in H^1(X,A_m)$ over $X$, the image $\chi([Y])\in \Hom_k(\hat{A}_m, KD'(X))$ is called the \emph{extended type} of $Y$ (\cite[Def. 8.2]{HS13}).
By \cite[Prop. 2.2]{C5}, if $X(\RA_k)^{\Br_{2/3}(X)}\neq \emptyset$, then $\chi$ is surjective, 
i.e., every $\lambda\in \Hom_k(\hat{A}_m, KD'(X))$ is the extended type of a certain $A_m$-torsor. 
In particular, by \cite[Def.~2.1]{C4}, a universal torsor of $m$-torsion is exactly an $A_m$-torsor corresponding to a certain extended type (the type $\psi_{\leq 0}$ in \cite{C4}), hence it exists. The proof is complete.
\end{proof}

\begin{prop}\label{desc-result-1}
Let $X$ be a quasi-projective smooth geometrically integral $k$-variety.
Assume $\pi_1(X_{\bk})$ is abelian. 
Then, for any subgroup $H\subset H^2(X_{\bk},\mu_{\infty})^{\Gamma_k} $,  
\begin{enumerate}
\item we have $ X(\RA_k)^{\Phi_X(H)}=X(\RA_k)^{\Phi_X(H+\Im(\bar{\epsilon}_X))}$;
\item for any finite commutative $k$-group $A$ and any $A$-torsor $f\colon Y\to X$ with $Y$ geometrically integral,
we have 
$$ X(\RA_k)^{\Phi_X(H)} =
\bigcup_{\sigma\in H^1(k,A)}f_{\sigma} Y_{\sigma}(\RA_k)^{\Phi_{Y_{\sigma}}(\bar{f}_{\sigma}^* (H)+ \Im(\bar{\epsilon}_{Y_{\sigma}}))}.$$
\end{enumerate}
\end{prop}

\begin{proof}
It is clear that $ X(\RA_k)^{\Phi_X(H)}\supset X(\RA_k)^{\Phi_X(H+\Im(\bar{\epsilon}_X))}$.
To prove the converse, we may assume $$ X(\RA_k)^{\Phi_X(H)}\neq\emptyset.$$
Therefore, $\Br_{2/3}(X)=\Phi_X(0)\subset \Phi(H)$ and $ X(\RA_k)^{\Br_{2/3}(X)}\supset  X(\RA_k)^{\Phi_X(H)}\neq\emptyset $.
 By Lemma \ref{desc-result-1-lem}, for any integer $m>0$, there exists a universal torsor of $m$-torsion.

For any finite subgroup $B\subset \Im(\bar{\epsilon}_X)\subset H^2(X_{\bk},\mu_{\infty})$, by Proposition \ref{Prop-Hur-1},
 there exists a finite \'etale Galois cover $Z/X_{\bk}$ such that $B|_Z=0$.
Since $\pi_1(X_{\bk})$ is abelian, the corresponding surjective homomorphism $\pi_1(X_{\bk})\to \Gal(Z/X_{\bk})$ factors through $\pi_1(X_{\bk})\to \pi_1(X_{\bk})/m$ for some $m$.
Let $f_m\colon  \CT_{m}\to X$ be a universal torsor of $m$-torsion under the $k$-group $A_{m}$. 
By (\ref{desc-result-1-e-pi}), $\CT_{m,\bk}$ is a Galois cover of $Z$, and hence $B|_{\CT_{m,\bk}}=0$.
Therefore, for any $\sigma\in H^1(k,A_{m})$, we have $\bar{f}_{m,\sigma}^*(H)=\bar{f}_{m,\sigma}^*(H+B)$.
 By Proposition \ref{desc-prop2}, 
\begin{equation*}
\begin{split}
X(\RA_k)^{\Phi_X(H)} &=  \bigcup_{\sigma\in H^1(k,A_{m})}f_{m,\sigma}(\CT_{m,\sigma}(\RA_k)^{\Phi_{\CT_{m,\sigma}}(\bar{f}_{m,\sigma}^* (H)) })\\
&= \bigcup_{\sigma\in H^1(k,A_{m})}f_{m,\sigma}(\CT_{m,\sigma}(\RA_k)^{\Phi_{\CT_{m,\sigma}}(\bar{f}_{m,\sigma}^* (H+B)) }) = X(\RA_k)^{\Phi_X(H+B) }.    
\end{split}
\end{equation*}
 
This implies (1).

The statement (2) follows from (1) and Proposition \ref{desc-prop2}. 
\end{proof}

\begin{proof}[Proof of Proposition \ref{desc-result-1-prop}]
Since $\Phi_X(0)=\Br_{2/3}(X)\subset \Br_1(X)$, we have $X(\RA_k)^{\Br_1(X)}\subset X(\RA_k)^{\Phi_X(0)}$.
On the other hand, $\Im(\epsilon'_X)\subset  \Phi_X(\Im(\bar{\epsilon}_X) )$ 
and hence $X(\RA_k)^{\Im(\epsilon'_X)}\supset X(\RA_k)^{\Phi_X(\Im(\bar{\epsilon}_X) )}$.
Similarly, the same conclusion holds for $Y_{\sigma}$.
Now we apply Proposition \ref{desc-result-1} to the case $H=0$ and the result follows.
\end{proof}

Now, we prove our Theorem \ref{thm-1}.
We need the following calculation about the Brauer group.

\begin{lem}\label{lem-desc-sym}
Let $X$ be a quasi-projective smooth geometrically integral $k$-variety. 
We consider the $\tSym^nX $ as in (\ref{notation-sym-eq}) 
and the homomorphisms $\tr^2_X$, $\tr^2_{X_{\bk}}$ and $\tr_{X,\Br}$ as in (\ref{eq-def-transfer}).
Let $H\coloneqq  \tr_{X_{\bk}}^2 (H^2(X_{\bk},\mu_{\infty})^{\Gamma_k})$.
Then we have
\begin{enumerate}
\item $\Phi_{\tSym^nX }(H)= \Im(\tr_{X,\Br}) + \Im(\Br(k))$;
\item $\Phi_{\tSym^n(X)}(H+\Im(\epsilon_{\tSym^nX_{\bar{k}}}))= \Br(\tSym^nX)$.
\end{enumerate}
\end{lem}

\begin{proof}
For (1), it is clear that we already have $\Im(\tr_{X,\Br}) \subset \Phi_{\tSym^nX }(H)$
and $\Im(\Br(k))\subset \Phi_{\tSym^nX }(H)$. Hence, it suffices to show that 
\begin{equation*}
 \Phi_{\tSym^nX }(H)\subset \Im(\tr_{X,\Br}) + \Im(\Br(k)).  
\end{equation*}
By the following commutative diagram
\begin{equation*}
    \xymatrix{H^2(\tSym^nX,\mu_\infty)\ar@{->>}[rrr]^{\Phi_{Ku,\tSym^nX}}&&& \Br(\tSym^nX)\\
    H^2(X,\mu_\infty)\ar@{->>}[rrr]^{\Phi_{Ku,X}}\ar[u]^{\tr^2_X}&&&\Br(X),\ar[u]^{\tr_{X,\Br}}
    }
\end{equation*}
it suffices to show that
\begin{equation*}
   \Phi_{r,\tSym^nX }^{-1}(H)\subset \Im(\tr^2_{X}) + \Im(H^2(k,\mu_\infty)). 
\end{equation*}

Picking an element $\alpha\in \Phi_{r,\tSym^nX }^{-1}(H)\subset H^2(\tSym^nX,\mu_\infty)$, 
one could find $\beta\in H^2(X_{\bar{k}},\mu_\infty)^{\Gamma_k}$ such that 
\begin{equation*}
    \tr^2_{X_{\bar{k}}}(\beta)=\Phi_{r,\tSym^nX }(\alpha)\in H^2(\tSym^nX_{\bar{k}},\mu_{\infty})^{\Gamma_k}.
\end{equation*}
Consider the following commutative diagram of exact sequences induced by Hochschild--Serre spectral sequences
\begin{equation*}
\xymatrix{ H^2_1(X,\mu_\infty) \ar@{^(->}[r] \ar[d]^{\tr_X^2|_{H^2_1}} &  H^2(X,\mu_\infty)\ar[r]\ar[d]^{\tr_X^2} &
H^2(X_{\bk},\mu_\infty)^{\Gamma_k}\ar[r]\ar[d]^{\tr_{X_{\bk}}^2}
& H^2(k,H^1(X_{\bk},\mu_\infty) )  \ar[d]_{\sim}^{H^2(k, \tr_{X_{\bk}}^1)}   \\
H^2_1(\tSym^nX,\mu_\infty)\ar@{^(->}[r] & H^2(\tSym^nX,\mu_\infty)\ar[r] &  H^2(\tSym^nX_{\bk})^{\Gamma_k}\ar[r]  
&  H^2(k,H^1(\tSym^nX_{\bk},\mu_\infty)  ) ,
}
\end{equation*}
where the last vertical arrow is an isomorphism by Theorem \ref{thm-DT-1} (iii). By a diagram chasing, one could lift $\beta$ to $\beta'\in H^2(X,\mu_\infty)$ and $\tr^2_X(\beta')-\alpha\in H^2_1(\tSym^nX,\mu_\infty) $. Moreover, we have another commutative diagram of exact sequences induced by Hochschild--Serre spectral sequences
\begin{equation*}
    \xymatrix{H^2(k,\mu_\infty)\ar[r] \ar[d]^{H^2(k,\tr_{X_{\bk}}^{n,0})}  &
H^2_1(X,\mu_\infty) \ar[r]\ar[d]^{\tr_X^2|_{H^2_1}}
& H^1(k,H^1(X_{\bk},\mu_\infty) ) \ar[r]\ar[d]_{\sim}^{H^1(k,\tr_{X_{\bk}}^1)} &0 \\
H^2(k,\mu_\infty)\ar[r]_-{\pi^* }   &   H^2_1(\tSym^nX,\mu_\infty) \ar[r] & H^1(k, H^1(\tSym^nX_{\bk},\mu_\infty) ) \ar[r]  &  0.
}
\end{equation*}
Thus, we have a surjection $H^2_1(X,\mu_\infty) \oplus H^2(k,\mu_\infty)\xrightarrow{(\tr_X^2|_{H^2_1},\pi^*)}H^2_1(\tSym^nX,\mu_\infty)$ and further there exist $\gamma\in H^2_1(X,\mu_\infty)\subset H^2(X,\mu_\infty)$ and $\omega\in H^2(k,\mu_\infty)$ such that
\begin{equation*}
    \tr^2_X(\beta')-\alpha=\tr_X^2(\gamma)+\pi_X^*(\omega)\in H^2_1(\tSym^nX,\mu_\infty) \subset H^2(\tSym^nX,\mu_\infty),
\end{equation*}
which implies that 
\begin{equation*}
    \alpha=\tr^2_X(\beta'-\gamma)-\pi_X^*(\omega)\in H^2(\tSym^nX,\mu_\infty).
\end{equation*}
Therefore, $\alpha\in \Im(\tr^2_{X}) + \Im(H^2(k,\mu_\infty))$.

For (2), it is clear that it suffices to show that $\Phi_{\tSym^n(X)}\left(H+\Im(\epsilon_{\tSym^nX_{\bar{k}}})\right)\supset  \Br(\tSym^nX)$. 
For every $\alpha\in \Br(\tSym^nX)\subset \Br(\tilde{V}(X))$, we could lift it to $\alpha'\in H^2(\tilde{V}(X),\mu_\infty)$. Hence, by Theorem \ref{thm-DT-2}, we have
\begin{equation*}
    \Phi_{r,\tilde{V}(X)}(\alpha')=\epsilon_{\tilde{V}(X)}(\beta)+\tr^2_X(\gamma)+i_{D_{\Bar{X}}}(\omega)\in 
    H^2(\tilde{V}(X)_{\bk},\mu_\infty)^{\Gamma_k}
\end{equation*}
for some $\beta\in H_{cont}^2(\pi_1(\tilde{V}(X))_{\bk},\mu_\infty)\cong H_{cont}^2(\pi_1(\tSym^n X)_{\bk},\mu_\infty)$, 
$\gamma\in H^2(X_{\bk},\mu_\infty)$ and $\omega\in \BQ/\BZ$. 
By Corollary \ref{cor-iso for Q/Z}, $\gamma\in H^2(X_{\bk},\mu_\infty)^{\Gamma_k}$.
Note that $i_{D_{\bar{X}}}(\omega)=\Phi_{r,\tilde{V}(X)}(i_{D_X}(\omega))\in H^2(\tilde{V}(X)_{\bk},\mu_\infty)$ (Proposition \ref{notation-sym-prop}). 
Thus, by diagram chasing, one could get that
\begin{equation*}
  \alpha\in \Phi_{\tSym^n(X)}(\epsilon_{\tilde{V}(X)}(\beta)+\tr^2_X(\gamma))\subset   \Phi_{\tSym^n(X)}(\Im(\tr^2_{X_{\bk}})+\Im(\epsilon_{\tSym^nX_{\bar{k}}}))\subset\Br(\tSym^nX).
\end{equation*}
The proof is complete.
\end{proof}

\begin{rem}
From the above proof, the results for $\tSym^n(X)$ also hold for $\tilde{V}(X)  $, i.e., we have
\begin{enumerate}
\item $\Phi_{\tilde{V}(X)}(H)= \Im(\tr_{X,\Br}) + \Im(\Br(k))$;
\item $\Phi_{\tilde{V}(X)}(H+\Im(\epsilon_{\tilde{V}(X)}))= \Br(\tilde{V}(X))$.
\end{enumerate}
\end{rem}

\begin{proof}[Proof of Theorem \ref{thm-1}]
We may assume $n\geq 2$. 
By Theorem \ref{thm-2}, $\pi_1(\tSym^nX_{\bk}) $ is abelian.
We apply Proposition \ref{desc-result-1} (1) to the case $H=\tr^2_{X_{\bk}}(H^2(X_{\bk},\mu_{\infty})^{\Gamma_k})$, and the result follows from Lemma \ref{lem-desc-sym}.
\end{proof}

\section{Proof of Corollary \ref{thm-1-cor}}\label{section6}

Let $k$ be a number field and $X$ a projective smooth geometrically integral variety over $k$.
Recall our notations $\Sym^nX$, $\tSym^nX$ and $[-]$ in Corollary \ref{thm-1-cor}
and $U(X)$, $U^n(X)$ in (\ref{notation-sym-eq}).

\begin{lem}\label{thm-1-cor-lem}
Let $X$ be a projective smooth geometrically integral variety over $k$.
Then, for any $a\in \Br(X)$ and any $(x_v)\in \tSym^nX(\RA_k)$, the Brauer--Manin pairing satisfies
 $$\langle [(x_v)],a\rangle_{\BM} = \langle (x_v), \tr_{X,\Br}(a)\rangle_{\BM}  . $$
\end{lem}

\begin{proof}
Since the Brauer--Manin pairing over $k$ is the sum of the Brauer--Manin pairing over all local field $k_v$,
it suffices to establish the statement over local field $k_v$.
More precisely, it suffices to show that for any $a\in \Br(X)$ and any $x_v\in \tSym^nX(k_v) $, one has
$$ \langle [x_v]_v,a\rangle_{\BM} = \langle x_v, \tr_{X,\Br}(a)\rangle_{\BM},$$
where $[-]_v \colon \tSym^nX(k_v)\to \CH_0(X_{k_v})$.

We consider two maps $$\phi_1,\phi_2\colon \tSym^nX(k_v)\to \BQ/\BZ$$
given by $\phi_1(x_v)\coloneqq  \langle x_v, \tr_{X,\Br}(a)\rangle_{\BM}$ and $\phi_2(x_v)\coloneqq  \langle [x_v]_v,a\rangle_{\BM} $.
By \cite[Prop. 10.5.2]{CTS}, $\phi_1 $ is locally constant.
Since the Brauer group is torsion, there exists an integer $m>0$ such that $m\cdot a=0$.
Therefore, the homomorphism $\CH_0(X_{k_v})\to \BQ/\BZ\colon z\mapsto \langle z,a\rangle_{\BM} $ factors through $ \CH_0(X_{k_v})/m$.
Since the composition $ \Sym^n X(k_v)\to \CH_0(X_{k_v})\to  \CH_0(X_{k_v})/m$ is continuous (\cite[Lem. 1.8]{W12}), 
$\phi_2$ is continuous and hence locally constant.
On the other hand, since $U(X)(k_v)\subset \tSym^nX(k_v)$ is open and dense, for any $x_v\in \tSym^nX(k_v)$, 
there exists $u_v\in U(X)(k_v)$ such that $\phi_i(x_v)=\phi_i(u_v)$ for $i=1,2$.
Therefore, we reduce to the case $x_v\in U(X)(k_v) $.

Let $W\subset \tSym^nX\times X $ be the universal family in (\ref{equ-def-W}) and let $p\colon W\to \tSym^nX$ be the projection.
For any $x_v\in U(X)(k_v)$, the fibre $W_{x_v}$ of $p$ is a reduced closed subscheme of $X_{k_v}$ of dimension $0$ and degree $n$,
and we have $[x_v]_v=[W_{x_v}]\in \CH_0(X_{k_v})$. 
Write $W_{x_v}=\{P_1,\cdots ,P_l\} $.
By the definition of Brauer--Manin pairing (\cite[\S 13.3.1 and \S 6.1]{CTS}), 
we have 
$$\langle x_v, \tr_{X,\Br}(a)\rangle_{\BM} \coloneqq \mathrm{inv}_v(x_v^*(\tr_{X,\Br}(a)) )=\mathrm{inv}_v(x_v^*(W^*(a)) ) $$
and 
$$ \langle [x_v]_v,a\rangle_{\BM} \coloneqq  \mathrm{inv}_v\left(\sum_i \mathrm{Cores}_{k_v(P_i)/k_v} P_i^*(a)\right)=
\mathrm{inv}_v (W_{x_v}^*(a)),$$
where 
\begin{itemize}
    \item $x_v^*:\Br(\tSym^nX)\to \Br(k_v)$ and $P_i^*\colon  \Br(X)\to \Br(k_v(P_i))$ are the natural pullbacks;
    \item $W_{x_v}^*: \Br(X)\to \Br(k_v)$ is the transfer map induced by the correspondence $W_{x_v}\in Cor(k_v,X)$;
    \item $\mathrm{inv}_v\colon \Br(k_v)\hookrightarrow \BQ/\BZ$ is the local invariant map in local class field theory.
\end{itemize}
Since $W\circ x_v=W_{x_v} \in Cor(k_v,X)$, we have $x_v^*(W^*(a)) = W_{x_v}^*(a) $ and the proof is complete.
\end{proof}

\begin{proof}[Proof of Corollary \ref{thm-1-cor}]
By Theorem \ref{thm-1}, the adelic point $(x_v)\in \tSym^nX(\RA_k)^{\Im(\tr_{X,\Br})}$ if and only if $(x_v)\in \tSym^nX(\RA_k)^{\Br}$.
Then the result follows from Lemma \ref{thm-1-cor-lem}.
\end{proof}

\bigskip

\textbf{Acknowledgments.} The authors would like to thank Xing Gu, Tam\'as Szamuely, Zhiyu Tian and Dingxin Zhang for helpful discussions.
The first author is partially supported by National Key R\&D Program of China Grant No. 2025YFA1017300 
(subprogram No. 2025YFA1017303).
The second author is supported by the European Union’s Horizon 2020
research and innovation programme under grant agreement No. 948066 (ERC-StG
RationAlgic).

\textbf{AI Declaration.}
The ideas and demonstrations in this article are the product of pure human thought.

\bibliographystyle{alpha}
\bibliography{MyBib}

\end{document}